\documentclass[12pt,twoside]{article}

\usepackage{amsgen,amsmath,amstext,amsbsy,amsopn,amsfonts,amssymb}

\usepackage{upref}

\usepackage{jltmac2e}

\renewcommand{\version}{{\rm Version of~}\date} 
\renewcommand{\version}{{\rm Submitted to the 
 \emph{Journal of Lie Theory} (July 2000).}} 
\renewcommand{\firstheadline}{\hfil \version} 
\renewcommand{\vol}{??}                             
\renewcommand{\annum}{??}                           
\renewcommand{\title}{\centerline{%
 Cartan-decomposition subgroups of $\operatorname{\rm SU(2,n)}$}} 
\renewcommand{\author}{Alessandra Iozzi and Dave Witte}
\renewcommand{\lastname}{Iozzi and Witte}

\renewcommand{\entries}{
 \nextref{Benoist}
 {Benoist,~Y.}
 {\emph{Actions propres sur les espaces homog\`enes r\'eductifs}}
 {Ann. Math. {\bf 144} (1996), 315--347}

\nextref{BorelCK}
 {Borel,~A.}
 {\emph{Compact Clifford-Klein forms of symmetric spaces}}
 {Topology {\bf 2} (1963), 111--122}

\nextref{HelgasonBook}
 {Helgason,~S.}
 {``Differential Geometry, Lie Groups, and Symmetric Spaces"}
 {Academic Press, New York, 1978}

\nextref{IozziWitte-tess}
 {Iozzi,~A., and D.~Witte}
 {Tessellations of homogeneous spaces of classical groups of real rank two}
 {(in preparation)}

\nextref{Kobayashi-isotropy}
 {Kobayashi,~T.}
 {\emph{On discontinuous groups acting on homogeneous spaces with non-compact
isotropy groups}}
 {J. Geom. Physics {\bf 12} (1993), 133--144}

\nextref{Kobayashi-properaction}
 {\sameauthor}
 {\emph{Proper action on a homogeneous space of reductive type}}
 {Math. Ann. {\bf 285} (1989), 249--263}

\nextref{Kobayashi-survey}
 {\sameauthor}
 {\emph{Discontinuous groups and Clifford-Klein forms of
pseudo-Riemannian homogeneous manifolds}}
 {in: B.~\O rsted and H.~Schlichtkrull, eds.,
 ``Algebraic and Analytic Methods in Representation Theory",
 Academic Press, New York, 1997, pp.~99--165}

\nextref{Kobayashi-criterion}
 {\sameauthor}
 {\emph{Criterion of proper actions on homogeneous spaces of reductive groups}}
 {J. Lie Th. {\bf 6} (1996), 147--163}

\nextref{Kostant}
 {Kostant,~B.}
 {\emph{On convexity, the Weyl group, and the Iwasawa decomposition}}
 {Ann. Sc. ENS. {\bf 6} (1973), 413--455}

\nextref{KulkarniCK}
 {Kulkarni,~R.}
 {\emph{Proper actions and pseudo-Riemannian space forms}}
 {Adv. Math. {\bf 40} (1981) 10--51}

\nextref{Labourie-survey}
 {Labourie,~F.}
 {\emph{Quelques r\'esultats r\'ecents sur les espaces localement
homog\`enes compacts}}
 {in: P.~de Bartolomeis, F.~Tricerri and E.~Vesentini, eds.,
 ``Manifolds and Geometry",
 Symposia Mathematica, v.~XXXVI,
 Cambridge U.~Press, 1996}

\nextref{MargulisCK}
 {Margulis,~G.~A.}
 {\emph{Existence of compact quotients of homogeneous
spaces, measurably proper actions, and decay of matrix coefficients}}
 {Bull. Soc. Math. France {\bf 125} (1997) 447--456}

\nextref{MilnorStasheff}
 {Milnor,~J.~W., and J.~D.~Stasheff}
 {``Characteristic Classes"}
 {Princeton U.\ Press, Princeton, 1974}

\nextref{OhWitte-CDS}
 {Oh,~H., and D.~Witte}
 {\emph{Cartan-decomposition subgroups of $\SO(2,n)$}}
 {Trans. Amer. Math. Soc. (to appear)}

\nextref{OhWitte-eg}
 {\sameauthor}
 {\emph{New examples of compact Clifford-Klein forms of homogeneous spaces of
$\SO(2,n)$}}
 {Internat. Math. Res. Not. {\bf 2000} (8~March 2000), no. 5, 235-251}

\nextref{OhWitte-CK}
 {\sameauthor}
 {\emph{Compact Clifford-Klein forms of homogeneous spaces of $\SO(2,n)$}}
 {(preprint)}

\nextref{Raghunathan}
 {Raghunathan,~M.~S.}
 {``Discrete Subgroups of Lie Groups"}
 {Springer-Verlag, New York, 1972}

}

\countdef\revised=100
\def\rec{Author: {\tt$\backslash$def$\backslash$rec$\{$??$\}$}}  
\def\fin{Author: {\tt$\backslash$def$\backslash$fin$\{$??$\}$}}  

 \renewcommand{\recdate}{}

\renewcommand{\address}
 {Department of Mathematics

 University of Maryland

 College Park, MD 20910 USA

\medskip

\emph{Current address:}

FIM

ETH Zentrum

CH--8092 Z\"urich Switzerland

\medskip

\emph{email:}~{\tt iozzi@math.ethz.ch}

}

\renewcommand{\addresstwo}
 {Department of Mathematics

 Oklahoma State University

 Stillwater, OK 74078 USA

\medskip

\emph{email:}~{\small\tt dwitte@math.okstate.edu}

}

\newcommand{\pref}[1]{{\upshape(}\ref{#1}{\upshape)}}
\newcommand{\see}[1]{{\upshape(}see~\ref{#1}{\upshape)}}
\newcommand{\fullref}[2]{\ref{#1}\pref{#1-#2}}
\renewcommand{\eqref}[1]{Eq.~\pref{#1}}

\renewcommand{\thesection}{\arabic{section}}

\numberwithin{equation}{section}

\newcommand{\Lie}[1]{\mathfrak{#1}}
\newcommand{\so}{\operatorname{\Lie{so}}}
\newcommand{\SL}{\operatorname{SL}}
\newcommand{\GL}{\operatorname{GL}}
\newcommand{\SO}{\operatorname{SO}}
\newcommand{\SU}{\operatorname{SU}}
\newcommand{\SP}{\operatorname{Sp}}
\newcommand{\Id}{\operatorname{Id}}
\newcommand{\real}{\mathord{\mathbb{R}}}
\newcommand{\complex}{\mathord{\mathbb{C}}}
\newcommand{\integer}{\mathord{\mathbb{Z}}}

\newcommand{\diag}{\operatorname{diag}}

\newcommand{\muH}[2]{\bigl[ #1, #2 \bigr]}

\newcommand{\Rrank}{\operatorname{\hbox{$\real$-rank}}}

\newcommand{\bigset}[2]{\left\{\, #1 
 \mathrel{\left| \vphantom {\left\{ #1 \mid #2 \right\} } \right.}
 #2 \,\right\} }

 \newcounter{case}
 \newenvironment{case}[1][\unskip]{\refstepcounter{case}\em
 \medskip \noindent Case \thecase\ #1.\ }{\unskip\upshape}
 \renewcommand{\thecase}{\arabic{case}}

 \newcounter{subcase}
 \newenvironment{subcase}[1][\unskip]{\refstepcounter{subcase}\em
 \medskip \noindent Subcase \thesubcase\ #1.\ }{\unskip\upshape}
\numberwithin{subcase}{case}

 \newcounter{subsubcase}
 \newenvironment{subsubcase}[1][\unskip]{\refstepcounter{subsubcase}\em
 \medskip \noindent Subsubcase \thesubsubcase\ #1.\ }{\unskip\upshape}
\numberwithin{subsubcase}{subcase}

\newcommand{\hypertarget}[2]{}
\newcommand{\hyperlink}[2]{#2}

\renewcommand{\Lie}[1]{\mathfrak{\lowercase{#1}}}

\renewcommand{\Re}{\operatorname{Re}}
\renewcommand{\Im}{\operatorname{Im}}
\newcommand{\cjg}[1]{\overline{#1}}

\newcommand{\fullsref}[2]{\ref{#1}\sref{#1-#2}}
\newcommand{\sref}[1]{{\upshape(}\ref{#1}*{\upshape)}}
\newcommand{\s}{\vphantom{\vrule height 15pt depth 5pt}}

\newcommand{\xx}{{\mathord{\mathsf{x}}}}
\newcommand{\yy}{{\mathord{\mathsf{y}}}}

\newcommand{\notargetlabel}{}
\let\notargetlabel=\label
 \renewcommand{\label}[1]{\hypertarget{#1}{}\notargetlabel{#1}}

\newcommand{\nolinkref}{}
\let\nolinkref=\ref
 \renewcommand{\ref}[1]{\hyperlink{#1}{\nolinkref{#1}}}

      \newtheorem{notate}[Theorem]{Notation} 
\newenvironment{Notation}{\begin{notate}\rm}{\end{notate}}

      \newtheorem{acknow}[Theorem]{Acknowledgments} 
\newenvironment{acknowledgments}{\begin{acknow}\rm}{\end{acknow}}

      \newtheorem{conjec}[Theorem]{Conjecture} 

\begin{document}
\firstpage

\begin{abstract}
 We give explicit, practical conditions that determine whether or not a closed,
connected subgroup~$H$ of $G = \SU(2,n)$ has the property that there exists a
compact subset~$C$ of~$G$ with $CHC = G$. To do this, we fix a Cartan
decomposition $G = K A^+ K$ of~$G$, and then carry out an approximate
calculation of $(KHK) \cap A^+$ for each closed, connected subgroup~$H$ of~$G$.
This generalizes the work of H.~Oh and D.~Witte for $G = \SO(2,n)$.
 \end{abstract}

\section{Introduction}

\begin{Definition}{\cite[Defn.~1.2]{OhWitte-CDS}}
 Let $H$ be a closed subgroup of a connected, simple, linear, real
Lie group~$G$. We say that $H$ is a
\emph{Cartan-decomposition subgroup} of~$G$ if 
 \begin{itemize}
 \item $H$ is connected, and
 \item there is a compact
subset~$C$ of~$G$, such that $CHC = G$. 
 \end{itemize}
 (Note that $C$ is only assumed to
be a sub\emph{set} of~$G$; it need not be a sub\emph{group}.)
 \end{Definition}

\begin{Example}
 The Cartan decomposition $G = KAK$ shows that the maximal split torus~$A$
is a Cartan-decomposition subgroup of~$G$.

It is known that $G = KNK$ \cite[Thm.~5.1]{Kostant}, so the maximal
unipotent subgroup~$N$ is also a Cartan-decomposition subgroup.

 If $\Rrank G = 0$ (that is, if $G$ is compact), then every (closed, connected)
subgroup of~$G$ is a Cartan-decomposition subgroup. 

 If $\Rrank G = 1$, then it not difficult to see that every (closed, connected)
noncompact subgroup of~$G$ is a Cartan-decomposition subgroup (cf.\
\cite[Lem.~3.2]{Kobayashi-isotropy}).
 \end{Example}

It is more difficult to characterize the Cartan-decomposition subgroups when
$\Rrank G = 2$, but H.~Oh and D.~Witte \cite{OhWitte-CDS} studied two examples
in detail. Namely, they described all the Cartan-decomposition subgroups of
$\SL(3,\real)$ and of $\SO(2,n)$, and they also explicitly described the
closed, connected subgroups that are \emph{not} Cartan-decomposition subgroups.
 Here, we obtain similar results for $\SU(2,n)$. Unfortunately, the results are
rather complicated to state.

\begin{Notation}
 Let $G = \SU(2,n)$ and fix an Iwasawa decomposition  $G = KAN$ and a
corresponding Cartan decomposition $G = K A^+ K$, where $A^+$ is the (closed)
positive Weyl chamber of~$A$ in which the roots occurring in the Lie algebra
of~$N$ are positive. Thus, $K$ is a maximal compact subgroup, $A$ is the
identity component of a maximal split torus, and $N$ is a maximal unipotent
subgroup.
 \end{Notation}

 To simplify, let us restrict our attention here to subgroups of~$N$.

\begin{Theorem}[{(cf.\ \ref{CDS<>h_m})}]
 Let $G = \SU(2,n)$ and let $H$ be a closed, connected subgroup of~$N$.
 Then $H$ is a Cartan-decomposition subgroup of~$G$ if and only if
 \begin{enumerate}
 \item $H$ satisfies at least one of the eight conditions in
Proposition~\ref{HinN-square}; and
 \item $H$ satisfies at least one of the five conditions in
Proposition~\ref{HinN-linear}.
 \end{enumerate}
 \end{Theorem}

\begin{Theorem}
 Let $G = \SU(2,n)$ and let $H$ be a closed, connected, nontrivial subgroup
of~$N$.
 Then $H$ is \textbf{not} a Cartan-decomposition subgroup of~$G$ if and only if
$H$ belongs to one of the eleven types of subgroups explicitly described in
Theorem~\ref{HinN-notCDS}.
 \end{Theorem}

For subgroups~$H$ that are not contained in~$N$, there is no loss of generality
in assuming that $H \subset AN$ \see{H'inAN}, and that $H$ satisfies the
additional technical condition of being compatible with~$A$
\see{conj-to-compatible}. Under these assumptions, Theorem~\ref{TU-proj},
Proposition~\ref{HnotTUproj}, and Lemma~\ref{dim1-mu(H)}, taken together, list
the possibilities for~$H$ and, in each case, determine whether $H$ is a
Cartan-decomposition subgroup or not.

Our results require an effective method to determine whether a subgroup is a
Cartan-decomposition subgroup or not. This is provided by the Cartan projection.

\begin{Definition} (Cartan projection)
 For each element~$g$ of~$G$, the Cartan decomposition $G = K A^+ K$ implies
that there is an element~$a$ of~$A^+$ with $g \in K a K$. In fact, the
element~$a$ is unique, so there is a well-defined function
 $$ \mbox{$\mu \colon G \to A^+$ given by $g \in K \, \mu(g) \, K$.} $$
 The function $\mu$ is continuous and proper (that is, the inverse image of any
compact set is compact). Some properties of the Cartan projection are discussed
in~\cite{Benoist} and~\cite{Kobayashi-survey}.
 \end{Definition}

We have $\mu(H) = A^+$ if and only if $KHK = G$. This immediately implies
that if $\mu(H) = A^+$, then $H$ is a Cartan-decomposition subgroup.
Y.~Benoist and T.~Kobayashi proved the deeper statement that, in the
general case, $H$ is a Cartan-decomposition subgroup if and only if
$\mu(H)$ comes within a bounded distance of every point in~$A^+$.

\begin{Notation}
 For subsets $U$ and~$V$ of~$A^+$, we write $U \approx V$ if there is a compact
subset~$C$ of~$A$, such that $U \subset V C$ and $V \subset U C$. This is
an equivalence relation.
 \end{Notation}

\begin{Theorem}[(Benoist {\cite[Prop.~5.1]{Benoist}, Kobayashi
\cite[Thm.~1.1]{Kobayashi-criterion}})] \label{CDSvsmu}
 A closed, connected subgroup~$H$ of~$G$ is a Cartan-decomposition subgroup if
and only if $\mu(H) \approx A^+$.
 \end{Theorem}

\begin{Remark}
 We may consider $\SO(2,n)$ to be the subgroup of $\SU(2,n)$ consisting of the
real matrices. Then, because $A \subset \SO(2,n)$, we see that $\SO(2,n)$ is a
Cartan-decomposition subgroup of $\SU(2,n)$. More generally, a subgroup of
$\SO(2,n)$ is a Cartan-decomposition subgroup of $\SO(2,n)$ if and only if it
is a Cartan-decomposition subgroup of $\SU(2,n)$. (For example, this follows
from the fact that the Cartan projection for $\SO(2,n)$ is the restriction of
the Cartan projection for $\SU(2,n)$.) Thus, our results generalize those
theorems of H.~Oh and D.~Witte \cite{OhWitte-CDS} that are directed toward
$\SO(2,n)$.
 \end{Remark}

\begin{Remark}
 One may define a partial order $\ll$ on the set of closed, connected subgroups
of~$G$ by
 $$ \mbox{$H_1 \prec H_2$ if there is a compact subset~$C$ of~$G$, such that
$H_1 \subset C H_2 C$} .$$
 (So $H$ is a Cartan-decomposition subgroup of~$G$ if
and only if $G \prec H$.) We see from \cite[Prop.~5.1]{Benoist} that $H_1 \prec
H_2$ if and only if there is a compact subset~$C$ of~$A$, such that $\mu(H_1)
\subset \mu(H_2) C$. Thus, it is of interest to calculate $\mu(H)$, for each
subgroup~$H$ of~$G$. Our results solve this problem: for each (closed,
connected) subgroup~$H$, we give an explicit subset~$U$ of~$A^+$, such that
$\mu(H) \approx U$. For the cases where $\mu(H) \not\approx A^+$, these results
are summarized in Tables \ref{CprojN-summaryfigure},
\ref{Cprojsemi-summaryfigure}, and~\ref{Cprojnotsemi-summaryfigure} of
Section~\ref{CDSsummary}, and the subset~$U$ is given in a standard form that
makes it easy to determine whether $H_1 \prec H_2$. Thus, we determine the
order structure of the relation~$\prec$, and also determine precisely where
each subgroup lies in this partial order.
 \end{Remark}

The interest in  Cartan-decomposition subgroups is largely due to
the following basic observation that, to construct nicely behaved actions
on homogeneous spaces, one must find subgroups that are \emph{not}
Cartan-decomposition subgroups. (See \cite[\S3]{Kobayashi-survey} for
some historical background on this result.)

\begin{Proposition}[({Calabi-Markus phenomenon, cf.~\cite[pf.\ of
Thm.~A.1.2]{KulkarniCK}})] \label{Calabi-Markus}
 \ \\  
  If $H$ is a Cartan-decomposition subgroup of~$G$, then no closed,
noncompact subgroup of~$G$ acts properly on $G/H$.
 \end{Proposition}

H.~Oh and D.~Witte \cite{OhWitte-eg, OhWitte-CK} used this proposition as a
starting point to study the existence of tessellations. (A homogeneous space
$G/H$ is said to have a tessellation if there is a discrete subgroup~$\Gamma$
of~$G$, such that $\Gamma$ acts properly on $G/H$, and $\Gamma \backslash G/H$
is compact.) In particular, when $n$~is even, they determined exactly which
homogeneous spaces $\SO(2,n)/H$ have a tessellation (under the assumption that
$H$ is connected). These results depend not only on the characterization of
Cartan-decomposition subgroups, but also on the calculation of $\mu(H)$ for
each subgroup~$H$, and on the maximum possible dimension of subgroups with a
given image under the Cartan projection. In \cite{IozziWitte-tess} we use some
of the results of the current paper to study tessellations of homogeneous
spaces of $\SU(2,n)$.

Here is an outline of the paper.
 Section~\ref{SU2n-section} describes the notation we use to specify elements
of $\SU(2,n)$.
 Section~\ref{CDS-section} recalls some general results on Cartan-decomposition
subgroups, and defines a representation~$\rho$.
 Section~\ref{HinNsquare-section} determines
whether $H$ contains large elements with $\|\rho(h)\|$ approximately equal to
$\|h\|^2$.
 Similarly,  Section~\ref{HinNlinear-section} determines
whether $H$ contains large elements with $\|\rho(h)\|$ approximately equal to
$\|h\|$.
 By combining the calculations of the preceding two sections,
Section~\ref{CDSinN-section} determines which subgroups of~$N$ are
Cartan-decomposition subgroups.
 Then Section~\ref{CDSnotinN-section} determines which other subgroups of~$G$
are Cartan-decomposition subgroups.
 Section~\ref{CDSsummary} determines the maximum possible
dimension of a subgroup of~$H$ with any given image under the Cartan projection.

\begin{acknowledgments}
 This research was partially supported by a grant from the National Science
Foundation (DMS-9801136).
 Much of the work was carried out during productive visits to the
University of Bielefeld (Germany) and the Isaac Newton Institute for
Mathematical Sciences (Cambridge, U.K.). We would like to thank the
German-Israeli Foundation for Research and Development for financial support
that made the visit to Bielefeld possible.
 D.W.\ would also like to thank the mathematics department of the University of
Maryland for its hospitality during the visit that initiated this project.
 \end{acknowledgments}

\section{Explicit coordinates in $\SU(2,n)$} \label{SU2n-section}

\begin{Notation} \label{SU2n-defn}
 We realize $\SU(2,n)$ as isometries of the indefinite Hermitian form 
 $$\langle v
\mid w \rangle
 = v_1 \cjg{w_{n+2}} + v_2 \cjg{w_{n+1}}
 + \sum_{i=3}^n v_i \cjg{w_i}
 + 
v_{n+1} \cjg{w_2} + v_{n+2} \cjg{w_1}$$
 on $\complex^{n+2}$.
 The virtue of this particular realization is that we may choose $A$ to
consist of the diagonal matrices in $\SU(2,n)$ that have nonnegative real
entries, and $N$ to consist of the upper-triangular matrices in
$\SU(2,n)$ with only $1$'s on the diagonal.
 Thus, the Lie algebra of $AN$ is
 \begin{equation}
 \label{LieAN}
 \Lie a + \Lie n =
 \bigset{
 \begin{pmatrix}
 t_1 & \phi & x & \eta & i\xx \\
 0   & t_2  & y & i\yy & -\cjg{\eta} \\
 0   & 0    & 0 & -y^{\dagger} & -x^{\dagger} \\
 0   & 0   & 0 & -t_2 & -\cjg{\phi} \\
 0   & 0   & 0 & 0 & -t_1 \\
 \end{pmatrix}
  }{ 
 \begin{matrix}
 t_1,t_2 \in \real, \\
 \phi,\eta \in \complex, \\
 x,y \in \complex^{n-2}, \\
 \xx, \yy \in \real
 \end{matrix}}
 ,
 \end{equation}
 where $\cjg{\phi}$ or~$\cjg{\eta}$ denotes the conjugate of a complex
number~$\phi$ or~$\eta$, and $x^{\dagger}$ or~$y^{\dagger}$ denotes the
conjugate-transpose of a row vector~$x$ or~$y$.
 Note that  the first two rows of any element of $\Lie a + \Lie n$ are
sufficient to determine the entire matrix.
 \end{Notation}

\begin{Notation}
 Because the exponential map is a diffeomorphism from~$\Lie N$ to~$N$, each
element of~$N$ has a unique representation in the form $\exp u$ with $u \in
\Lie N$. Thus, each element~$h$ of~$N$ determines corresponding values of
$\phi$, $x$, $y$, $\eta$, $\xx$ and~$\yy$ (with $t_1 = t_2 = 0$). We write
 $$ \phi_h, x_h, y_h, \eta_h, \xx_h, \yy_h $$
 for these values.
 \end{Notation}

\begin{Notation}
 We let $\alpha$ and~$\beta$ be the simple real roots of $\SU(2,n)$,
defined by $\alpha(a) = a_1/a_2$ and $\beta(a) = a_2$, for an element~$a$
of~$A$ of the form
 $$ a= \diag (a_1, a_2 , 1,1,\ldots,1,1 , a_2^{-1} ,
a_1^{-1}) . $$
 Thus, 
 \begin{itemize}
 \item the root space $\Lie u_\alpha$ is the $\phi$-subspace in~$\Lie n$, 
 \item the root space $\Lie u_\beta$ is the $y$-subspace in~$\Lie n$, 
 \item the root space $\Lie u_{\alpha+\beta}$ is the $x$-subspace in~$\Lie n$,
 \item  the root space $\Lie u_{\alpha+2\beta}$ is the $\eta$-subspace
in~$\Lie n$,
 \item  the root space $\Lie u_{2\beta}$ is the $\yy$-subspace
in~$\Lie n$, and
 \item  the root space $\Lie u_{2\alpha+2\beta}$ is the $\xx$-subspace
in~$\Lie n$.
 \end{itemize}
 \end{Notation}

\begin{Notation}
 For a given Lie algebra~$\Lie h \subset \Lie n$, we use $\Lie z$ to
denote $\Lie h \cap (\Lie U_{\alpha+2\beta} + \Lie U_{2\alpha+2\beta} +
\Lie U_{2\beta})$. In other words, 
 $$ \Lie z = \{\, u \in \Lie H \mid \mbox{$\phi_u = 0$ and $x_u = y_u = 0$}
\,\}. $$
 (We remark that if $\phi_u = 0$ for every $u \in \Lie H$, then $[\Lie
H, \Lie H] \subset \Lie Z$ and $\Lie Z$ is contained in the center of~$\Lie H$.)
 \end{Notation}

\begin{Notation}
 For $h \in \SU(2,n)$, define 
 $$\Delta(h) = 
 \det \begin{pmatrix}
 h_{1,n+1} & h_{1,n+2} \cr
 h_{2,n+1} & h_{2,n+2} \cr
 \end{pmatrix} .$$
 \end{Notation}

The following results collect some straightforward calculations that will be
used repeatedly throughout the paper.

\begin{Remark}
 For 
 $$ u = 
 \begin{pmatrix}
 0 & \phi & x & \eta & i \xx \\
 0 & 0 & y & i \yy & -\cjg{\eta} \\
 0 & 0 & 0 & -y^{\dagger} & -x^{\dagger} \\
 0 & 0 & 0 & 0 & -\cjg{\phi} \\
 0 & 0 & 0 & 0 & 0 \\
 \end{pmatrix}
 \in \Lie N
 \hbox{\qquad and\qquad $h = \exp u \in N$,}
 $$
 we have
 $$  \exp(u)= 
 \begin{pmatrix}
 \s 1 & \phi & x+\frac{1}{2} \phi y 
      &  
           \begin{matrix}
              \eta -\frac{1}{2} x y^{\dagger}
                \vphantom{\vrule height 8pt depth 8pt} \\
              {} + \frac{1}{2} i \phi \yy - \frac{1}{6} \phi |y|^2
                \vphantom{\vrule height 8pt depth 8pt} 
           \end{matrix}
      &
 \begin{matrix}
 - \frac{1}{2} |x|^2
 - \Re(\phi \cjg{\eta})
 +\frac{1}{24} |\phi|^2 |y|^2
                \vphantom{\vrule height 8pt depth 8pt}  \\
 {} + i \left(
 \xx
 - \frac{1}{6} |\phi|^2 \yy
 + \frac{1}{3} \Im(\cjg{\phi} x y^{\dagger}) 
 \right)
                \vphantom{\vrule height 8pt depth 8pt} 
 \end{matrix} \\
 \s 0 & 1 & y
    & i \yy -\frac{1}{2} |y|^2
    & -\cjg{\eta} -\frac{1}{2} y x^{\dagger}-\frac{1}{2} i \cjg{\phi}
\yy +\frac{1}{6} \cjg{\phi} |y|^2 \\
 \s 0 & 0 &\Id& -y^{\dagger} & -x^{\dagger}+\frac{1}{2} \cjg{\phi}
y^{\dagger} \\
 \s 0 & 0 & 0 & 1 & -\cjg{\phi} \\
 \s 0 & 0 & 0 & 0 & 1 \\
 \end{pmatrix}
 $$
 and
 $$ \Delta(h)
 = \begin{matrix}
  - |\eta|^2
 +\xx \yy
 - \frac{1}{4} |x|^2 |y|^2
 + \frac{1}{4} |x y^{\dagger}|^2
 -\frac{1}{6} |y|^2 \Re( \eta \cjg{\phi} )
  \vphantom{\vrule height 10pt depth 10pt}
 \\
 {} -\frac{1}{6} \yy \Im( xy^{\dagger} \cjg{\phi} )
 +\frac{1}{12} \yy ^2 |\phi|^2
 -\frac{1}{144} |y|^4 |\phi|^2
  \vphantom{\vrule height 10pt depth 10pt}
 \\
 {} + i \left(
 \frac{1}{24} \yy |\phi|^2 |y|^2
 +\Im(xy^{\dagger}  \cjg{\eta} )
 +\frac{1}{2} \xx |y|^2
 +\frac{1}{2} \yy |x|^2 
 \right)
 .
  \vphantom{\vrule height 10pt depth 10pt}
 \end{matrix}
 $$

 When $\phi = 0$, these simplify to:
 $$  \exp(u)= 
 \begin{pmatrix}
 \s 1 & 0 & x &  \eta -\frac{1}{2} x y^{\dagger}
      & i \xx - \frac{1}{2} |x|^2 \\
 \s 0 & 1 & y
    & i \yy -\frac{1}{2} |y|^2
    & -\cjg{\eta} -\frac{1}{2} y x^{\dagger} \\
 \s 0 & 0 &\Id& -y^{\dagger} & -x^{\dagger} \\
 \s 0 & 0 & 0 & 1 & 0 \\
 \s 0 & 0 & 0 & 0 & 1 \\
 \end{pmatrix}
 $$
 and
 $$ \Delta(h)
 = \begin{matrix}
  - |\eta|^2
 +\xx \yy
 - \frac{1}{4} |x|^2 |y|^2
 + \frac{1}{4} |x y^{\dagger}|^2
     \vphantom{\vrule height 8pt depth 8pt}  \\
 {} + i \left(
 \Im(xy^{\dagger}  \cjg{\eta} )
 +\frac{1}{2} \xx |y|^2
 +\frac{1}{2} \yy |x|^2 
 \right)
 .
     \vphantom{\vrule height 8pt depth 8pt}  \end{matrix}
 $$

Similarly, when $y = 0$, we have
 $$  \exp(u)= 
 \begin{pmatrix}
 \s 1 & \phi & x
      &  \eta +  \frac{1}{2} i \phi \yy
      & \vphantom{\vrule height 20pt depth 20pt}
 \begin{matrix}
 - \frac{1}{2} |x|^2 
 - \Re(\phi \cjg{\eta})
   \\
 {} + i \left(
 \xx
 - \frac{1}{6} |\phi|^2 \yy
 \right)
 \end{matrix} \\
 \s 0 & 1 & 0
    & i \yy
    & -\cjg{\eta} -\frac{1}{2} i \cjg{\phi} \yy \\
 \s 0 & 0 &\Id& 0 & -x^{\dagger} \\
 \s 0 & 0 & 0 & 1 & -\cjg{\phi} \\
 \s 0 & 0 & 0 & 0 & 1 \\
 \end{pmatrix}
 ,$$
 and
 \begin{equation} \label{topright-y=0}
 \Delta(h)
 = \left( \xx \yy + \tfrac{1}{12} |\phi|^2 \yy^2 - |\eta|^2 \right)
 + i \left( \tfrac{1}{2} |x|^2 \yy \right)
 .
 \end{equation}
 \end{Remark}

\begin{Remark}
 For 
 \begin{equation} \label{umatrix}
 u = 
 \begin{pmatrix}
 0 & \phi & x & \eta & i \xx \\
   & 0 & y & i \yy & -\cjg{\eta} \\
   &   & \cdots \\
 \end{pmatrix}
 \mbox{\qquad and\qquad}
\tilde u = 
 \begin{pmatrix}
 0 & \tilde \phi & \tilde x & \tilde \eta & i \tilde \xx \\
   & 0 & \tilde y & i \tilde \yy & -\cjg{\tilde \eta} \\
   &   & \cdots \\
 \end{pmatrix}
 ,
 \end{equation}
 we have
 $$ [u, \tilde u] =
 \begin{pmatrix}
 0 & 0 & \phi \tilde y - \tilde \phi y
 & - x \tilde y^{\dagger} + \tilde x y^{\dagger} + i \phi \tilde \yy - i
\tilde \phi \yy
 & -2i \Im(x \tilde x^{\dagger} + \phi \cjg{\tilde \eta}
 - \tilde \phi \cjg{\eta}) \\
   & 0 & 0 & -2 i \Im (y \tilde y^{\dagger})
 &  \tilde y x^{\dagger} -  y \tilde x^{\dagger} + i \cjg{\phi} \tilde
\yy - i \cjg{\tilde \phi} \yy \\
  &  &  & \cdots \\
 \end{pmatrix}
 , $$ 
 and
 \begin{equation} \label{[u,u,u]}
 \bigl[ [u,\tilde u], \hat u \bigr] = 
 \begin{pmatrix}
 0 & 0 & 0 &
 - (\phi \tilde y - \tilde \phi y) \hat y^{\dagger}
 + 2 i \hat \phi \Im( y \tilde y^{\dagger} )
 & * \\
   & 0 & 0 & 0 & * \\
 &&& \cdots \\
 \end{pmatrix}
 .
 \end{equation}
 \end{Remark}

\section{Preliminaries on Cartan-decomposition subgroups} \label{CDS-section}

\begin{Notation}
 We employ the usual Big Oh and little oh notation: for functions
$f_1,f_2$ on~$H$, and a subset~$Z$ of~$H$, we say \emph{$f_1 = O(f_2)$
for $z \in Z$} if there is a constant~$C$, such that, for all large $z
\in Z$, we have $\|f_1(z)\| \le C \|f_2(z)\|$. (The values of each~$f_i$
are assumed to belong to some finite-dimensional normed vector space,
typically either~$\complex$ or a space of complex matrices. Which particular
norm is used does not matter, because all norms are equivalent up to a
bounded factor.) We say \emph{$f_1 = o(f_2)$ for $z \in Z$} if
$\|f_1(z)\|/\|f_2(z)\| \to 0$ as $z \to \infty$. Also, we write $f_1
\asymp f_2$ if $f_1 = O(f_2)$ and $f_2 = O(f_1)$.
 \end{Notation}

\begin{Definition}
 Define $\rho \colon \SU(2,n) \to \GL(\complex^{n+2} \wedge \complex^{n+2})$ by
$\rho(h) = h \wedge h$, so $\rho$ is the second exterior power of the standard
representation of $\SU(2,n)$. Thus, we may define $\| \rho(h) \|$ to be the
maximum absolute value among the determinants of all the $2 \times 2$
submatrices of the matrix~$h$.
 \end{Definition}

We now introduce convenient notation for describing the image of a
subgroup under the Cartan projection~$\mu$.

\begin{Notation} \label{approx}
 For functions $f_1,f_2 \colon \real^+ \to \real^+$, and a subgroup~$H$
of~$\SU(2,n)$, we write 
 $\mu(H) \approx \muH{f_1(\|h\|)}{f_2(\|h\|)}$ if, for every
sufficiently large $C > 1$, we have
 $$ \mu(H) \approx \bigset{a \in A^+}
 { C^{-1} f_1 \bigl( \|a\| \bigr) \le \|\rho(a)\| \le C f_2 \bigl( \|a\|
\bigr)} .$$
 (If $f_1$~and~$f_2$ are monomials, or other very tame functions, then
it does not matter which particular norm is used.)
 \end{Notation}

We have $A^+ = \{\, a \in A \mid a_{1,1} \ge a_{2,2} \ge 1 \,\}$,
 so, for $a \in A^+$, we have
 $$ \|a\| = a_{1,1} \le a_{1,1} \, a_{2,2} = \| \rho(a) \| \le a_{1,1}^2 =
\|a\|^2 .$$
 Thus $A^+ \approx \muH{\|h\|}{\|h\|^2}$, so, from Theorem~\ref{CDSvsmu}, we see
that $H$ is a Cartan-decomposition subgroup of~$G$ if and only if $\mu(H)
\approx \muH{\|h\|}{\|h\|^2}$. This observation, which is essentially due to
Y.~Benoist (in a much more general context, cf.\ \cite[Lem.~2.4]{Benoist}),
leads to the following result.

\begin{Proposition}[{(cf.\ \cite[Prop.~3.24]{OhWitte-CDS})}] \label{CDS<>h_m}
 \ \ 
 A closed, connected subgroup~$H$ of $\SU(2,n)$ is a
Cartan-decomposition subgroup if and only if
 \begin{enumerate}
 \item there is a sequence $\{h_m\}$ in~$H$, such that $h_m \to \infty$ as $n
\to \infty$, and $\rho(h_m) \asymp \|h_m\|^2$; and
 \item there is a sequence $\{h_m\}$ in~$H$, such that $h_m \to \infty$ as $n
\to \infty$, and $\rho(h_m) \asymp h_m$.
 \end{enumerate}
 \end{Proposition}

The following result allows us to replace $H$ by a conjugate subgroup whenever
it is convenient.

\begin{Lemma}[{(cf.\ \cite[Prop.~1.5]{Benoist},
\cite[Cor.~3.5]{Kobayashi-criterion})}]
 Let $H$ be any closed, connected subgroup of $\SU(2,n)$. For every $g \in G$,
we have $\mu( g^{-1} H g) \approx \mu(H)$.

In particular, $H$ is a Cartan-decomposition subgroup if and only if $g^{-1} H
g$ is a Cartan-decomposition subgroup.
 \end{Lemma}

\section{When is the size of $\rho(h)$ quadratic?} \label{HinNsquare-section}

In this section, Proposition~\ref{HinN-square} is a list of subgroups that
contain a sequence $\{h_m\}$ with $\rho(h_m) \asymp \|h_m\|^2$, and
Proposition~\ref{HinN-nosquareproj} is a list of subgroups that do not contain
such a sequence. Then Proposition~\ref{HinN-squareornot} shows that both lists
are complete.

\begin{Proposition} \label{HinN-square}
 Assume that $G = \SU(2,n)$. Let $H$ be a closed,
connected subgroup of~$N$.
 There is a sequence $h_m \to \infty$ in~$H$ with
$\rho(h_m) \asymp \|h_m\|^2$ if either
 \begin{enumerate}
 \item \label{HinN-square-indep}
 there is an element~$u$ of~$\Lie h$ with $\phi_u = 0$, such that the vectors
$x_u$ and~$y_u$ are linearly independent over~$\complex$; or
 \item \label{HinN-square-eta2neq}
 there is an element~$z$ of~$\Lie z$, such that $|\eta_z|^2 \neq \xx_z
\yy_z$; or
 \item \label{HinN-square-neq0}
 there are elements $u$ of~$\Lie h$ and $z$ of~$\Lie z$, such that $\phi_u =
0$, and 
 $\xx_z |y_u|^2 + \yy_z |x_u|^2 + 2 \Im ( x_u
y_u^{\dagger} \cjg{\eta_z}) \neq 0$; or
 \item \label{HinN-square-y=0+yy=0}
 there is an element~$u$ of~$\Lie h$, such that $\phi_u \neq 0$, $y_u = 0$,
$\yy_u = 0$, and $|x_u|^2 + 2 \Re(\phi_u \cjg{\eta_u}) = 0$; or 
 \item \label{HinN-square-yu=0} 
 $\Lie U_{2\alpha + 2\beta} \subset \Lie H$ and
 there is an element~$u$ of~$\Lie h$, such that $\phi_u \neq 0$, $\yy_u \neq
0$, and $y_u = 0$; or
 \item  \label{HinN-square-xvyu}
 there are elements $u$ and~$v$ of~$\Lie h$, such that $\phi_u \neq 0$, $y_u
\neq 0$,  $\phi_v = 0$, $y_v = 0$, $x_v \neq 0$, $\yy_v = 0$, and $x_v
y_u^{\dagger} = 0$; or
 \item \label{HinN-square-xx}
 $\Lie U_{2\alpha+2\beta} \subset \Lie Z$, and
 there are nonzero elements $u$ and~$v$ of~$\Lie h$, satisfying $\phi_u \neq
0$, $y_u \neq 0$, $\phi_v = 0$, $y_v = 0$, $\yy_v \neq 0$, and $x_v
y_u^{\dagger} = - i \phi_u \yy_v$; or
 \item \label{HinN-square->0}
 $\dim \Lie H = 3$, $\Lie Z = \Lie U_{2\alpha+2\beta}$, there exist $u,v \in
\Lie H \setminus \Lie Z$, such that $y_u \neq 0$, $y_v = 0$, $\yy_v = 0$,
$|x_v|^2 + 2\Re(\phi_v \cjg{\eta_v}) > 0$, and we have $\phi_h
\neq 0$ for every $h \in \Lie H \setminus \Lie Z$.
 \end{enumerate}
 \end{Proposition}

\begin{Remark}
 In Conclusions~\pref{HinN-square-xvyu} and~\pref{HinN-square-xx}, the
restriction on $x_v y_u^{\dagger}$ is not necessary; it was included to avoid
overlap with Conclusion~\pref{HinN-square-eta2neq}. Namely, if $x_v
y_u^{\dagger} \neq - i \phi_u \yy_v$, then $[u,v]$ satisfies $\yy = 0$ and
$\eta \neq 0$, so Conclusion~\pref{HinN-square-eta2neq} holds. Also, it is not
necessary to assume $\yy_v \neq 0$ in Conclusion~\pref{HinN-square-xx}, because
Conclusion~\pref{HinN-square-xvyu} holds if $\yy_v = 0$ (and $x_v \neq 0$).
Thus, \pref{HinN-square-xvyu} and~\pref{HinN-square-xx} may be replaced with
the following: 
 \begin{enumerate}
 \item[{\sref{HinN-square-xvyu}}]
 there are elements $u$ and~$v$ of~$\Lie h$, such that $\phi_u \neq 0$, $y_u
\neq 0$,  $\phi_v = 0$, $y_v = 0$, $x_v \neq 0$, and $\yy_v = 0$; or
 \item[{\sref{HinN-square-xx}}] 
 $\Lie U_{2\alpha+2\beta} \subset \Lie Z$, and
 there are nonzero elements $u$ and~$v$ of~$\Lie h$, satisfying $\phi_u \neq
0$, $y_u \neq 0$, $\phi_v = 0$, $y_v = 0$, and $x_v \neq 0$.
 \end{enumerate}
 \end{Remark}

\begin{proof} 
 We separately consider each of the eight cases in the statement of the
proposition.

\setcounter{case}{0}

\pref{HinN-square-indep}
 Let $h^t = \exp(t u)$.
 Replacing $H$ by a conjugate under~$U_\alpha$, we may
assume that $x_u$ is orthogonal to~$y_u$; that is, $x_u y_u^{\dagger} = 0$. 
Then it is clear that $\rho(h^t) \asymp \Delta(h^t) \asymp t^4 \asymp
\|h^t\|^2$.

\pref{HinN-square-eta2neq}
 Let $h^t = \exp(t z)$. We have $h^t \asymp t$ and
 $$ \Delta(h^t) = \xx_{tz} \yy_{tz} - |\eta_{tz}|^2
 = t^2( \xx_z \yy_z - |\eta_z|^2) \asymp t^2 .$$
 Therefore
 $\rho(h^t) \asymp \Delta(h^t) \asymp t^2 \asymp \|h^t\|^2$.

\pref{HinN-square-neq0}
 For any large~$t$, let $h = \exp(tu + t^2 z)$. Clearly, we have
$|x_h| + |y_h| = O(t)$ and $|\xx_h| + |\yy_h| + |\eta_h| = O(t^2)$, so $h
= O(t^2)$. 

We have
 \begin{eqnarray*}
 \Im \Delta(h^t)
 =
 t^4 \biggl[ \frac{1}{2} \bigl( 2 \Im (x_u y_u^{\dagger} \cjg{\eta_z} +
\xx_z |y_u|^2 + \yy_z |x_u|^2)
 \biggr]
 + O(t^3)
 \asymp t^4
 .
 \end{eqnarray*}
 Therefore, $\rho(h^t) \asymp t^4 \asymp \|h^t\|^2$.

\pref{HinN-square-y=0+yy=0}
 For any large~$t$, let $h = \exp(tu)$. Then $h_{1,n+2} = i t \xx_u$, so it is
easy to see that $h \asymp t$. We have $\rho(h) \asymp t^2 \asymp \|h\|^2$.

\pref{HinN-square-yu=0}
 Replacing $H$ by a conjugate (under a diagonal matrix), we may assume that
$\phi_u = \yy_u$. Then, by renormalizing, we may assume that $\phi_u = \yy_u
= 1$. Let $z$ be the element of~$\Lie U_{2\alpha + 2\beta}$ with $\xx_z = 1$.
By subtracting a multiple of~$z$ from~$u$, we may assume $\xx_u = 0$. For any
large~$t$, let 
 $h = \exp(6tu + 36 t^3 z)$, so $h_{1,n+2}$ is real. We have
 $$\Re \Delta(h) = (36 t^3) ( 6t) + \frac{1}{12}(6t)^2 (6t)^2 + O(t^2)
 \asymp t^4 ,$$
 so $\rho(h) \asymp t^4 \asymp \|h\|^2$.

\pref{HinN-square-xvyu}
 For each large~$t$, let $h$ be an element of~$\exp(tu + \real v)$, such that
$h_{1,n+2}$ is pure imaginary. (This exists because the sign of $-\frac{1}{2}
|x|^2$ is opposite that of $\frac{1}{24} |\phi|^2 |y|^2$.) We note that $x_h
\asymp t^2$ and $|\eta_h| + |\xx_h| = O(t^2)$, but
 $\phi_h \asymp y_h \asymp t$ and $|\yy_h| + |x_h y_h^{\dagger}| = O(t)$.
 Thus $h = O(t^3)$ and
 $$ \rho(h) \asymp \Re \Delta(h) = -\frac{1}{4} |x_h|^2 |y_h|^2
 - \frac{1}{144} |y_h|^4 |\phi_h|^2 + O(t^5)
 \asymp t^6 \asymp \|h\|^2 .$$

\pref{HinN-square-xx}
 Because $x_v y_u^{\dagger} = - i \phi_u \yy_v$, we have $x_v \neq 0$, so, for
any large~$t$, we may choose $h \in \exp( tu + \real v + \Lie U_{2\alpha +
2\beta})$, such that $h_{1,n+2} = 0$.
 Thus $\phi_h \asymp y_h \asymp t$, but $x_h \asymp \yy_h \asymp t^2$ and
$|\eta_h| + |\xx_h| = O(t^2)$. Then (because $h_{1,n+2} = 0$) it is easy to
verify that $h = O(t^3)$. However
 $$ \Im \Delta(h)
 = \frac{1}{24} \yy_h |\phi_h|^2 |y_h|^2 + \frac{1}{2} \yy_h |x_h|^2 + O(t^5)
 \asymp t^6 .$$
 So $\rho(h) \asymp \|h\|^2$.

\pref{HinN-square->0}
 For any large~$t$, choose $s = O(1)$, such that $\Re \bigl( \exp(su +
tv)_{1,n+2} \bigr) = 0$. (This is possible, because $-\frac{1}{2} |x_v|^2 -
\Re(\phi_v \cjg{\eta_v}) < 0$.) Then we may choose $h \in \exp(s u + tv +
\Lie Z)$, such that $h_{1,n+2} = 0$. Then $\phi_h \asymp t$, $|x_h| + |\eta_h| =
O(t)$, and $|y_h| + |\yy_h| = O(1)$, so we have
 $\rho(h) \asymp t^2 \asymp \|h\|^2$.
 \end{proof}

\begin{Proposition} \label{HinN-nosquareproj}
 Assume that $G = \SU(2,n)$. Let $H$ be a closed,
connected, nontrivial subgroup of~$N$.
 \begin{enumerate}
 \item \label{HinN-nosquareproj-H=Z}
 If $\dim \Lie H = 1$, $\Lie H = \Lie Z$, and we have $|\eta_h|^2 = \xx_h
\yy_h$ for every $h \in H$, then $\rho(h) \asymp h$ for every $h \in H$.
 \item \label{HinN-nosquareproj-phi&y=0}
 If $\phi_h = 0$ and $y_h = 0$ for every $h \in \Lie H$, $\Lie Z \subset \Lie
U_{2\alpha + 2\beta}$, and there is some $u \in \Lie H$, such that $\yy_u \neq
0$, then $\mu(H) \approx \muH{\|h\|}{\|h\|^{3/2}}$, unless $\dim H = 1$, in
which case $\rho(h) \asymp \|h\|^{3/2}$ for every $h \in H$.
 \item \label{HinN-nosquareproj-phi=0}
 Suppose $\phi_h = 0$ for every $h \in \Lie H$, and there is some $\lambda \in
\complex$, such that $x_h = \lambda y_h$ for every $h \in H$, and we have
$\eta_z = i\lambda \yy_z$ and $\xx_z = |\lambda|^2 \yy_z$ for every $z \in \Lie
Z$.
 \begin{enumerate}
 \item \label{HinN-nosquareproj-phi=0-neq0}
 If there is some $u \in \Lie H$, such that $\xx_u + |\lambda|^2 \yy_u
+ 2 \Im( \lambda \cjg{\eta_u}) \neq 0$, then $\mu(H) \approx
\muH{\|h\|}{\|h\|^{3/2}}$, unless $\dim H = 1$, in which case $\rho(h)
\asymp \|h\|^{3/2}$ for every $h \in H$.
 \item Otherwise, $\rho(h) \asymp h$ for every $h \in H$.
 \end{enumerate}
 \item \label{HinN-nosquareproj-y&yy=0}
 If $y_h = 0$, $\yy_h = 0$, and $|x_h|^2 + 2 \Re(\phi_h \cjg{\eta_h})
\neq 0$ for every $h \in \Lie H \setminus \Lie U_{2\alpha+2\beta}$ {\upshape(}so
$\Lie Z \subset \Lie U_{2\alpha+2\beta}${\upshape)}, then $\rho(h) \asymp h$
for every $h \in H$.
 \item \label{HinN-nosquareproj-phi=yy}
 If $\Lie Z = 0$, there is some $u \in \Lie H$ and some nonzero $\phi_0 \in
\complex$, such that $\phi_u \neq 0$, and we have $\phi_h = \phi_0 \yy_h$ and
$y_h = 0$, for every $h \in \Lie H$, then $\mu(H) \approx
\muH{\|h\|}{\|h\|^{4/3}}$, unless $\dim H = 1$, in which case, $\rho(h) \asymp
\|h\|^{4/3}$ for every $h \in H$.
 \item \label{HinN-nosquareproj-onlyphiyyxx}
 If $\dim \Lie H \le 3$, $\Lie Z = 0$, we have $\phi_v \asymp y_v$ and $v = O
\bigl( |\phi_v| + |\yy_v| \bigr)$ for every $v \in \Lie H$, and there exists $u
\in \Lie H$, such that $\phi_u \neq 0$, then $\rho(h) \asymp \|h\|^{3/2}$ for
every $h \in H$.
 \item \label{HinN-nosquareproj-dim3}
 If $\dim \Lie H = 2$, $\Lie Z = \Lie U_{2\alpha+2\beta}$, and $\phi_h
\neq 0$ and $y_h \neq 0$ for every $h \in \Lie H \setminus \Lie Z$, then
$\mu(H) \approx \muH{\|h\|}{\|h\|^{3/2}}$.
 \end{enumerate}
 \end{Proposition}

\begin{proof*}
 We separately consider each of the seven cases in the statement of the
proposition.

\pref{HinN-nosquareproj-H=Z}
 Because $\Delta(h) = 0$ for every $h \in H$, it is clear that $\rho(h) \asymp
h$ for every $h \in H$.

\pref{HinN-nosquareproj-phi&y=0}
 We have $|\eta_h| + |\yy_h| = O(x_h)$, so $h_{1,n+2} \asymp |x_h|^2 +
|\xx_h|$ and $h_{i,j} = O(x_h) = O \bigl( |h_{1,n+2}|^{1/2} \bigr)$ whenever
$(i,j) \neq (1,n+2)$. Thus, $\rho(h) = O \bigl( \|h\|^{3/2} \bigr)$.

We have $\rho \bigl( \exp(tu) \bigr) \asymp \Im \Delta \bigl( \exp(tu)
\bigr) \asymp t^3 \asymp \| \exp(tu) \|^{3/2}$. If $\dim H > 1$, then there
is some nonzero $v \in \Lie H$, such that $\yy_v = 0$. Then, for $h \in \exp(
\real v)$, we have $\rho(h) \asymp |x_h|^2 + |\xx_h| \asymp h$.

\pref{HinN-nosquareproj-phi=0}
 Replacing $H$ by a conjugate under~$U_\alpha$, we may assume that $\lambda =
0$, so $x_h = 0$ for every $h \in H$, and $\eta_z = \xx_z = 0$ for every $z
\in \Lie Z$ (which means $\Lie Z \subset \Lie U_{2\beta}$). Therefore, the Weyl
reflection corresponding to the root~$\alpha$ conjugates~$\Lie H$ to a
subalgebra either of type~\pref{HinN-nosquareproj-phi&y=0} or of
type~\pref{HinN-nosquareproj-y&yy=0}, depending on whether or not there is some
$u \in \Lie H$, such that $\xx_u + |\lambda|^2 \yy_u + 2 \Im( \lambda
\cjg{\eta_u}) \neq 0$.

\pref{HinN-nosquareproj-y&yy=0}
 By assumption, the quadratic form $|x|^2 + 2 \Re(\phi \cjg{\eta})$ is definite
on $\Lie H / \Lie Z$, so
 $|x|^2 + |\phi|^2 + |\eta|^2 = O \bigl( |x|^2 + 2 \Re(\phi \cjg{\eta})
\bigr)$. Therefore, $h_{i,j} = O( |h_{1,n+2}|^{1/2} )$ whenever $(i,j) \neq
(1,n+2)$. Furthermore, $h_{i,j} = O(1)$ whenever $i \neq 1$ and $j \neq
n+2$. Thus, $\rho(h) \asymp h$.

\pref{HinN-nosquareproj-phi=yy}
 For any sequence $\{h_m\} \to \infty$ in~$H$, we write
$\phi_m,x_m, y_m, \yy_m, \eta_m, \xx_m$ for $\phi_{h_m}$, etc.

 We have $\phi_m \asymp \yy_m$. If $x_m = O(|\yy_m|^{3/2})$, then
 $\rho(h_m) \asymp \Re \Delta(h_m) \asymp \yy_m^4 \asymp \|h_m\|^{4/3}$. (This
completes the proof if $\dim H = 1$.)
 If $|\yy_m|^{3/2} = o(x_m)$, then $h_m \asymp h_{1,n+2} \asymp |x_m|^2$, but
 $h_{i,j} = O \bigl( |x_m| + \yy_m^2 \bigr) = O \bigl( |x_m|^{4/3} \bigr)$
 whenever $(i,j) \neq (1,n+2)$, and $h_{i,j} = O(\yy_m) = O \bigl(
|x_m|^{2/3} \bigr)$ whenever $i \neq 1$ and $j \neq n+2$.
 Therefore
 $$ \rho(h_m) = O \bigl( |x_m|^2 |x_m|^{2/3} +  |x_m|^{4/3}  |x_m|^{4/3} \bigr)
 = O \bigl(  |x_m|^{8/3} \bigr)
 = O \bigl( \|h_m\|^{4/3} \bigr) .$$

If $\dim H > 1$, then there is some (large) $h \in H$ with $\yy_h = 0$
(and hence $\phi_h = 0$). Thus $\rho(h) \asymp |x_h|^2 \asymp h$.

\pref{HinN-nosquareproj-onlyphiyyxx}
  For any sequence $\{h_m\} \to \infty$ in~$H$, we show that $\rho(h_m)
\asymp \Delta(h_m) \asymp \|h_m\|^{3/2}$. 
 We write $\phi_m,x_m, y_m, \yy_m, \eta_m, \xx_m$ for $\phi_{h_m}$, etc.

 If $\yy_m = o(\phi_m^2)$, then $h_{1,n+2} \asymp \phi_m^4$, but $h_{i,j} =
O(\phi_m^3)$ whenever $(i,j) \neq (1,n+2)$, and $h_{i,j} = O(\phi_m^2)$ whenever
$i \neq 1$ and $j \neq n+2$. Thus, $\rho(h_m) \asymp \Re \Delta(h_m) \asymp
\phi_m^6 \asymp \|h_m\|^{3/2}$.

We may now assume that $\phi_m^2 = O(\yy_m)$. Thus, there is some $v \in \Lie
H$, such that $\phi_v = 0$ and $\yy_v = 1$. (Note that, because $y_v
\asymp \phi_v$, we have $y_v = 0$.) Because $[u,v] \in \Lie Z = 0$, we must
have $\eta_{[u,v]} = 0$, so $x_v y_u^{\dagger} = - i \phi_u \yy_v \neq 0$. In
particular, $x_v \neq 0$, so $x_m \asymp \yy_m$.

We have $h_{1,n+2} = O \bigl( |x_m|^2 \bigr) = O( \yy_m^2 )$, but $h_{i,j} =
O \bigl( |\phi_m \yy_m| + |\yy_m| \bigr) = O \bigl( |\yy_m|^{3/2} \bigr)$
whenever $(i,j) \neq (1,n+2)$, and $h_{i,j} = O(\yy_m)$ whenever $i \neq 1$
and $j \neq n+2$. Thus, $h_m = O( \yy_m^2)$ and $\rho(h_m) =O(\yy_m^3)$. 

Furthermore, we have
 $$ \Im \Delta(h_m) = \frac{1}{24} \yy_m |\phi_m|^2 |y_m|^2 + \frac{1}{2}
\yy_m |x_m|^2 + O(\yy_m^2 \phi_m)
 \asymp \yy_m ^3 ,$$
 because $\yy_m |x_m|^2 \asymp \yy_m ^3$, and the terms $\frac{1}{24} \yy_m
|\phi_m|^2 |y_m|^2$ and $\frac{1}{2} \yy_m |x_m|^2$ cannot cancel (since they
both have the same sign
as~$\yy_m$). We conclude that $\rho(h_m) \asymp \Delta(h_m) \asymp \yy_m^3$.

All that remains is to show $\yy_m^2 = O( h_m )$.
 If $\phi_m^2 = o(\yy_m)$, then 
 $$\Re h_{1,n+2} = - \frac{1}{2} |x|^2 + O (
\phi_m^2 \yy) \asymp \yy_m^2 ,$$
 as desired.
 If $\yy_m = o(\phi_m^2)$, then 
 $$ \Re h_{1,n+2} \asymp o(\phi_m^4) + o(\phi_m^3) + |\phi_m^4|
 \asymp \phi_m^4 ,$$
 so $\yy_m = o(\phi_m^2) = o(\phi_m^4) = o(h_m)$,
 as desired.
 Thus, we may assume that $\yy_m \asymp \phi_m^2$. Because $x_m = \yy_m
x_v + O(\phi_m)$ and
 $x_v y_m^{\dagger} = - i \phi_m \yy_v = - i \phi_m$, we have
 \begin{eqnarray*}
 \Im(h_{1,n+2})
 &=& O(\yy_m) - \frac{1}{6} |\phi_m|^2 \yy_m
 + \left[ \frac{1}{3} \Im \bigl( \cjg{\phi_m} (\yy_m x_v) y_m^{\dagger} \bigr)
 + O(\phi_m^3) \right] \\
 &=& -\frac{1}{6} |\phi_m|^2 \yy_m
 - \frac{1}{3} |\phi_m|^2 \yy_m
 + O(\phi_m^3)
 \asymp \yy_m^2,
 \end{eqnarray*}
 as desired.

\pref{HinN-nosquareproj-dim3}
 For $z \in \Lie Z$, we have $\rho(z) \asymp z$. For $u \in \Lie H \setminus
\Lie Z$ with $y_u \neq 0$, we have $\rho \bigl( \exp(tu) \bigr) \asymp t^6
\asymp \| \exp(tu)\|^{3/2}$. All that remains is to show $\rho(h) =
O \bigl( \|h\|^{3/2} \bigr)$ for every $h \in H$.

Note that $\phi_h \asymp y_h$, and $|x_h| + |\eta_h| + |\yy_h| = O(\phi_h)$. If
$\phi_h = O(1)$, then it is obvious that $\rho(h) \asymp h$. Thus, we may assume
$|\phi_h| \to \infty$. Then, because 
 $\Re h_{1,n+2} \asymp
|\phi_h|^2 |y_h|^2 \asymp \phi_h^4$,
 but
 $h_{i,j} = O \bigl( \phi_h | y_h|^2 \bigr)
= O(\phi_h^3)$ whenever $(i,j) \neq (1,n+2)$, and 
 $h_{i,j} = O(\phi_h^2)$ whenever $i \neq 1$ and $j \neq n+2$, we have
  $$\rho(h) = O\bigl[ |\phi_h|^4 |\phi_h|^2  + \bigl(
|\phi_h|^3 \bigr)^2 \bigr]
 = O \bigl( |\phi_h|^6 \bigr)
 = O \bigl( |h_{1,n+2}|^{3/2} \bigr)
 = O \bigl( \|h\|^{3/2} \bigr) .
 \qeddis 
  \end{proof*}

\begin{Proposition} \label{HinN-squareornot}
 Assume that $G = \SU(2,n)$. Let $H$ be a closed, connected, nontrivial
subgroup of~$N$. 
 \begin{enumerate}
 \item \label{HinN-squareornot-yes}
 There is a sequence $h_m \to \infty$ in~$H$
with $\rho(h_m) \asymp \|h_m\|^2$ if and only if $H$ is one of the subgroups
described in Proposition~\ref{HinN-square}.
 \item \label{HinN-squareornot-no}
 There is \textbf{not} a sequence $h_m \to \infty$ in~$H$
with $\rho(h_m) \asymp \|h_m\|^2$ if and only if $H$ is one of the subgroups
described in Proposition~\ref{HinN-nosquareproj}.
 \end{enumerate}
 \end{Proposition}

\begin{proof} 
 It suffices to show that $H$ is described in either
Proposition~\ref{HinN-square} or Proposition~\ref{HinN-nosquareproj}.

 We may assume 
 \begin{equation} \label{nosquare-eta2=xxyy}
 \text{\rm $|\eta_z|^2 = \xx_z \yy_z$ for every $z \in \Lie Z$}
 \end{equation}
(otherwise, \fullref{HinN-square}{eta2neq} holds). Because $|\eta|^2 - \xx
\yy$ is a quadratic form of signature $(3,1)$ on $\Lie U_{2\beta} + \Lie
U_{\alpha+2\beta} + \Lie U_{2\alpha+2\beta}$, then we must have $\dim \Lie Z
\le 1$. Thus, we may assume $\Lie H \neq \Lie Z$ (otherwise
\fullref{HinN-nosquareproj}{H=Z} holds).

\setcounter{case}{0}

\begin{case}
 Assume $\phi_h = 0$ and $y_h = 0$ for every $h \in H$ {\upshape(}and $\Lie H
\neq \Lie Z${\upshape)}.
 \end{case}
 We may assume $\yy_z = 0$ for every $z \in \Lie Z$, for, otherwise,
\fullref{HinN-square}{neq0} holds. Then, from~\eqref{nosquare-eta2=xxyy}, we
have $\eta_z = 0$ for every $z \in \Lie Z$. Thus, $\Lie Z \subset \Lie
U_{2\alpha+2\beta}$. We may assume $\yy_h = 0$ for every $h \in H$, for
otherwise Conclusion~\fullref{HinN-nosquareproj}{phi&y=0} holds. We conclude
that \fullref{HinN-nosquareproj}{y&yy=0} holds.

\begin{case}
 Assume $\phi_h = 0$ for every $h \in H$, and there is some $u \in \Lie H$
with $y_u \neq 0$.
 \end{case}
 We may assume that $x_h$ and~$y_h$ are linearly dependent over~$\complex$
for every $h \in H$ (otherwise \fullref{HinN-square}{indep} holds). In
particular, there exists $\lambda \in \complex$, such that $x_u = \lambda
y_u$.

\begin{subcase}
 Assume $\Lie Z = 0$.
 \end{subcase}

\begin{subsubcase} \label{HinN-squareornot-subsub<xy>not<y>}
 Assume there exists $v \in \Lie H$, such that either $x_v \notin \complex
y_u$ or $y_v \notin \complex y_u$.
 \end{subsubcase}
 We may assume there exists $w \in\Lie H$, such that $x_w \neq \lambda y_w$
(otherwise \fullref{HinN-nosquareproj}{phi=0} holds). Furthermore, by adding
a small linear combination of~$u$ and~$v$ to~$w$, we may assume that $y_w
\neq 0$ and that either $x_w \notin \complex y_u$ or $y_w \notin \complex
y_u$. Because $x_w$ and~$y_w$ are linearly dependent, there exists
$\lambda_1$ ($\neq \lambda$) such that $x_w = \lambda_1 y_w$. (Then note
that we must have $y_w \notin \complex y_u$.) Then
 $$ x_{u+w} = x_u + x_w = \lambda y_u + \lambda_1 y_w \notin \complex (y_u +
y_w) = \complex y_{u+w} $$
 (because $\lambda \neq \lambda_1$ and $\{y_u,y_w\}$ is linearly independent
over~$\complex$). This contradicts the fact that $x_{u+w}$ and~$y_{u+w}$ are
linearly dependent over~$\complex$.

\begin{subsubcase}
 Assume $x_h,y_h \in \complex y_u$, for every $h \in \Lie H$. 
 \end{subsubcase}
 For each $h \in \Lie H$, there exist $\lambda_x,\lambda_y \in \complex$,
such that $x_h = \lambda_x y_u$ and $y_h = \lambda_y y_u$. Because $\Lie Z =
0$, we must have $\yy_{[h,u]} = 0$, so $\Im(y_h y_u^{\dagger}) = 0$, which
means that $\lambda_y$ is real. We must also have $\eta_{[h,u]} = 0$, so
 $$ 0 = -x_h y_u^{\dagger} + x_u y_h^{\dagger} = (-\lambda_x + \lambda
\cjg{\lambda_y}) |y_u|^2
 = (-\lambda_x + \lambda \lambda_y) |y_u|^2 .$$
 Thus $\lambda_x = \lambda \lambda_y$, so
 $$x_h = \lambda_x y_u = \lambda \lambda_y y_u 
 = \lambda y_h .$$
 Therefore \fullref{HinN-nosquareproj}{phi=0} holds.

\begin{subcase}
 Assume $\Lie Z \neq 0$.
 \end{subcase}
 We show that either \fullref{HinN-square}{eta2neq},
\fullref{HinN-square}{neq0} or~\fullref{HinN-nosquareproj}{phi=0} holds. 
 Straightforward calculations show that conditions
\fullref{HinN-square}{eta2neq}, \fullref{HinN-square}{neq0}
and~\fullref{HinN-nosquareproj}{phi=0} are invariant under conjugation
by~$U_\alpha$, so we may assume that $\lambda = 0$; that is, $x_u = 0$. Thus,
we may assume $\xx_z = 0$ for every $z \in \Lie Z$, for, otherwise,
\fullref{HinN-square}{neq0} holds. Then we may assume $\eta_z = 0$ for every
$z \in \Lie Z$, for, otherwise, \fullref{HinN-square}{eta2neq} holds;
therefore $\Lie Z = U_{2\beta}$. We may now assume $x_h = 0$ for every $h \in
\Lie H$, for, otherwise,  \fullref{HinN-square}{neq0} holds. Thus,
\fullref{HinN-nosquareproj}{phi=0} holds (with $\lambda = 0$).

\begin{case}
 Assume there exists $u \in \Lie H$ with $\phi_u \neq 0$.
 \end{case}
 We claim that $\Lie Z \subset \Lie U_{2\alpha+2\beta}$. If not, then there is
some $z \in \Lie Z$, such that either $\eta_z \neq 0$ or $\yy_z \neq 0$. If
$\yy_z = 0$, then $|\eta_z|^2 \neq 0 = \xx_z \yy_z$, so
\fullref{HinN-square}{eta2neq} holds. On the other hand, if $\yy_z \neq 0$,
then, letting $z' = [u,z]$, we have $\yy_{z'} = 0$ and $\eta_{z'} \neq 0$, so
\fullref{HinN-square}{eta2neq} holds once again.

\begin{subcase}
 Assume $y_h = 0$ for every $h \in \Lie H$.
 \end{subcase}
 We may assume that there is some $v \in \Lie H$, such that $\yy_v \neq 0$
(otherwise, either \fullref{HinN-square}{y=0+yy=0}
or~\fullref{HinN-nosquareproj}{y&yy=0} holds). Then we may assume $\Lie Z = 0$
(otherwise, \fullref{HinN-square}{yu=0} holds). 

We claim that
\fullref{HinN-nosquareproj}{phi=yy} holds. If not, then there is some $w \in
\Lie H$, such that $\phi_w \neq 0$ and $\yy_w = 0$. Then $\eta_{[v,w]} \neq
0$, which contradicts the assumption that $\Lie Z = 0$.

\begin{subcase} 
 Assume there is some $v \in \Lie H$, such that $y_v \neq 0$.
 \end{subcase}

\begin{subsubcase}
 Assume $\Lie Z = \Lie U_{2\alpha+2\beta}$.
 \end{subsubcase}
 Suppose, for the moment, that there exists $w \in \Lie H \setminus \Lie Z$ with
$\phi_w = 0$. We may assume that $y_w = 0$ (otherwise,
\fullref{HinN-square}{neq0} holds). Therefore $x_w \neq 0$, so
\fullsref{HinN-square}{xx} holds.

We may now assume that $\phi_w \neq 0$ for every $w \in \Lie H \setminus \Lie
Z$. This implies that $x$, $y$, $\eta$, and~$\yy$ are functions of~$\phi$;
in particular, $\dim \Lie H \le 3$. Also, because $\Lie Z \neq 0$ and $u,v
\notin \Lie Z$, we must have $\dim \Lie H \ge 2$.

We claim $\dim \Lie H = 2$ (so \fullref{HinN-nosquareproj}{dim3} holds). If
not, then $\dim \Lie H = 3$, so there exist $u,w \in \Lie H$, such that
$\phi_u = 1$ and $\phi_w = i$. Because $\phi_{[u,w]} = 0$, we must have
$[u,w] \in \Lie U_{2\alpha+2\beta}$.
 Therefore
 $ 0 = x_{[u,w]} = y_w - i y_u$, so $y_w = i y_u$.
 Furthermore,
 $$ 0 = \yy_{[u,w]} = -2 i \Im(y_u y_w^{\dagger}) = -2i \Im \bigl( -i
|y_u|^2 \bigr) = -2 i |y_u|^2 ,$$
 so $y_u = 0$. Then $y_w = i y_u$ is also~$0$.
 This implies $y_h = 0$ for every $h \in H$. This contradicts the fact that
$y_v \neq 0$.

\begin{subsubcase}
 Assume $\Lie Z = 0$.
 \end{subsubcase}
 Lemma~\ref{HinN-Z=0} below implies that either
\fullref{HinN-nosquareproj}{onlyphiyyxx} or~\fullsref{HinN-square}{xvyu}
holds.
 \end{proof}

\begin{Lemma} \label{HinN-Z=0}
 Let $H$ be a closed, connected subgroup of~$N$, such that $\Lie Z = 0$, and
assume there exist $u,v \in \Lie H$, such that $\phi_u \neq 0$ and $y_v \neq
0$. Then either $H$ is described
in~\fullref{HinN-nosquareproj}{onlyphiyyxx} {\upshape(}and
in~\fullref{HinN-nolinearproj}{onlyphiyyxx}, which is the same{\upshape)}, or
$H$ is a a Cartan-decomposition subgroup {\upshape(}and is described
in~\fullsref{HinN-square}{xvyu}
and~\fullref{HinN-linear}{phi=0&=0}{\upshape)}.
 \end{Lemma}

\begin{proof}
 Let us begin by establishing that $\phi_h \asymp y_h$ for $h \in \Lie H$.
 If not, then we may assume either that $y_u = 0$ or that $\phi_v = 0$.
Then, because $\bigl[ [u,v],v \bigr] \in \Lie Z = 0$, we see
from~\eqref{[u,u,u]} that
 $$ 0 = -(\phi_u y_v - \phi_v y_u) y_v^{\dagger} + 2i \phi_v \Im (y_u
y_v^{\dagger})
 = - \phi_u |y_v|^2 - 0 + 0 \neq 0.$$
 This contradiction establishes the claim.

\setcounter{case}{0}

\begin{case}
 Assume there is a nonzero $w \in \Lie H$, such that $\phi_w = 0$ and $\yy_w
= 0$. 
 \end{case}
 Note, from the preceding paragraph, that $y_w = 0$. Then, because $\Lie Z =
0$, we must have $x_w \neq 0$. Therefore, \fullsref{HinN-square}{xvyu}
and~\fullref{HinN-linear}{phi=0&=0} hold, so $\mu(H) \approx
\muH{\|h\|}{\|h\|^2}$, so $H$ is a Cartan-decomposition subgroup.

\begin{case}
 Assume there does not exist such an element $w \in \Lie H$.
 \end{case}
 Then $H$ is described in~\fullref{HinN-nosquareproj}{onlyphiyyxx} and
in~\fullref{HinN-nolinearproj}{onlyphiyyxx}.
 \end{proof}

\section{When is the size of $\rho(h)$ linear?} \label{HinNlinear-section}

In this section, Proposition~\ref{HinN-linear} is a list of subgroups that
contain a sequence $\{h_m\}$ with $\rho(h_m) \asymp h_m$, and
Proposition~\ref{HinN-nolinearproj} is a list of subgroups that do not contain
such a sequence. Then Proposition~\ref{HinN-linearornot} shows that both lists
are complete.

\begin{Proposition} \label{HinN-linear}
 Assume that $G = \SU(2,n)$. Let $H$ be a closed,
connected subgroup of~$N$. There is a sequence $h_m \to \infty$ in~$H$ with
$\rho(h_m) \asymp h_m$ if either
 \begin{enumerate}
 \item \label{HinN-linear-phi=0&eta2=}
 there is a nonzero element~$z$ of~$\Lie Z$ with $|\eta_z|^2 = \xx_z \yy_z$;
or
 \item \label{HinN-linear-phi=0&=0}
there is an element~$u$ of~$\Lie h$, such that $\phi_u = 0$, $\dim_{\complex}
\langle x,y \rangle = 1$, and 
 $$ \xx_u |y_u|^2 + \yy_u |x_u|^2 + 2 \Im ( x_u
y_u^{\dagger} \cjg{\eta_u}) = 0 ;$$
 or
 \item \label{HinN-linear-y=0&yy=0}
there is an element~$h$ of~$H$
with $y_h = 0$, $\yy_h = 0$ and
 $|x_h|^2 + 2 \Re(\phi_h \cjg{\eta_h}) \neq 0$; or
 \item \label{HinN-linear-y=0&etaneq0}
there are elements~$u$ of~$\Lie h$ and~$z$ of~$\Lie z$,
such that $\phi_u \neq 0$, $y_u = 0$, $\yy_u \neq 0$, $\eta_z \neq 0$, and
$\yy_z = 0$; or
 \item \label{HinN-linear-phi&y}
there are nonzero elements~$u$ of~$\Lie h$ and $z$ of~$\Lie Z$,
such that
 $\phi_u \neq 0$, $y_u \neq 0$, $\yy_z = 0$, $\phi_u \cjg{\eta_z}$ is
real, and
 $$\xx_z |y_u|^2 - \phi_u \yy_u \cjg{\eta_z} + 2 \Im \bigl(
\cjg{\eta_z} x_u y_u^{\dagger} \bigr) = 0 .$$
 \end{enumerate}
 \end{Proposition}

\begin{proof}
  We separately consider each of the five cases in the statement of the
proposition.

 \pref{HinN-linear-phi=0&eta2=}
 From \fullref{HinN-nosquareproj}{H=Z}, we have $\rho(h) \asymp h$ for all $h
\in \exp(\real z)$.

\pref{HinN-linear-phi=0&=0}
 Replacing $H$ by a conjugate under~$\langle U_\alpha, U_{-\alpha} \rangle$, we
may assume that $y_u = 0$ (and $x_u \neq 0$). Then, from the assumption of this
case, we know that $\yy_u$ is also~$0$. Therefore,
\fullref{HinN-nosquareproj}{y&yy=0} implies that $\rho(h) \asymp h$ for all $h
\in \exp(\real u)$.

 \pref{HinN-linear-y=0&yy=0}
 From \fullref{HinN-nosquareproj}{y&yy=0}, we have $\rho(h) \asymp h$ for all $h
\in \exp(\real u)$.

 \pref{HinN-linear-y=0&etaneq0}.
 For any large~$t$, choose $h \in \exp( tu + \Lie z )$, such that 
 $\xx_h \yy_h + \tfrac{1}{12} |\phi_h|^2 \yy_h^2 - |\eta_h|^2 = 0$. Note that
$\eta_h \asymp |\phi_h y_h| \asymp t^2$, so $h \asymp \Re h_{1,n+2} \asymp t^3$,
but $h_{i,j} = O(t^2)$ whenever $(i,j) \neq (1,n+2)$, and $h_{i,j} = O(t)$
whenever $i \notin \{1,2\}$ or $j \notin \{n+1,n+2\}$. From the choice
of~$h$, we have 
 $$ \Delta(h)
 = 0 + i \left( \tfrac{1}{2} |x_h|^2 \yy_h \right) = O(t^3) = O(h) ,
 $$
 so it is not difficult to see that $\rho(h) \asymp h$.

\pref{HinN-linear-phi&y}
 Replacing $\Lie H$ by a conjugate, we may assume $u \in \Lie U_{\alpha} +
\Lie U_{\beta}$. (First, conjugate by an element of~$U_\beta$ to make $\yy_u =
0$. Then conjugate by an element of~$U_{\alpha}$ to make $x_u$ orthogonal
to~$y_u$. Then conjugate by an element of~$U_\beta$ that centralizes~$y_u$, to
make $x_u = 0$. Then conjugate by an element of~$U_{\alpha+\beta}$ to make
$\eta_u = 0$. Then conjugate by an element of~$U_{\alpha+2\beta}$ to
make~$\xx_u = 0$.) Then, by assumption, we must have $\xx_z = 0$, because
$\yy_u = 0$ and $x_u = 0$.

Furthermore, replacing $\Lie H$ by a conjugate under a diagonal matrix (that
belongs to~$G$), we may assume that $\phi_u$ and~$y_u$ are real. Then
$\eta_z$ must also be real (because $\phi_u \cjg{\eta_z}$ is real).
Thus, we see that $u,z \in \so(2,n)$. So \cite[Thm.~5.3(1)]{OhWitte-CDS}
implies that $H$ is a Cartan-decomposition subgroup.
 \end{proof}

\begin{Proposition} \label{HinN-nolinearproj}
 Assume that $G = \SU(2,n)$. Let $H$ be a closed,
connected, nontrivial subgroup of~$N$ such that
 \begin{equation} \label{HinN-nolinearproj-Z}
 \text{$|\eta_z|^2 \neq \xx_z \yy_z$,
for every nonzero $z \in \Lie Z$.}
 \end{equation}
 \begin{enumerate}
 \item \label{HinN-nolinearproj-H=Z}
 If $\Lie H = \Lie Z$ {\upshape(}so $\dim H \le 3${\upshape)}, then $\rho(h)
\asymp \|h\|^2$ for every $h \in H$.
 \item \label{HinN-nolinearproj-phi=0(dim2)}
 If $\phi_h = 0$ and $\dim_{\complex} \langle x_u, y_u \rangle \neq 1$ for
every $h \in \Lie H$, then $\rho(h) \asymp \|h\|^2$ for every $h \in H$.
 \item \label{HinN-nolinearproj-phi=0(dim1&2)}
 If $\phi_h = 0$ for every $h \in \Lie H$, there exist nonzero $u$
and~$v$ in~$\Lie H$, such that  $\dim_{\complex} \langle x_u, y_u
\rangle \neq 1$ and  $\dim_{\complex} \langle x_v, y_v
\rangle = 1$, and
 $\xx_v |y_v|^2 + \yy_v |x_v|^2 + 2 \Im(x_v y_v^{\dagger}
\cjg{\eta_v}) \neq 0$
 for every such $v \in \Lie H$, then $\mu(H) \approx
\muH{\|h\|^{3/2}}{\|h\|^2}$.
 \item \label{HinN-nolinearproj-onlyphiyyxx}
 If $\dim \Lie H \le 3$, $\Lie Z = 0$, we have $\phi_v \asymp y_v$ and $v = O
\bigl( |\phi_v| + |\yy_v| \bigr)$ for every $v \in \Lie H$, and there exists $u
\in \Lie H$, such that $\phi_u \neq 0$, then $\rho(h) \asymp \|h\|^{3/2}$ for
every $h \in H$.
 \item \label{HinN-nolinearproj-alpha}
 If $\dim \Lie H \le 2$ and $\phi_h \neq 0$, $y_h = 0$, $\yy_h = 0$,
and $|x_h|^2 + 2 \Re(\phi_h \cjg{\eta_h}) = 0$ for every nonzero $h \in
\Lie H$, then $\rho(h) \asymp \|h\|^2$ for every $h \in H$.
 \item \label{HinN-nolinearproj-u+z}
 If $\dim \Lie H = 2$ and there exist nonzero $u \in \Lie H$ and $z \in \Lie
Z$, such that $\phi_u \neq 0$, $y_u \neq 0$, $\yy_z = 0$, $\phi_u
\cjg{\eta_z}$ is real, and 
 $\xx_z |y_u|^2 - \phi_u \yy_u \cjg{\eta_z} + 2 \Im \bigl(
\cjg{\eta_z} x_u y_u^{\dagger} \bigr) \neq 0$, then $\mu(H) \approx
\muH{\|h\|^{5/4}}{\|h\|^2}$.
 \item \label{HinN-nolinearproj-phi=0(dim1)}
 If $\dim \Lie H = 1$, and we have $\phi_h = 0$,  $\dim_{\complex} \langle
x_h, y_h \rangle = 1$, and
 $$\xx_h |y_h|^2 + \yy_h |x_h|^2 + 2 \Im(x_h y_h^{\dagger}
\cjg{\eta_h}) \neq 0$$
 for every nonzero $h \in \Lie H$, then $\rho(h) \asymp \|h\|^{3/2}$ for
every $h \in H$.
 \item \label{HinN-nolinearproj-dim1(h4/3)}
 If $\dim \Lie H = 1$, and $\phi_h \neq 0$, $y_h = 0$, and $\yy_h \neq 0$,
for every nonzero $h \in \Lie H$, then $\rho(h) \asymp \|h\|^{4/3}$ for
every $h \in H$.
 \end{enumerate}
 \end{Proposition}

\begin{proof}
  We separately consider each of the eight cases in the statement of the
proposition.

 \pref{HinN-nolinearproj-H=Z}
 From~\eqref{HinN-nolinearproj-Z}, we know that the quadratic form $|\eta|^2
- \xx \yy$ is anisotropic on~$\Lie Z = \Lie H$, so 
 $$\Delta(h) = |\eta_h|^2 -
\xx_h \yy_h
 \asymp |\eta_h|^2 + \xx_h^2 + \yy_h^2 
 \asymp \|h\|^2 .$$

\pref{HinN-nolinearproj-phi=0(dim2)}
 Because $\dim_{\complex} \langle x_u, y_u \rangle \neq 1$, we have 
 $$ |x_h|^2 |y_h|^2 - |x_h y_h^{\dagger}|^2 \asymp |x_h|^4 + |y_h|^4 ,$$
 so Lemma~\ref{HinN-phi=0(dim<xy>)} implies $\rho(h) \asymp \|h\|^2$.

\pref{HinN-nolinearproj-phi=0(dim1)}
 From either Proposition~\fullref{HinN-nosquareproj}{phi&y=0}
or~\fullref{HinN-nosquareproj}{phi=0-neq0} (depending on whether $y_h$
is~$0$ or not), we have $\rho(h) \asymp \|h\|^{3/2}$ for every $h \in H$.

\pref{HinN-nolinearproj-phi=0(dim1&2)}
 From Lemma~\ref{HinN-phi=0(dim<xy>)}, we have $\|h\|^{3/2} = O \bigl(
\rho(h) \bigr)$.

From~\pref{HinN-nolinearproj-phi=0(dim2)}, we see that $\rho(h) \asymp
\|h\|^2$ for $h \in \exp( \real u )$.

From~\pref{HinN-nolinearproj-phi=0(dim1)}, we see that $\rho(h) \asymp
\|h\|^{3/2}$ for $h \in \exp( \real v )$.

\pref{HinN-nolinearproj-onlyphiyyxx}
 See Proposition~\fullref{HinN-nosquareproj}{onlyphiyyxx}.

\pref{HinN-nolinearproj-alpha}
 Because $\Re h_{1,n+2} = 0$, it is easy to see that
 $\rho(h) \asymp \phi_h^2 \asymp \|h\|^2$.

\pref{HinN-nolinearproj-u+z}
 Replacing $H$ by a conjugate, we may assume $x_u = 0$ and $\yy_u = 0$.
Therefore, $x_h = 0$ and $\yy_h = 0$ for every $h \in H$. Thus 
 $$ \\x_z |y_u|^2 = \xx_z |y_u|^2 - \phi_u \yy_u \cjg{\eta_z} + 2 \Im
\bigl( \cjg{\eta_z} x_u y_u^{\dagger} \bigr) \neq 0 ,$$
 so $\xx_z \neq 0$. From~\eqref{HinN-nolinearproj-Z}, we know $\eta_z \neq
0$.

 We have $\rho(t z) \asymp \|tz\|^2$ (see \fullref{HinN-nolinearproj}{H=Z}). 

 Because $\phi_u$ is a real multiple of $\eta_z$, we may let $h$ be a large
element of~$H$, such that $\eta_h = - |y_h|^2 \phi_h/12 + O(\phi_h)$. (So
$y_h \asymp \phi_h$ and $\xx_h \asymp \eta_h \asymp \phi_h^3$.) Then
 \begin{eqnarray*}
 \Delta(h) 
 &=& \left( -|\eta_h|^2 - \frac{1}{6} |y_h|^2 \eta_h \cjg{\phi_h} -
\frac{1}{144} |y_h|^4 |\phi_h|^2 \right)
 + i \left( \frac{1}{2} \xx_h |y_h|^2 \right) \\
 &=& O(\phi_h^4) +  i \left( \frac{1}{2} \xx_h |y_h|^2 \right) 
 \asymp \phi_h^5
 .
 \end{eqnarray*}
 It is clear that all other matrix entries of $\rho(h)$ are $O(\phi_h^5)$.
 Thus, we have $\rho(h) \asymp \phi_h^5 \asymp \|h\|^{5/4}$.

 Now suppose there is a sequence $h_m \to \infty$ in~$H$ with $\rho(h_m) = o
\bigl( \|h_m\|^{5/4} \bigr)$.

 \setcounter{case}{0}

\begin{case}
 Assume $\eta_m = o(\phi_m^3)$.
 \end{case}
 We have $h_m \asymp \phi_m^4$, so 
 $$ \phi_m^6 \asymp \Re \Delta(h_m) = O \bigl( \rho(h_m) \bigr) =
o(\|h_m\|^{5/4}) = o(\phi_m^5) .$$
 This is a contradiction.

\begin{case}
 Assume $\phi_m^3 = o(\eta_m)$.
 \end{case}
 We have $h_m \asymp \Re h_{1,n+2} \asymp \phi_m \eta_m$, so
 $$ \eta_m^2 \asymp \Re \Delta(h_m) = O \bigl( \rho(h_m) \bigr) =
o(\|h_m\|^{5/4}) = o(\|h_m\|^{3/2}) = o( |\phi_m \eta_m|^{3/2} )
 = o(\eta_m^2)
.$$
 This is a contradiction.

\begin{case}
 Assume $\eta_m \asymp \phi_m^3$.
 \end{case}
 We have $h_m = O(\phi_m^4)$, so
 $$ \phi_m^5 \asymp \xx_m |y_m|^2 \asymp \Im \Delta(h_m) = O \bigl(
\rho(h_m) \bigr) = o \bigl( \|h_m\|^{5/4} \bigr) = o(\phi_m^5)
.$$
 This is a contradiction.

\pref{HinN-nolinearproj-dim1(h4/3)}
 See Proposition~\fullref{HinN-nosquareproj}{phi=yy}.
 \end{proof}
 
\begin{Proposition} \label{HinN-linearornot}
 Assume that $G = \SU(2,n)$. Let $H$ be a closed, connected, nontrivial
subgroup of~$N$.
 \begin{enumerate}
 \item \label{HinN-linearornot-yes}
 There is a sequence $h_m \to \infty$ in~$H$
with $\rho(h_m) \asymp h_m$ if and only if $H$ is one of the subgroups
described in Proposition~\ref{HinN-linear}.
 \item \label{HinN-linearornot-no}
 There is \textbf{not} a sequence $h_m \to \infty$ in~$H$
with $\rho(h_m) \asymp \|h_m\|^2$ if and only if $H$ is one of the subgroups
described in Proposition~\ref{HinN-nolinearproj}.
 \end{enumerate}
 \end{Proposition}

\begin{proof}
 It suffices to show that $H$ is described in either
Proposition~\ref{HinN-linear} or Proposition~\ref{HinN-nolinearproj}.

 We may assume \pref{HinN-nolinearproj-Z} holds (otherwise,
Conclusion~\fullref{HinN-linear}{phi=0&eta2=} holds).

\setcounter{case}{0}

\begin{case}
 Assume $\phi_h = 0$ for every $h \in H$.
 \end{case}
 We may assume there exists $v \in \Lie H$, such that $\dim_{\complex}
\langle x_v, y_v \rangle = 1$ (otherwise
\fullref{HinN-nolinearproj}{phi=0(dim2)} holds). Furthermore, we may assume
$\xx_v |y_v|^2 + \yy_v |x_v|^2 + 2 \Im(x_v y_v^{\dagger} \cjg{\eta_v}) \neq 0$
 for every such~$v$ (otherwise \fullref{HinN-linear}{phi=0&=0} holds). Then
we may assume $\dim_{\complex} \langle x_u, y_u \rangle = 1$ for every
nonzero $u \in \Lie H$ (otherwise
\fullref{HinN-nolinearproj}{phi=0(dim1&2)} holds).

The argument in Subsubcase~\ref{HinN-squareornot-subsub<xy>not<y>} of the proof
of Proposition~\ref{HinN-nosquareproj} implies there exists $\lambda \in
\complex$, such that, for every $h \in H$, we have $x_h = \lambda y_h$
(or vice-versa: for every $h$, we have $y_h = \lambda x_h$). Thus, replacing
$H$ by a conjugate under $\langle U_\alpha, U_{-\alpha} \rangle$, we may
assume $x_h = 0$ for every $h \in H$. 

If $\dim H > 1$, then there is some nonzero $u \in \Lie H$, such that $\xx_h
= 0$. This contradicts the fact that $\xx_v |y_v|^2 + \yy_v |x_v|^2 + 2 \Im(x_v y_v^{\dagger}
\cjg{\eta_v}) \neq 0$. Thus, we conclude that $\dim H = 1$, so
\fullref{HinN-nolinearproj}{phi=0(dim1)} holds.

\begin{case}
 Assume the projection of~$\Lie H$ to~$\Lie U_\alpha$ is one-dimensional.
 \end{case}
 Replacing~$H$ by a conjugate under~$A$, we may assume $\phi_h$ is real for
every $h \in H$. Fix some $u \in \Lie H$, such that $\phi_u \neq 0$.

 We may assume that $\Lie U_{2\alpha+2\beta} \not\subset \Lie H$
(otherwise Conclusion~\fullref{HinN-linear}{phi=0&eta2=} holds).
Therefore $[\Lie H, u]$ must be zero, so $\yy_z = 0$ and $\eta_z$ is a
nonzero real, for every nonzero $z \in \Lie Z$. (This implies $\dim \Lie Z
\le 1$.)

\begin{subcase}
 Assume $y_h = 0$ for every $h \in H$.
 \end{subcase}
 We may assume Conclusion~\fullref{HinN-linear}{phi=0&=0} does not hold.

We claim that $\Lie H = \real u + \Lie Z$. Suppose not. Then there is some
$v \in \Lie H$, such that $\phi_v = 0$ and $x_v \neq 0$. Because
Conclusion~\fullref{HinN-linear}{phi=0&=0} does not hold, we must have $\yy_v
\neq 0$. Then $[v,u,u]$ is a nonzero element of~$\Lie U_{2\alpha+2\beta}$.
(This can be seen easily by replacing $H$ with a conjugate, so that $u \in
\Lie U_\alpha$.) This contradicts our assumption that $\Lie
U_{2\alpha+2\beta} \not\subset \Lie H$. 

If $\yy_u \neq 0$, then either
Conclusion~\fullref{HinN-nolinearproj}{dim1(h4/3)} or
\fullref{HinN-linear}{y=0&etaneq0} holds (depending on whether $\Lie Z$
is~$0$ or not). If $\yy_u \neq 0$, then \fullref{HinN-linear}{y=0&yy=0} or
\fullref{HinN-nolinearproj}{u+z} holds.

\begin{subcase}
 Assume the projection of~$\Lie H$ to $\Lie U_{\beta}$ is nontrivial.
 \end{subcase}
 Then we may assume $y_u \neq 0$.

\begin{subsubcase}
 Assume there are nonzero $v \in \Lie H$ and $z \in \Lie Z$, such that
$\phi_v = 0$, $y_v = 0$, and $x_v \neq 0$.
 \end{subsubcase}
 We may assume that Conclusion~\fullref{HinN-linear}{phi&y} does not hold.
Therefore, for every real~$t$, we must have
 \begin{eqnarray*}
 0 &\neq& \xx_z |y_u|^2 - \phi_u (\yy_u + t \yy_v) \cjg{\eta_z} + 2
\Im \bigl( \cjg{\eta_z} (x_u + t x_v) y_u^{\dagger} \bigr) \\
 &=& t \bigl[ - \phi_u \yy_v \cjg{\eta_z} + 2
\Im \bigl( \cjg{\eta_z} x_v y_u^{\dagger} \bigr) \bigr]
 + \text{constant}
 .
 \end{eqnarray*}
 Thus, the coefficient of~$t$ must vanish, which (using the fact that
$\eta_z$ is real and nonzero) means
 \begin{equation} \label{coeff(t)vanish}
 0 = - \phi_u \yy_v + 2 \Im \bigl( x_v y_u^{\dagger} \bigr) .
 \end{equation}

We have $[u,v] \in \Lie Z$, so $\eta_{[u,v]}$ is real. Thus,
 $$ 0 = \Im \eta_{[u,v]} = \Im \bigl( x_v y_u^{\dagger} + i \phi_u \yy_v
\bigr)
 = \Im \bigl( x_v y_u^{\dagger} \bigr) + \phi_u \yy_v .$$
 Comparing this with \eqref{coeff(t)vanish}, we conclude that $\phi_u \yy_v
= 0$. 
 Therefore $\yy_v = 0$, so Conclusion~\fullref{HinN-linear}{phi=0&=0} holds
(for the element~$v$).

\begin{subsubcase}
 Assume there do not exist nonzero $v \in \Lie H$ and $z \in \Lie Z$, such
that $\phi_v = 0$, $y_v = 0$, and $x_v \neq 0$.
 \end{subsubcase}
 We must have 
 \begin{equation} \label{phi=0->y=0}
 \text{$y_w = 0$ for every $w \in \Lie H$, such that $\phi_w = 0$.}
 \end{equation}
 (Otherwise, we obtain a contradiction by setting $v = [u,w]$ and $z =
[u,w,w]$.)
 We may assume 
 \begin{equation} \label{yy=0->y=0}
 \text{$\yy_v \neq 0$ for every $v \in \Lie H$ such that $\phi_v =
0$, $y_v = 0$, and $x_v \neq 0$.}
 \end{equation}
 (Otherwise, Conclusion~\fullref{HinN-linear}{phi=0&=0} holds.) 

We claim $\dim \Lie H \le 2$. If not, then there exist linearly
independent $v,w \in \Lie H$, such that $\phi_v = \phi_w = 0$.
From~\pref{phi=0->y=0}, we know that $y_v = y_w = 0$. By replacing with a
linear combination, we may assume $\yy_w = 0$. Then, from~\pref{yy=0->y=0},
we know that $x_w = 0$, so $w \in \Lie Z$. Because $\Lie Z$ is (at most)
one-dimensional, but $v$ and~$w$ are linearly independent, we know that $v
\notin \Lie Z$, so $x_v \neq 0$. This contradicts the assumption of this
subsubcase.

We may now assume $\dim \Lie H = 2$ (otherwise
Conclusion~\fullref{HinN-nolinearproj}{alpha} holds). Choose a nonzero
$v \in \Lie H$, such that $\phi_v = 0$. If $x_v \neq 0$, then
Conclusion~\fullref{HinN-nolinearproj}{alpha} holds. If $x_v = 0$, then $v
\in \Lie Z$, so either Conclusion~\fullref{HinN-linear}{phi&y}
or~\fullref{HinN-nolinearproj}{u+z} holds.

\begin{case}
 Assume the projection of~$\Lie H$ to~$\Lie U_\alpha$ is two-dimensional.
 \end{case}
 We may assume $\Lie Z = 0$ (otherwise, $\Lie U_{2\alpha+2\beta} \subset
\Lie H$, so Conclusion~\fullref{HinN-linear}{phi=0&eta2=} holds). We may
assume $y_h = 0$ for every $h \in H$ (otherwise Lemma~\ref{HinN-Z=0} implies
that either \fullref{HinN-nolinearproj}{alpha}
or~\fullref{HinN-linear}{phi=0&=0} applies. Therefore $[\Lie H, \Lie H]
\subset \Lie Z = 0$, so $\Lie H$ is abelian.

Let $u,v \in \Lie H$ with $\phi_u = 1$ and $\phi_v = i$. Then
 $$ 0 = \eta_{[u,v]} = i \yy_v + \yy_u ,$$ 
 so $\yy_u = \yy_v = 0$. Then, for every $w \in \Lie H$, we have 
 $0 = \eta_{[u,w]} = i \yy_w$, so $\yy_w = 0$.
 We may assume 
 \begin{equation} \label{x^2+2phieta=0}
 |x_h|^2 + 2 \Re(\phi_h \cjg{\eta_h}) = 0
 \end{equation}
 for every $h \in
\Lie h$ (otherwise Conclusion~\fullref{HinN-linear}{y=0&yy=0} holds). This
implies $\dim \Lie H = 2$ (otherwise, there is some $w \in \Lie H$ such that
$\phi_w = 0$ and $x_w \neq 0$, and then \eqref{x^2+2phieta=0} does not hold
for $h = u + tw$ when $t$ is sufficiently large).
 Thus, Conclusion~\fullref{HinN-nolinearproj}{alpha} holds.
 \end{proof}

\begin{Lemma} \label{HinN-phi=0(dim<xy>)}
 Let $H$ be a closed, connected, nontrivial subgroup of~$N$.
  Assume $\phi_h = 0$ for every $h \in \Lie H$, that
\pref{HinN-nolinearproj-Z} holds, and that
 $\xx_v |y_v|^2 + \yy_v |x_v|^2 + 2 \Im(x_v y_v^{\dagger}
\cjg{\eta_v}) \neq 0$
 for every $v \in \Lie H$ such that $\dim_{\complex} \langle x_v, y_v
\rangle = 1$.
 Then $\|h\|^{3/2} = O \bigl( \Delta(h) \bigr)$ for every $h \in H$.

Furthermore, $\Delta(h) \asymp \|h\|^2$ whenever $|x_h|^2 |y_h|^2 - |x_h
y_h^{\dagger}|^2 \asymp |x_h|^4 + |y_h|^4$.
 \end{Lemma}

\begin{proof}
 We have $h \asymp |x_h|^2 + |y_h|^2 + |\xx_h| + |\yy_h| + |\eta_h|$. 
 Also, from \eqref{HinN-nolinearproj-Z}, we have $|\eta_z|^2 - \xx_z \yy_z
\asymp \bigl( |\xx_z| + |\yy_z| + |\eta_z| \bigr)^2$ for every $z \in \Lie
Z$. Also, 
 $\xx_v |y_v|^2 + \yy_v |x_v|^2 + 2 \Im(x_v y_v^{\dagger}
\cjg{\eta_v}) \asymp |v|^3$ whenever $\dim_{\complex} \langle x_v,y_v
\rangle = 1$.

\setcounter{case}{0}

\begin{case}
 Assume $|x_h|^2 |y_h|^2 - |x_h y_h^{\dagger}|^2 = o \bigl( |x_h|^4 + |y_h|^4
\bigr)$.
 \end{case}
 Then there is some $v \in \Lie H$ such that $v - \log h = o \bigl( |x_h| +
|y_h| \bigr)$ and $|x_v|^2 |y_v|^2 - |x_v y_v^{\dagger}|^2 = 0$. We have
$\dim_{\complex} \langle x_v,y_v \rangle = 1$. Therefore
 \begin{eqnarray*}
 \Im \Delta(h)
 &\asymp& \xx_h |y_h|^2 + \yy_h |x_h|^2 + 2 \Im(x_h y_h^{\dagger}
\cjg{\eta_h}) \\
 &=& \xx_v |y_v|^2 + \yy_v |x_v|^2 + 2 \Im(x_v y_v^{\dagger} \cjg{\eta_v}) \\
    && {\ } + o \bigl( |\eta_h|^3 + |\xx_h|^3 + |\yy_h|^3 + |x_h|^3 + |y_h|^3
\bigr) \\
 &\asymp& |v|^3 + o \bigl( |\eta_h|^3 + |\xx_h|^3 + |\yy_h|^3 + |x_h|^3 +
|y_h|^3 \bigr) \\
 &\asymp& |\eta_h|^3 + |\xx_h|^3 + |\yy_h|^3 + |x_v|^3 + |y_v|^3 \\
 &\neq& o \bigl( \|h\|^{3/2} \bigr) .
 \end{eqnarray*}
 Thus, $\|h\|^{3/2} = O \bigl( \rho(h) \bigr)$.

\begin{case}
 Assume $|x_h|^2 |y_h|^2 - |x_h y_h^{\dagger}|^2 \asymp |x_h|^4 + |y_h|^4$.
 \end{case}
 We may assume $\Re \Delta(h) = o \bigl( |x_h|^4 + |y_h|^4 \bigr)$ for
otherwise it is clear that $\Re \Delta(h) \asymp \|h\|^2$. (So we have 
 $\|h\| \asymp |\eta_h| + |\xx_h| + |\yy_h| \asymp |x_h|^2 + |y_h|^2$.)
 Thus, there is some $z \in \Lie Z$, such that $z - \log h = o \bigl( \log h
\bigr)$ and 
 $$ |\eta_z|^2 - \xx_z \yy_z 
 = - \frac{1}{4} \bigl( |x_h|^2 |y_h|^2 - |x_h y_h^{\dagger}|^2 \bigr)
 + o
\bigl( |x_h|^4 + |y_h|^4 \bigr)
 < 0 .$$
 (This implies that $\xx_z$ and~$\yy_z$ must have the same sign.)
 From~\pref{HinN-nolinearproj-Z}, we conclude that $|\eta_z|^2 - \xx_z \yy_z
< 0$ for every $z \in \Lie Z$. Thus, there is a constant $\epsilon < 1$, such
that $|\eta_z| \le \epsilon \sqrt{\xx_z \yy_z}$ for every $z \in
\Lie Z$. Then
  $$|\Im(x_h y_h^{\dagger} \cjg{\eta_z})|
 \le |\eta_z| |x_h| |y_h|
 \le \frac{\epsilon}{2} \bigl|
\xx_z |y_h|^2 + \yy_z |x_h|^2 \bigr| ,$$
 so
 $$ \Im(x_h y_h^{\dagger} \cjg{\eta_z}) + \frac{1}{2} 
\xx_z |y_h|^2 + \frac{1}{2} \yy_z |x_h|^2 
 \asymp \frac{1}{2}  \xx_z
|y_h|^2 + \frac{1}{2} \yy_z |x_h|^2 . $$
 Therefore
 \begin{eqnarray*}
 \Im \Delta(h) 
 &=& \Im(x_h y_h^{\dagger} \cjg{\eta_h}) + \frac{1}{2}  \xx_h
|y_h|^2 + \frac{1}{2} \yy_h |x_h|^2 \\
 &=& \Im(x_h y_h^{\dagger} \cjg{\eta_z}) + \frac{1}{2}  \xx_z
|y_h|^2 + \frac{1}{2} \yy_z |x_h|^2
 + o \bigl( ( |x_h|^2 + |y_h|^2 ) \log h \bigr) \\
 &\asymp& \frac{1}{2}  \xx_z
|y_h|^2 + \frac{1}{2} \yy_z |x_h|^2 \\
 &\asymp& |x_h|^4 + |y_h|^4 \\
 &\asymp& \|h\|^2 . 
 \end{eqnarray*}
 \end{proof}

\section{Non-Cartan-decomposition subgroups contained in~$N$}
\label{CDSinN-section}

\begin{Theorem} \label{HinN-notCDS}
 Assume that $G = \SU(2,n)$. Here is a complete list of the closed,
connected, nontrivial subgroups~$H$ of~$N$, such that $H$ is {\bf not} a
Cartan-decomposition subgroup.
 \begin{enumerate}
 \item \label{HinN-notCDS-H=Z}
 If $\dim \Lie H = 1$, $\Lie H = \Lie Z$, and we have $|\eta_h|^2 = \xx_h
\yy_h$ for every $h \in H$, then $\rho(h) \asymp h$ for every $h
\in H$.
 \item \label{HinN-notCDS-phi&y=0}
 If $\phi_h = 0$ and $y_h = 0$ for every $h \in \Lie H$, there is some $u \in
\Lie H$, such that $\yy_u \neq 0$, and $\Lie Z \subset \Lie U_{2\alpha +
2\beta}$, then $\mu(H) \approx \muH{\|h\|}{\|h\|^{3/2}}$, unless $\dim H =
1$, in which case $\rho(h) \asymp \|h\|^{3/2}$ for every $h \in H$.
 \item \label{HinN-notCDS-phi=0}
 Suppose $\phi_h = 0$ for every $h \in \Lie H$, and there is some $\lambda
\in \complex$, such that $x_h = \lambda y_h$ for every $h \in H$, and we
have $\eta_z = i\lambda \yy_z$ and $\xx_z = |\lambda|^2 \yy_z$ for every $z
\in \Lie Z$.
 \begin{enumerate}
 \item \label{HinN-notCDS-phi=0-neq0}
 If there is some $u \in \Lie H$, such that $\xx_u + |\lambda|^2 \yy_u
+ 2 \Im( \lambda \cjg{\eta_u}) \neq 0$, then $\mu(H) \approx
\muH{\|h\|}{\|h\|^{3/2}}$, unless $\dim H = 1$, in which case $\rho(h)
\asymp \|h\|^{3/2}$ for every $h \in H$.
 \item \label{HinN-notCDS-phi=0-=0}
 Otherwise, $\rho(h) \asymp h$ for every $h \in H$.
 \end{enumerate}
 \item \label{HinN-notCDS-y&yy=0}
 If $y_h = 0$, $\yy_h = 0$, and $|x_h|^2 + 2 \Re(\phi_h \cjg{\eta_h})
\neq 0$ for every $h \in \Lie H \setminus \Lie U_{2\alpha+2\beta}$ {\upshape(}so
$\Lie Z \subset \Lie U_{2\alpha+2\beta}${\upshape)}, then $\rho(h) \asymp h$ for
every $h \in H$.
 \item \label{HinN-notCDS-phi=yy}
 If $\Lie Z = 0$, there is some $u \in \Lie H$ and some nonzero $\phi_0 \in
\complex$, such that $\phi_u \neq 0$, and we have $\phi_h = \phi_0 \yy_h$ and
$y_h = 0$, for every $h \in \Lie H$, then $\mu(H) \approx
\muH{\|h\|}{\|h\|^{4/3}}$, unless $\dim H = 1$, in which case, $\rho(h) \asymp
\|h\|^{4/3}$ for every $h \in H$.
 \item \label{HinN-notCDS-phi=0(dim2)}
 If $\phi_h = 0$ and $\dim_{\complex} \langle x_u, y_u \rangle \neq 1$ for
every $h \in \Lie H$, and $|\eta_z|^2 \neq \xx_z \yy_z$,
for every nonzero $z \in \Lie Z$, then $\rho(h) \asymp \|h\|^2$ for every $h
\in H$.
 \item \label{HinN-notCDS-phi=0(dim1&2)}
 If $\phi_h = 0$ for every $h \in \Lie H$, there exist nonzero $u$
and~$v$ in~$\Lie H$, such that  $\dim_{\complex} \langle x_u, y_u
\rangle \neq 1$ and $\dim_{\complex} \langle x_v, y_v
\rangle = 1$, and we have
 $\xx_v |y_v|^2 + \yy_v |x_v|^2 + 2 \Im(x_v y_v^{\dagger}
\cjg{\eta_v}) \neq 0$
 for every such $v \in \Lie H$,
 and $|\eta_z|^2 \neq \xx_z \yy_z$,
for every nonzero $z \in \Lie Z$ then $\mu(H) \approx
\muH{\|h\|^{3/2}}{\|h\|^2}$.
 \item \label{HinN-notCDS-onlyphiyyxx}
 If $\dim \Lie H \le 3$, $\Lie Z = 0$, we have $\phi_v \asymp y_v$ and $v = O
\bigl( |\phi_v| + |\yy_v| \bigr)$ for every $v \in \Lie H$, and there exists $u
\in \Lie H$, such that $\phi_u \neq 0$, then $\rho(h) \asymp \|h\|^{3/2}$ for
every $h \in H$.
 \item \label{HinN-notCDS-dim3}
 If $\dim \Lie H = 2$, $\Lie Z = \Lie U_{2\alpha+2\beta}$, $\phi_h
\neq 0$ and $y_h \neq 0$ for every $h \in \Lie H \setminus \Lie Z$, then
$\mu(H) \approx \muH{\|h\|}{\|h\|^{3/2}}$.
 \item \label{HinN-notCDS-alpha}
 If $\dim \Lie H \le 2$ and $\phi_h \neq 0$, $y_h = 0$, $\yy_h = 0$,
and $|x_h|^2 + 2 \Re(\phi_h \cjg{\eta_h}) = 0$ for every nonzero $h \in
\Lie H$, then $\rho(h) \asymp \|h\|^2$ for every $h \in H$.
 \item \label{HinN-notCDS-u+z}
 If $\dim \Lie H = 2$ and there exist nonzero $u \in \Lie H$ and $z \in \Lie
Z$, such that $\phi_u \neq 0$, $y_u \neq 0$, $\yy_z = 0$, $\phi_u
\cjg{\eta_z} \neq 0$ is real, and 
 $\xx_z |y_u|^2 - \phi_u \yy_u \cjg{\eta_z} + 2 \Im \bigl(
\cjg{\eta_z} x_u y_u^{\dagger} \bigr) \neq 0$, then $\mu(H) \approx
\muH{\|h\|^{5/4}}{\|h\|^2}$.
 \end{enumerate}
 \end{Theorem}

\begin{proof}
 The theorem is obtained by merging the statement of
Proposition~\ref{HinN-nosquareproj} with the statement of
Proposition~\ref{HinN-nolinearproj}, and eliminating some redundancy
\see{CDS<>h_m}. Specifically:

\begin{itemize}
 \item \fullref{HinN-nosquareproj}{H=Z}  appears here as
\fullref{HinN-notCDS}{H=Z}.
 \item \fullref{HinN-nosquareproj}{phi&y=0} appears here as
\fullref{HinN-notCDS}{phi&y=0}.
 \item \fullref{HinN-nosquareproj}{phi=0} appears here as
\fullref{HinN-notCDS}{phi=0}.
 \item \fullref{HinN-nosquareproj}{y&yy=0} appears here as
\fullref{HinN-notCDS}{y&yy=0}.
 \item \fullref{HinN-nosquareproj}{phi=yy} appears here as
\fullref{HinN-notCDS}{phi=yy}.
 \item \fullref{HinN-nosquareproj}{onlyphiyyxx} appears here as
\fullref{HinN-notCDS}{onlyphiyyxx}.
 \item \fullref{HinN-nosquareproj}{dim3} appears here as
\fullref{HinN-notCDS}{dim3}.

 \item \fullref{HinN-nolinearproj}{H=Z} is a special case
of~\fullref{HinN-notCDS}{phi=0(dim2)}.
 \item \fullref{HinN-nolinearproj}{phi=0(dim2)} appears here as
\fullref{HinN-notCDS}{phi=0(dim2)}.
 \item \fullref{HinN-nolinearproj}{phi=0(dim1&2)} appears here as
\fullref{HinN-notCDS}{phi=0(dim1&2)}.
 \item \fullref{HinN-nolinearproj}{onlyphiyyxx} appears here as
\fullref{HinN-notCDS}{onlyphiyyxx}.
 \item \fullref{HinN-nolinearproj}{alpha} appears here as
\fullref{HinN-notCDS}{alpha}.
 \item \fullref{HinN-nolinearproj}{u+z} appears here as
\fullref{HinN-notCDS}{u+z}.
 \item \fullref{HinN-nolinearproj}{phi=0(dim1)} is a special case of
\fullref{HinN-notCDS}{phi=0-neq0} (with $\dim H = 1$).
 \item \fullref{HinN-nolinearproj}{dim1(h4/3)} is a special case of
\fullref{HinN-notCDS}{phi=yy} (with $\dim H = 1$).
 \end{itemize}
 \end{proof}

\begin{Corollary} \label{HinN-N_A(H)}
 Assume that $G = \SU(2,n)$. Here is a complete list of the closed,
connected, nontrivial subgroups~$H$ of~$N$, such that $H$ is {\bf not} a
Cartan-decomposition subgroup, {\bf and} $N_A(H)$ is nontrivial.
 \begin{enumerate}

 \item \label{HinN-N_A(H)-H=Z}
 Suppose $\dim \Lie H = 1$, $\Lie H = \Lie Z$, and we have $|\eta_h|^2 = \xx_h
\yy_h$ for every $h \in H$.
 \begin{enumerate}
 \item \label{HinN-N_A(H)-H=Z-A}
 If $\Lie H = \Lie U_{2\beta}$ or $\Lie H = \Lie U_{2\alpha+2\beta}$, then
$N_A(H) = A$.
 \item \label{HinN-N_A(H)-H=Z-a}
 Otherwise, $N_A(H) = \ker(\alpha)$.
 \end{enumerate}

 \item \label{HinN-N_A(H)-phi&y=0}
 Suppose $\phi_h = 0$ and $y_h = 0$ for every $h \in \Lie H$, there is some $u \in
\Lie H$, such that $\yy_u \neq 0$, and $\Lie Z \subset \Lie U_{2\alpha +
2\beta}$.
 If $\Lie H = \bigl( \Lie H \cap (\Lie U_{\alpha+\beta} + \Lie U_{2\beta})
\bigr) + \Lie Z$, then $N_A(H) = \ker(\alpha-\beta)$.

 \item \label{HinN-N_A(H)-phi=0}
 Suppose $\phi_h = 0$ for every $h \in \Lie H$, and there is some
nonzero $\lambda \in \complex$, such that $x_h = \lambda y_h$ for every $h
\in H$, and we have $\eta_z = i\lambda \yy_z$ and $\xx_z = |\lambda|^2 \yy_z$
for every $z \in \Lie Z$.
 If $\Lie H = \bigl( \Lie H \cap (\Lie U_{\beta} + \Lie U_{\alpha+\beta})
\bigr) + \Lie Z \neq \Lie Z$, then $N_A(H) = \ker(\alpha)$.

 \item \label{HinN-N_A(H)-phi=0&x=0}
 Suppose $\phi_h = 0$ and $x_h = 0$ for every $h \in \Lie H$, we have
$\Lie Z \subset \Lie U_{2 \beta}$, and $\Lie H \neq \Lie
Z$.
 \begin{enumerate}
 \item \label{HinN-N_A(H)-phi=0&x=0-A}
 If $\Lie H = ( \Lie H \cap \Lie U_{\beta}) + \Lie Z$, then $N_A(H) =
A$.
 \item Otherwise:
 \begin{enumerate}
 \item \label{HinN-N_A(H)-phi=0&x=0-notA-a+b}
 If $\Lie H = \bigl( \Lie H \cap
(\Lie U_{\beta} + \Lie U_{\alpha+2\beta}) \bigr) + \Lie Z$, then $N_A(H) =
\ker(\alpha+\beta)$.
 \item \label{HinN-N_A(H)-phi=0&x=0-notA-2a+b}
  If $\Lie H = \bigl( \Lie H \cap (\Lie U_{\beta} + \Lie
U_{2\alpha+2\beta}) \bigr) + \Lie Z$, then $N_A(H) = \ker(2\alpha+\beta)$.
 \item \label{HinN-N_A(H)-phi=0&x=0-notA-b}
  If $\Lie Z = 0$ and $\Lie H \subset \Lie U_{\beta} + \Lie
U_{2\beta}$, then $N_A(H) = \ker(\beta)$.
 \end{enumerate}
 \end{enumerate}

 \item \label{HinN-N_A(H)-y&yy=0}
 Suppose $y_h = 0$, $\yy_h = 0$, and $|x_h|^2 + 2 \Re(\phi_h \cjg{\eta_h})
\neq 0$ for every $h \in \Lie H \setminus \Lie U_{2\alpha+2\beta}$.
 \begin{enumerate}
 \item \label{HinN-N_A(H)-y&yy=0-A}
  If $\Lie H = ( \Lie H \cap \Lie U_{\alpha+\beta}) + \Lie Z$, then
$N_A(H) = A$.
 \item \label{HinN-N_A(H)-y&yy=0-a+b}
  If $\Lie H \subset \Lie U_{\alpha+\beta} + \Lie
U_{2\alpha+2\beta}$, but $\Lie H \neq ( \Lie H \cap \Lie
U_{\alpha+\beta}) + \Lie Z$, then $N_A(H) = \ker(\alpha+\beta)$.
 \item \label{HinN-N_A(H)-y&yy=0-b}
  If $\Lie H = \bigl( \Lie H \cap (\Lie U_{\alpha} + \Lie
U_{\alpha+\beta} + \Lie
U_{\alpha+2\beta}) \bigr) + \Lie Z$, but $\Lie H \not\subset \Lie
U_{\alpha+\beta} + \Lie U_{2\alpha+2\beta}$, then $N_A(H) = \ker(\beta)$.
 \end{enumerate}

 \item \label{HinN-N_A(H)-phi=yy}
 Suppose $\Lie Z = 0$, there is some nonzero $\phi_0 \in \complex$, such that
$\phi_h = \phi_0 \yy_h$ and $y_h = 0$, for every $h \in \Lie H$, and  there
is some $u \in \Lie H$, such that $\phi_u \neq 0$.
 If $\Lie H = \bigl( \Lie H \cap ( \Lie U_\alpha + \Lie
U_{2\beta}) \bigr) + (\Lie H \cap \Lie U_{\alpha+\beta})$, then $N_A(H) =
\ker(\alpha-2\beta)$.

 \item \label{HinN-N_A(H)-phi=0(dim2)}
 Suppose $\phi_h = 0$ and $\dim_{\complex} \langle x_u, y_u \rangle \neq 1$ for
every $h \in \Lie H$, and $|\eta_z|^2 \neq \xx_z \yy_z$,
for every nonzero $z \in \Lie Z$.
 \begin{enumerate}
 \item \label{HinN-N_A(H)-phi=0(dim2)-A}
 If $\Lie H \subset \Lie U_{\alpha+2\beta}$, then $N_A(H) = A$.
 \item \label{HinN-N_A(H)-phi=0(dim2)-a}
 If $\Lie H \not\subset \Lie U_{\alpha+2\beta}$, and $\Lie
H = \bigl( \Lie H \cap ( \Lie U_{\beta} + \Lie U_{\alpha+\beta}) \bigr) + \Lie
Z$, then $N_A(H) = \ker(\alpha)$.
 \end{enumerate}

 \item \label{HinN-N_A(H)-phi=0(dim1&2)}
 Suppose $\phi_h = 0$ for every $h \in \Lie H$, there exist nonzero $u,v
\in \Lie H$, such that  $\dim_{\complex} \langle x_u, y_u \rangle \neq 1$
and  $\dim_{\complex} \langle x_v, y_v \rangle = 1$, 
 $\xx_v |y_v|^2 + \yy_v |x_v|^2 + 2 \Im(x_v y_v^{\dagger}
\cjg{\eta_v}) \neq 0$
 for every such $v \in \Lie H$,
 and $|\eta_z|^2 \neq \xx_z \yy_z$,
for every nonzero $z \in \Lie Z$.
 \begin{enumerate}
 \item \label{HinN-N_A(H)-phi=0(dim1&2)-x}
 If $\Lie H = \bigl( \Lie H \cap (\Lie
U_{\alpha+\beta} + \Lie U_{2\beta}) \bigr) + (\Lie H \cap \Lie
U_{\alpha+2\beta})$, then $N_A(H) = \ker(\alpha-\beta)$ (and $\dim H \le 3$).
 \item \label{HinN-N_A(H)-phi=0(dim1&2)-y}
 If $\Lie H = \bigl( \Lie H \cap (\Lie
U_{\beta} + \Lie U_{2\alpha+2\beta}) \bigr) + (\Lie H \cap \Lie
U_{\alpha+2\beta})$, then $N_A(H) = \ker(2\alpha+\beta)$ (and $\dim H \le
3$).
 \end{enumerate}

 \item \label{HinN-N_A(H)-onlyphiyyxx}
 Suppose $\dim \Lie H \le 3$, 
 $\Lie H = \bigl( \Lie H \cap ( \Lie U_\alpha + \Lie U_{\beta})
\bigr) +  \bigl( \Lie H \cap ( \Lie U_{\alpha+\beta} + \Lie U_{2\beta})
\bigr)$,
 $\Lie H \cap ( \Lie U_\alpha + \Lie U_{\beta}) \neq 0$, and
 we have $\phi_h \asymp y_h$ and $x_h \asymp \yy_h$ for
$h \in \Lie H$, then $N_A(H) = \ker(\alpha-\beta)$.

 \item \label{HinN-N_A(H)-dim3}
 Suppose $\dim \Lie H = 2$, $\Lie Z = \Lie U_{2\alpha+2\beta}$,
$\phi_h \neq 0$ and $y_h \neq 0$ for every $h \in \Lie H \setminus \Lie Z$.
 If $\Lie H = \bigl( \Lie H \cap ( \Lie U_\alpha + \Lie U_{\beta})
\bigr) + \Lie Z$, then $N_A(H) = \ker(\alpha-\beta)$.

 \item \label{HinN-N_A(H)-alpha}
 Suppose $\dim \Lie H \le 2$ and $\phi_h \neq 0$, $y_h = 0$, $\yy_h = 0$,
and $|x_h|^2 + 2 \Re(\phi_h \cjg{\eta_h}) = 0$ for every nonzero $h \in
\Lie H$.
 \begin{enumerate}
 \item \label{HinN-N_A(H)-alpha-A}
 If $\Lie H \subset \Lie U_{\alpha}$, then $N_A(H) = A$.
 \item \label{HinN-N_A(H)-alpha-b}
 If $\Lie H \subset \Lie U_{\alpha} + \Lie
U_{\alpha+\beta} + \Lie U_{\alpha+2\beta}$, but $\Lie H \not\subset \Lie
U_{\alpha}$, then $N_A(H) =
\ker(\beta)$.
 \item \label{HinN-N_A(H)-alpha-a+2b}
 If $\Lie H \subset \Lie U_{\alpha} + \Lie U_{2\alpha+2\beta}$, but
$\Lie H \not\subset \Lie U_{\alpha}$, then $N_A(H) = \ker(\alpha+2\beta)$.
 \end{enumerate}

 \end{enumerate}
 \end{Corollary}

\begin{proof} 
 It is clear that each of the given subgroups is normalized by the indicated
torus. We now show that the list is complete, and that no larger subtorus
of~$A$ normalizes~$H$.

Assume $N_A(H)$ is nontrivial. We proceed in cases, determined by
Theorem~\ref{HinN-notCDS}.

\setcounter{case}{0}

\begin{case} \label{HinN-N_A(H)-H=Zpf}
 Assume \fullref{HinN-notCDS}{H=Z}.
 \end{case}
 We may assume $\Lie H$ is neither $\Lie U_{2\beta}$ nor~$\Lie
U_{2\alpha+2\beta}$ (otherwise \pref{HinN-N_A(H)-H=Z-A} applies). Then,
because $|\eta_u|^2 = \xx_u \yy_u$ for every $u \in \Lie H$, we see that
$\eta_u \neq 0$ for every nonzero $u \in \Lie H$. Thus, the projection
of~$\Lie H$ to~$\Lie U_{\alpha+2\beta}$ is nontrivial. However, because 
$|\eta_u|^2 = \xx_u \yy_u$, we have $\Lie H \cap \Lie U_{\alpha+2\beta} =
0$. We know that $\Lie H \subset \Lie
U_{\alpha+2\beta} + \Lie U_{2\beta} + \Lie U_{2\alpha+2\beta}$ (because
$\Lie H = \Lie Z$), so, because each of $2\beta$ and $2\alpha+2\beta$ differs
from $\alpha+2\beta$ by~$\alpha$, we conclude that $N_A(H) =
\ker(\alpha)$,  so \pref{HinN-N_A(H)-H=Z-a} applies.

\begin{case}
 Assume \fullref{HinN-notCDS}{phi&y=0}.
 \end{case}
 Let $V$ be the projection of~$\Lie H$ to $\Lie U_{\alpha+\beta} + \Lie
U_{2\beta}$. Because $\yy_u \neq 0$, we know that $V$ projects nontrivially
to~$\Lie U_{2\beta}$. However, because $\Lie Z \subset \Lie
U_{2\alpha+2\beta}$, we also know that $V \cap \Lie U_{2\beta} = 0$.
Therefore $N_A(H) = \ker(\alpha-\beta)$. Then, because neither
$\alpha+2\beta$ nor $2\alpha+2\beta$ differs from $\alpha+\beta$ by a
multiple of $\alpha-\beta$, we conclude that $\Lie H = \bigl( \Lie H \cap
(\Lie U_{\alpha+\beta} + \Lie U_{2\beta}) \bigr) + \Lie Z$, so
\pref{HinN-N_A(H)-phi&y=0} applies.

\begin{case}
 Assume \fullref{HinN-notCDS}{phi=0}.
 \end{case}
  We may assume $\Lie H \neq \Lie Z$ (otherwise Case~\ref{HinN-N_A(H)-H=Zpf}
applies).

 \begin{subcase}
 Assume $\lambda \neq 0$.
 \end{subcase}
 Because $\Lie H \neq \Lie Z$, the projection of $\Lie H$ to $\Lie U_{\beta}
+ \Lie U_{\alpha+\beta}$ is nontrivial. However, because $\lambda \neq 0$,
this projection intersects neither  $\Lie U_{\beta}$ nor~$\Lie
U_{\alpha+\beta}$. Therefore $N_A(H) \subset \ker(\alpha)$. Then, because
neither $2\beta$, $\alpha+2\beta$, nor $2\alpha+2\beta$ differs from~$\beta$
by a multiple of~$\alpha$, we conclude that $\Lie
H = \bigl( \Lie H \cap (\Lie U_{\beta} + \Lie U_{\alpha+\beta}) \bigr) +
\Lie Z$, so \pref{HinN-N_A(H)-phi=0} applies.

\begin{subcase}
 Assume $\lambda = 0$.
 \end{subcase}
 This means $x_u = 0$ for every $u \in \Lie H$, and $\Lie Z \subset \Lie
U_{2\beta}$.  

Because $\Lie H \neq \Lie Z$, we know that $\Lie H$ projects nontrivially to
$\Lie U_{\beta}$. Because $\Lie Z \subset \Lie U_{2\beta}$, we know that
$\Lie H \cap \Lie U_{\alpha+2\beta} = \Lie H \cap \Lie U_{2\alpha+2\beta} =
0$. Thus, it is easy to see that if $\Lie H$ projects nontrivially to $\Lie
U_{\alpha+2\beta}$ or~$\Lie U_{2\alpha+2\beta}$ then either
\pref{HinN-N_A(H)-phi=0&x=0-notA-a+b}
or~\pref{HinN-N_A(H)-phi=0&x=0-notA-2a+b} applies. 

Thus, we may assume $\Lie H \subset \Lie U_{\beta} + \Lie
U_{2\beta}$.  If $\Lie Z \neq 0$, then $\Lie H = ( \Lie H \cap \Lie
U_{\beta}) + \Lie U_{2\beta}$, so \pref{HinN-N_A(H)-phi=0&x=0-A}
applies. Otherwise, \pref{HinN-N_A(H)-phi=0&x=0-notA-b}
applies.

\begin{case}
 Assume \fullref{HinN-notCDS}{y&yy=0}.
 \end{case}
 
\begin{subcase}
 Assume the projection of~$\Lie H$ to~$\Lie U_\alpha$ is trivial.
 \end{subcase}
  Because 
 $$|x_u|^2 = |x_u|^2 + 2 \Re(\phi_u \cjg{\eta_u}) \neq 0$$
 for every $u \in \Lie H \setminus \Lie U_{2\alpha+2\beta}$, we know that $x_u
\neq 0$ for every $u \in \Lie H \setminus \Lie U_{2\alpha+2\beta}$. Thus, if the
projection of~$\Lie H$ to~$\Lie U_{\alpha+2\beta}$ is nontrivial, then $N_A(H)
= \ker(\beta)$, and we see that \pref{HinN-N_A(H)-y&yy=0-b} applies. If not,
then $\Lie H \subset \Lie U_{\alpha+\beta} + \Lie U_{2\alpha+2\beta}$, so
either \pref{HinN-N_A(H)-y&yy=0-A} or \pref{HinN-N_A(H)-y&yy=0-a+b} applies.

\begin{subcase}
 Assume the projection of~$\Lie H$ to~$\Lie U_\alpha$ is nontrivial.
 \end{subcase}
 Let $V$ be the projection of~$\Lie H$ to $\Lie U_\alpha + \Lie
U_{\alpha+\beta} + \Lie U_{\alpha+2\beta}$. Because $|x_u|^2 + 2 \Re(\phi_u
\cjg{\eta_u}) \neq 0$ for every $u \in \Lie H \setminus \Lie
U_{2\alpha+2\beta}$, we know that $V \cap \Lie U_{\alpha} = 0$. Then,
because $\alpha$, $\alpha+\beta$, and~$\alpha+2\beta$ all differ by
multiples of~$\beta$, we conclude that $N_A(H) = \ker(\beta)$. Therefore
\pref{HinN-N_A(H)-y&yy=0-b} applies.

\begin{case}
 Assume \fullref{HinN-notCDS}{phi=yy}.
 \end{case}
 Let $V$ be the projection of~$\Lie H$ to $\Lie U_{\alpha} + \Lie
U_{2\beta}$. Because $\phi_h = \phi_0 \yy_h$, we see that $V \cap \Lie
U_{\alpha} = 0$ and $V \cap \Lie U_{2\beta} = 0$. Therefore $N_A(H) =
\ker(\alpha-2\beta)$.

Because no other roots differ by a multiple of $\alpha-2\beta$ (and $\Lie Z
= 0$), we conclude that $\Lie H = \bigl( \Lie H \cap ( \Lie U_\alpha + \Lie
U_{2\beta}) \bigr) + (\Lie H \cap \Lie U_{\alpha+\beta})$. Thus,
\pref{HinN-N_A(H)-phi=yy} applies.

\begin{case}
 Assume \fullref{HinN-notCDS}{phi=0(dim2)}.
 \end{case}
 
\begin{subcase}
 Assume $\Lie H \neq \Lie Z$.
 \end{subcase}
 Let $V$ be the projection of~$\Lie H$ to $\Lie U_\beta + \Lie
U_{\alpha+\beta}$. From the assumption of this subcase, we know $V \neq 0$.
However, because $\dim_{\complex} \langle x_u, y_u \rangle \neq 1$ for
every $u \in \Lie H$, we know that $V \cap \Lie U_\beta = 0$ and $\Lie H
\cap \Lie U_{\alpha+\beta} = 0$. Therefore $N_A(H) = \ker(\alpha)$, so
\pref{HinN-N_A(H)-phi=0(dim2)-a} applies.
 
\begin{subcase}
 Assume $\Lie H = \Lie Z$.
 \end{subcase}
 We may assume $\Lie H \not\subset \Lie U_{\alpha+2\beta}$ (otherwise
\pref{HinN-N_A(H)-phi=0(dim2)-A} applies). Therefore, $\Lie H$ projects
nontrivially to $\Lie U_{2\beta} + \Lie U_{2\alpha+2\beta}$. However,
because $|\eta_z|^2 \neq \xx_z \yy_z$,
for every nonzero $z \in \Lie Z$, we know that $V \cap \Lie U_{2\beta} = 0$
and $V \cap \Lie U_{2\alpha+2\beta} = 0$. Because $2\beta$, $\alpha+2\beta$,
and $2\alpha+2\beta$ all differ by multiples of~$\alpha$, we conclude that
$N_A(H) = \ker(\alpha)$, so \pref{HinN-N_A(H)-phi=0(dim2)-a} applies.

\begin{case}
 Assume \fullref{HinN-notCDS}{phi=0(dim1&2)}.
 \end{case}

\begin{subcase}
 Assume $N_A(H) = \ker(\alpha)$.
 \end{subcase}
 Because $\alpha+\beta$ is the only root that differs from $\beta$ by a
multiple of $\alpha$, we must have $\Lie H = \bigl( \Lie H \cap ( \Lie
U_{\beta} + \Lie U_{\alpha+\beta}) \bigr) + \Lie Z$. Thus, there is some $w
\in \Lie H$, such that $x_w = x_v$ and $y_w = y_v$, but the projection
of~$w$ to $\Lie U_{2\beta} + \Lie U_{\alpha+2\beta} + \Lie
U_{2\alpha+2\beta}$ is zero. This contradicts the fact that  $\xx_w |y_w|^2
+ \yy_w |x_w|^2 + 2 \Im(x_w y_w^{\dagger} \cjg{\eta_w}) \neq 0$.

\begin{subcase}
 Assume $N_A(H) \neq \ker(\alpha)$.
 \end{subcase}
 Because $2\beta$, $\alpha+2\beta$, and $2\alpha+2\beta$ all differ by
multiples of~$\alpha$, we must have $\Lie Z = (\Lie Z \cap \Lie U_{2\beta})
+ (\Lie Z \cap \Lie U_{\alpha + 2\beta}) + (\Lie Z \cap \Lie
U_{2\alpha+2\beta})$. Then, because $|\eta_z|^2 \neq \xx_z \yy_z$ for every
nonzero $z \in \Lie Z$, we conclude that $\Lie Z \subset \Lie
U_{\alpha+2\beta}$.

 Let $V$ be the projection of~$\Lie H$ to $\Lie U_{\beta} + \Lie
U_{\alpha+\beta}$. Because $\beta$ and $\alpha+\beta$ differ by~$\alpha$, we
know that $V = (V \cap \Lie U_{\beta}) + (V \cap \Lie U_{\alpha+\beta})$. 

\begin{subsubcase}
 Assume $x_v \neq 0$.
 \end{subsubcase}
 Because $V = (V \cap \Lie U_{\beta}) + (V \cap \Lie U_{\alpha+\beta})$,
there is some $w \in V$, such that $x_w \neq 0$ and $y_w = 0$. For every
such~$w$, because $\xx_w |y_w|^2 + \yy_w |x_w|^2 + 2 \Im(x_w y_w^{\dagger}
\cjg{\eta_w}) \neq 0$, we know that $\yy_w \neq 0$. Thus, we see that
$N_A(H) = \ker \bigl( (\alpha+\beta) - 2\beta \bigr) = \ker(\alpha-\beta)$.

We know that $\Lie H \cap \Lie U_\beta = 0$, that $\Lie H$ projects trivially
to~$\Lie U_\alpha$, and that $\alpha$ is the only root that differs
from~$\beta$ by a multiple of $\alpha-\beta$, so we conclude that $y_h = 0$
for every $h \in H$.

We now see that \pref{HinN-N_A(H)-phi=0(dim1&2)-x} applies.

\begin{subsubcase}
 Assume $y_v \neq 0$.
 \end{subsubcase}
 This is similar to the preceding subsubcase (indeed, they are conjugate under
the Weyl reflection corresponding to the root~$\alpha$); we see that
\pref{HinN-N_A(H)-phi=0(dim1&2)-y} applies.

\begin{case}
 Assume \fullref{HinN-notCDS}{onlyphiyyxx}.
 \end{case}
 By considering the projection of~$\Lie H$ to $\Lie U_{\alpha} + \Lie
U_{\beta}$, and noting that $\phi_h \asymp y_h$ for every $h \in H$, we see
that $N_A(H) = \ker(\alpha-\beta)$. The only other pair of roots that differ
by a multiple of $\alpha-\beta$ is $\{\alpha+\beta, 2\beta\}$. Thus, we
see that \pref{HinN-N_A(H)-onlyphiyyxx} applies.

\begin{case}
 Assume \fullref{HinN-notCDS}{dim3}.
 \end{case}
 By considering the projection of~$\Lie H$ to $\Lie U_{\alpha} + \Lie
U_\beta$, we see that $N_A(H) = \ker(\alpha-\beta)$. Because $\phi_u \neq 0$
for every $u \in \Lie H \setminus \Lie U_{2\alpha+2\beta}$, but $\beta$ is
the only root that differs from~$\alpha$ by a multiple of $\alpha - \beta$,
we conclude that $\Lie H$ projects trivially into every root space except
$\Lie U_{\alpha}$, $\Lie U_{\beta}$, and~$\Lie U_{2\alpha+2\beta}$. Thus
\pref{HinN-N_A(H)-dim3} applies.

\begin{case}
 Assume \fullref{HinN-notCDS}{alpha}.
 \end{case}
 We may assume $\Lie H \not\subset \Lie U_\alpha$ (otherwise
\pref{HinN-N_A(H)-alpha-A} applies). Thus, there is some root $\sigma \neq
\alpha$, such that the projection of~$\Lie H$ to~$\Lie U_\sigma$ is
nontrivial. However, because  $\phi_h \neq 0$ for every nonzero $h \in
\Lie H$, we know that $\Lie H \cap \Lie U_\sigma = 0$. Thus, $N_A(H) =
\ker(\alpha-\sigma)$.

Because $y_h = 0$ and $\yy_h = 0$ for every nonzero $h \in \Lie H$, we know
that $\sigma \neq \beta$ and $\sigma \neq 2\beta$. If $\sigma = \alpha +
\beta$ or $\sigma = \alpha + 2\beta$, we obtain \pref{HinN-N_A(H)-alpha-b}.
If $\sigma = 2\alpha+2\beta$, we obtain \pref{HinN-N_A(H)-alpha-a+2b}.

\begin{case}
 Assume \fullref{HinN-notCDS}{u+z}.
 \end{case}
 Because $\phi_u \neq 0$ and $y_u \neq 0$, we must have $N_A(H) =
\ker(\alpha-\beta)$. Then, because $\alpha+\beta$ does not differ
from~$\alpha$ by a multiple of $\alpha-\beta$, we conclude that $x_u = 0$.

Because $\eta_z \neq 0$, but no root differs
from~$\alpha+2\beta$ by a multiple of $\alpha-\beta$, we conclude that $\Lie
H \cap \Lie U_{\alpha+2\beta} \neq 0$. Because $\Lie Z$ is one-dimensional,
this implies $z \in \Lie U_{\alpha+2\beta}$, so $\xx_z = 0$.

Since $\xx_z = 0$ and $x_u = 0$, we conclude, from the inequality $\xx_z
|y_u|^2 - \phi_u \yy_u \cjg{\eta_z} + 2 \Im \bigl( \cjg{\eta_z} x_u
y_u^{\dagger} \bigr) \neq 0$, that $\yy_u \neq 0$. This is a contradiction,
because $2\beta$ does not differ from~$\alpha$ by a multiple of
$\alpha-\beta$, and $\Lie H \cap \Lie U_{2\beta} = 0$ (because, as shown
above, $\Lie Z \subset \Lie U_{\alpha+2\beta}$).
 \end{proof}

\section{Subgroups that are not contained in~$N$} \label{CDSnotinN-section}

Let $H$ be a closed, connected subgroup of~$G$ that is not contained in~$N$. In
this section, we determine whether $H$ is a Cartan-decomposition subgroup
or not (and, if not, we calculate $\mu(H)$).

Lemma~\ref{H'inAN} shows that we may assume $H \subset AN$,  and then
Lemma~\ref{conj-to-compatible} shows that we may assume $H$ satisfies the
technical condition of being compatible with~$A$. (Both of these lemmas are
well known.)
 Furthermore, we may assume that  $H \cap N$ is \textbf{not} a
Cartan-decomposition subgroup, and that $A \not\subset H$ (otherwise, it is
obvious that $H$ is a Cartan-decomposition subgroup). 

Theorem~\ref{TU-proj} describes $\mu(H)$ for every such subgroup that is a
semidirect product $(H \cap A) \ltimes (H \cap N)$; and
Proposition~\ref{HnotTUproj} describes $\mu(H)$ for the other subgroups (except
that the one-dimensional case appears in Lemma~\ref{dim1-mu(H)}).

\begin{Lemma}[{\cite[Lem.~2.9]{OhWitte-CDS}}] \label{H'inAN}
 Let $H$ be a closed, connected subgroup of a connected, almost simple, linear,
real Lie group~$G$. There is a closed, connected subgroup~$H'$ of~$G$ and a
compact subgroup~$C$ of~$G$, such that $CH = CH'$, and $H'$ is conjugate to a
subgroup of~$AN$.
 \end{Lemma}

\begin{Definition} \label{compatible}
 Let us say that a subgroup~$H$ of $AN$ is \emph{compatible} with~$A$ if
$H \subset T U C_N(T)$, where $T = A \cap (HN)$, $U = H \cap N$, and
$C_N(T)$ denotes the centralizer of~$T$ in~$N$. 
 \end{Definition}

\begin{Lemma}[{\cite[Lem.~2.3]{OhWitte-CDS}}] \label{conj-to-compatible}
 If $H$ is a closed, connected subgroup of~$AN$, then $H$ is conjugate,
via an element of~$N$, to a subgroup that is compatible with~$A$.
 \end{Lemma}

\begin{Theorem} \label{TU-proj}
 Assume that $G = \SU(2,n)$. Here is a list of every closed, connected,
nontrivial subgroup~$H$ of $AN$, such that $H$ is of the form $H = T \ltimes U$,
where $T$ is a one-dimensional subgroup of~$A$, and $U$ is a nontrivial subgroup
of~$N$ that is \textbf{not} a Cartan-decomposition subgroup.
 \begin{enumerate}

 \item \label{TU-proj-H=Z}
 Suppose $\dim \Lie U = 1$, $\Lie U = \Lie Z$, and we have $|\eta_h|^2 = \xx_h
\yy_h$ for every $h \in U$.
 \begin{enumerate}
 \item \label{TU-proj-H=Z-A}
 If $\Lie U = \Lie U_{2\beta}$ or $\Lie U = \Lie U_{2\alpha+2\beta}$, then
$\mu(H)$ is described in \cite[Prop.~3.17 or Cor.~3.18]{OhWitte-CDS}.
 \item \label{TU-proj-H=Z-a}
 Otherwise, $T = \ker(\alpha)$, and $H$ is a Cartan-decomposition subgroup.
 \end{enumerate}

 \item \label{TU-proj-phi&y=0}
 Suppose $\Lie U = \bigl( \Lie U \cap (\Lie U_{\alpha+\beta} + \Lie U_{2\beta})
\bigr) + \Lie Z$, $\Lie Z \subset \Lie
U_{2\alpha + 2\beta}$, there is some
$v \in \Lie U$, such that $\yy_v \neq 0$, and $T = \ker(\alpha-\beta)$. Then
$\mu(H) \approx \muH{\|h\|}{\|h\|^{3/2}}$, unless $\dim H = 2$, in which case
$\rho(h) \asymp \|h\|^{3/2}$ for every $h \in H$.

 \item \label{TU-proj-phi=0}
 Suppose $\Lie U = \bigl( \Lie U \cap (\Lie U_{\beta} + \Lie U_{\alpha+\beta})
\bigr) + \Lie Z$, $T = \ker(\alpha)$, and there is some nonzero $\lambda \in
\complex$, such that we have $x_u = \lambda y_u$ for every $u \in U$, and we
have $\eta_z = i\lambda \yy_z$ and $\xx_z = |\lambda|^2 \yy_z$ for every $z \in
\Lie Z$.
 Then $H$ is a Cartan-decomposition subgroup.

 \item \label{TU-proj-phi=0&x=0}
 Suppose $\phi_u = 0$ and $x_u = 0$ for every $u \in \Lie U$, we have $\Lie
Z \subset \Lie U_{2 \beta}$, and $\Lie U \neq \Lie Z$.
 \begin{enumerate}
 \item \label{TU-proj-phi=0&x=0-A}
 If $\Lie U = ( \Lie U \cap \Lie U_{\beta}) + \Lie Z$, then $\mu(H)$ is
described in \cite[Prop.~3.17 or Cor.~3.18]{OhWitte-CDS}.
 \item Otherwise:
 \begin{enumerate}
 \item \label{TU-proj-phi=0&x=0-notA-a+b}
 If $\Lie U = \bigl( \Lie U \cap
(\Lie U_{\beta} + \Lie U_{\alpha+2\beta}) \bigr) + \Lie Z$, then $T =
\ker(\alpha+\beta)$, and $\rho(h) \asymp h$ for every $h \in H$.
 \item \label{TU-proj-phi=0&x=0-notA-2a+b}
  If $\Lie U = \bigl( \Lie U \cap (\Lie U_{\beta} + \Lie
U_{2\alpha+2\beta}) \bigr) + \Lie Z$, then $T = \ker(2\alpha+\beta)$, and
$\mu(H) \approx \muH{\|h\|}{\|h\|^{3/2}}$, unless $\dim H = 2$, in which
case $\rho(h) \asymp \|h\|^{3/2}$ for every $h \in H$.
 \item \label{TU-proj-phi=0&x=0-notA-b}
  If $\Lie Z = 0$ and $\Lie U \subset \Lie U_{\beta} + \Lie
U_{2\beta}$, then $T = \ker(\beta)$, and $H$ is a Cartan-decomposition
subgroup.
 \end{enumerate}
 \end{enumerate}

 \item \label{TU-proj-y&yy=0}
 Suppose $y_u = 0$, $\yy_u = 0$, and $|x_u|^2 + 2 \Re(\phi_u \cjg{\eta_u})
\neq 0$ for every $u \in U \setminus U_{2\alpha+2\beta}$.
 \begin{enumerate}
 \item \label{TU-proj-y&yy=0-A}
  If $\Lie U = ( \Lie U \cap \Lie U_{\alpha+\beta}) + \Lie Z$, then $\mu(H)$
is described in \cite[Prop.~3.17 or Cor.~3.18]{OhWitte-CDS}.
 \item \label{TU-proj-y&yy=0-a+b}
  If $\Lie U \subset \Lie U_{\alpha+\beta} + \Lie
U_{2\alpha+2\beta}$, but $\Lie U \neq ( \Lie U \cap \Lie
U_{\alpha+\beta}) + \Lie Z$, then $T = \ker(\alpha+\beta)$, and $H$ is a
Cartan-decomposition subgroup.
 \item \label{TU-proj-y&yy=0-b}
  If $\Lie U = \bigl( \Lie U \cap (\Lie U_{\alpha} + \Lie U_{\alpha+\beta} +
\Lie U_{\alpha+2\beta}) \bigr) + \Lie Z$, but $\Lie U \not\subset \Lie
U_{\alpha+\beta} + \Lie U_{2\alpha+2\beta}$, then $T = \ker(\beta)$, and
$\rho(h) \asymp h$ for every $h \in H$.
 \end{enumerate}

 \item \label{TU-proj-phi=yy}
 Suppose $\Lie U = \bigl( \Lie U \cap ( \Lie U_\alpha + \Lie U_{2\beta}) \bigr)
+ (\Lie U \cap \Lie U_{\alpha+\beta})$,
 $T = \ker(\alpha-2\beta)$,
 $\Lie U \not\subset \Lie U_{\alpha+\beta}$, and there is some nonzero $\phi_0
\in \complex$, such that $\phi_u = \phi_0 \yy_u$ for every $u \in U$. Then
$\mu(H) \approx \muH{\|h\|}{\|h\|^{4/3}}$, unless $\dim H = 2$, in which case,
$\rho(h) \asymp \|h\|^{4/3}$ for every $h \in H$.

 \item \label{TU-proj-phi=0(dim2)}
 Suppose $\phi_u = 0$ and $\dim_{\complex} \langle x_u, y_u \rangle \neq 1$ for
every $u \in U$, and $|\eta_z|^2 \neq \xx_z \yy_z$,
for every nonzero $z \in \Lie Z$.
 \begin{enumerate}
 \item \label{TU-proj-phi=0(dim2)-A}
 If $\Lie U \subset \Lie U_{\alpha+2\beta}$, then $\mu(H)$ is described in
\cite[Prop.~3.17 or Cor.~3.18]{OhWitte-CDS}.
 \item \label{TU-proj-phi=0(dim2)-a}
 If $\Lie U \not\subset \Lie U_{\alpha+2\beta}$, and $\Lie U = \bigl(
\Lie U \cap ( \Lie U_{\beta} + \Lie U_{\alpha+\beta}) \bigr) + \Lie Z$, then
$T = \ker(\alpha)$, and $\rho(h) \asymp \|h\|^2$ for every $h \in H$.
 \end{enumerate}

 \item \label{TU-proj-phi=0(dim1&2)}
 Suppose $\phi_u = 0$ for every $u \in U$, there exist nonzero $v_1,v_2
\in \Lie U$, such that  $\dim_{\complex} \langle x_{v_1}, y_{v_1} \rangle \neq
1$ and  $\dim_{\complex} \langle x_{v_2}, y_{v_2} \rangle = 1$, and we have
 $\xx_{v_2} |y_{v_2}|^2 + \yy_{v_2} |x_{v_2}|^2 + 2 \Im(x_{v_2}
y_{v_2}^{\dagger} \cjg{\eta_{v_2}}) \neq 0$
 for every such ${v_2} \in \Lie U$,
 and $|\eta_z|^2 \neq \xx_z \yy_z$,
for every nonzero $z \in \Lie Z$.
 \begin{enumerate}
 \item \label{TU-proj-phi=0(dim1&2)-x}
 If $\Lie U = \bigl( \Lie U \cap (\Lie
U_{\alpha+\beta} + \Lie U_{2\beta}) \bigr) + (\Lie U \cap \Lie
U_{\alpha+2\beta})$, then $T = \ker(\alpha-\beta)$ and
$\mu(H) \approx \muH{\|h\|^{3/2}}{\|h\|^2}$.
 \item \label{TU-proj-phi=0(dim1&2)-y}
 If $\Lie U = \bigl( \Lie U \cap (\Lie
U_{\beta} + \Lie U_{2\alpha+2\beta}) \bigr) + (\Lie U \cap \Lie
U_{\alpha+2\beta})$, then $T = \ker(2\alpha+\beta)$ and $\mu(H) \approx
\muH{\|h\|^{3/2}}{\|h\|^2}$.
 \end{enumerate}

 \item \label{TU-proj-onlyphiyyxx}
 Suppose $\dim \Lie U \le 3$, 
 $\Lie U = \bigl( \Lie U \cap ( \Lie U_\alpha + \Lie U_{\beta})
\bigr) +  \bigl( \Lie U \cap ( \Lie U_{\alpha+\beta} + \Lie U_{2\beta})
\bigr)$,
 $\Lie U \cap ( \Lie U_\alpha + \Lie U_{\beta}) \neq 0$, and
 we have $\phi_u \asymp y_u$ and $x_u \asymp \yy_u$ for
$u \in U$.
 Then $T = \ker(\alpha-\beta)$, and $\rho(h) \asymp \|h\|^{3/2}$ for
every $h \in H$.

 \item \label{TU-proj-dim3}
 Suppose $\dim \Lie U = 2$, $\Lie Z = \Lie U_{2\alpha+2\beta}$, $\phi_u
\neq 0$ and  $y_u \neq 0$ for every $u \in  U \setminus Z$.
 If $\Lie U = \bigl( \Lie U \cap ( \Lie U_\alpha + \Lie U_{\beta})
\bigr) + \Lie Z$, then $T = \ker(\alpha-\beta)$, and $\mu(H) \approx
\muH{\|h\|}{\|h\|^{3/2}}$.

 \item \label{TU-proj-alpha}
 Suppose $\dim \Lie U \le 2$ and $\phi_u \neq 0$, $y_u = 0$, $\yy_u = 0$,
and $|x_u|^2 + 2 \Re(\phi_u \cjg{\eta_u}) = 0$ for every nontrivial $u \in U$.
 \begin{enumerate}
 \item \label{TU-proj-alpha-A}
 If $\Lie U \subset \Lie U_{\alpha}$, then $\mu(H)$ is described in
\cite[Prop.~3.17 or Cor.~3.18]{OhWitte-CDS}.
 \item \label{TU-proj-alpha-b}
 If $\Lie U \subset \Lie U_{\alpha} + \Lie U_{\alpha+\beta} + \Lie
U_{\alpha+2\beta}$, but $\Lie U \not\subset \Lie U_{\alpha}$, then $T =
\ker(\beta)$, and $H$ is a Cartan-decomposition subgroup.
 \item \label{TU-proj-alpha-a+2b}
 If $\Lie U \subset \Lie U_{\alpha} + \Lie U_{2\alpha+2\beta}$, but $\Lie U
\not\subset \Lie U_{\alpha}$, then $T = \ker(\alpha+2\beta)$, and $\rho(h)
\asymp \|h\|^2$ for every $h \in H$.
 \end{enumerate}

 \end{enumerate}
 \end{Theorem}

\begin{proof}
 For $h \in H$, we wish to approximately calculate $\|\rho(h)\|$. We write $h
= au$ with $a \in T$ and $u \in U$. 
 Writing $a = \diag(a_1,a_2,\ldots,a_{n+2})$, we always assume either that $a_1
> 1$ or that $a_1 = 1$ and $a_2 \ge 1$ (perhaps replacing $h$
with~$h^{-1}$---because $\|\rho(h)\| = \|\rho(h^{-1})\|$, this causes no
harm).

Because $T$ normalizes~$U$, we know that $U$ is a subgroup that is listed in
Corollary~\ref{HinN-N_A(H)}, and we have $T \subset N_G(U)$. This leads to the
various cases listed in the statement of the theorem. 

 \pref{TU-proj-H=Z-a} We have $\rho(u) \asymp u$ for $u \in U$ and $\rho(a)
\asymp \|a\|^2$ for $a \in T$, so $H$ is a Cartan-decomposition subgroup.

 \pref{TU-proj-phi&y=0} We have $|\phi_u| + |y_u| + |\eta_u| +
|\xx_u| = 0$ and $\yy_u = O(x_u)$, so $u_{i,j} = O \bigl( 1 + |x_u| \bigr)$
whenever $(i,j) \neq (1,n+2)$. Then, because $a_1 = a_2^2$, we see that 
 $$ u_{i,j} = O \bigl[ a_2 \bigl( 1+ |x_u| \bigr) \bigr]
 = O \bigl( |h_{1,1}|^{1/2} + | h_{1,n+2}|^{1/2} \bigr)
 = O \bigl( \|h\|^{1/2} \bigr) $$
 whenever $i>1$. Therefore $\rho(h) = O
\bigl( \|h\|^{3/2} \bigr)$. This completes the proof if $\dim H > 2$ (that
is, if $\dim U > 1$).

If $\dim U = 1$, then $\yy_u \asymp x_u$ and $\xx_u = 0$. We have 
 $\|h\| = a_1 \bigl( 1 + |x_u|^2 \bigr)$,
 $$ \Delta(h) = a_1 a_2 \left[ i \left( \frac{1}{2} |x_u|^2 \yy_u \right)
\right]
 \asymp \bigl(a_1 |x_u|^2 \bigr)^{3/2} $$
 and
 $$ \det \begin{pmatrix}
 h_{1,1} & h_{1,2} \\
 h_{2,1} & h_{2,2}
 \end{pmatrix}
 = a_1 a_2 = a_1^{3/2} .$$
 Thus, $\|h\|^{3/2} = O \bigl( \rho(h) \bigr)$. We conclude that $\rho(h)
\asymp \|h\|^{3/2}$.

 \pref{TU-proj-phi=0}
 Replacing $H$ by a conjugate under~$U_\alpha$, we may replace~$H$ with a
similar subgroup~$H'$ with $\lambda = 0$. Thus, $H' = T \ltimes U'$ with $U'
\subset U_\beta U_{2\beta}$. Then \cite[Prop.~3.17]{OhWitte-CDS} implies $H$
is a Cartan-decomposition subgroup.

 \pref{TU-proj-phi=0&x=0-notA-a+b} The Weyl reflection corresponding to the
root~$\alpha$ conjugates~$H$ to a subgroup of
type~\pref{TU-proj-y&yy=0-b}.

 \pref{TU-proj-phi=0&x=0-notA-2a+b} The Weyl reflection corresponding to the
root~$\alpha$ conjugates~$H$ to a subgroup of
type~\pref{TU-proj-phi&y=0}.

 \pref{TU-proj-phi=0&x=0-notA-b} \cite[Prop.~3.17]{OhWitte-CDS} implies $H$
is a Cartan-decomposition subgroup.

 \pref{TU-proj-y&yy=0-a+b} \cite[Prop.~3.17]{OhWitte-CDS} implies $H$
is a Cartan-decomposition subgroup.

 \pref{TU-proj-y&yy=0-b}
 We have
 $$ h_{i,j} = 
 \begin{cases}
 O(1) & \mbox{if $i \neq 1$ and $j \neq n+2$} \\
 O(a_1 x) & \mbox{if $i = 1$ and $j \neq n+2$} \\
 O(x) & \mbox{if $i \neq 1$ and $j = n+2$} 
 \end{cases}$$
 and $h_{1,n+2} \asymp a_1 \bigl( |x|^2 + |\xx| \bigr)$.
 We conclude that $\rho(h) \asymp h$.

\pref{TU-proj-phi=yy} From the proof of~\fullref{HinN-nosquareproj}{phi=yy},
we know
 that $u \asymp u_{1,n+2}$,
 that $u_{i,j} = O \bigl( \|u\|^{3/2} \bigr)$ whenever $(i,j) \neq (1,n+2)$,
 and that $u_{i,j} = O \bigl( \|u\|^{1/3} \bigr)$ whenever $i \neq 1$ and $j
\neq n+2$.
 (In particular, $h \asymp a_1 \bigl( 1 + u_{1,n+2} \bigr)$.)
 Furthermore, we have $a_1 = a_2^3$. Therefore
 $$ \rho(h) \asymp a_1 a_2 \rho(u) \asymp a_1^{4/3} \rho(u) .$$
 The desired conclusion follows.

 \pref{TU-proj-phi=0(dim2)-a}
 From Lemma~\ref{HinN-phi=0(dim<xy>)}, we know
 $\|u\|^2 = O \bigl( 1 + |\Delta(u)| \bigr)$.
 Then, because
  $$ \det \begin{pmatrix}
 h_{1,1} & h_{1,2} \\
 h_{2,1} & h_{2,2}
 \end{pmatrix}
 = a_1 a_2 = a_1^2 $$
 and $\Delta(au) = a_1^2 \Delta(u)$, we have
 $$ \|h\|^2 = O \bigl( a_1^2 \|u\|^2 \bigr)
 = O \bigl( a_1^2 + |\Delta(h)| \bigr)
 = O \bigl( \rho(h) \bigr) .$$

\pref{TU-proj-phi=0(dim1&2)} Assume \pref{TU-proj-phi=0(dim1&2)-x}. (The
other case, \pref{TU-proj-phi=0(dim1&2)-y}, is conjugate to this one by the
Weyl reflection corresponding to the root~$\alpha$.) From
Lemma~\ref{HinN-phi=0(dim<xy>)}, we have $\|u\|^{3/2} = O \bigl( 1 +
|\Delta(u)| \bigr)$. Then, because $a_1 = a_2^2$, we have
 $$ \|h\|^{3/2} = a_1 a_2 \|u\|^{3/2} 
 = O \bigl( a_1 a_2 + |\Delta(h)| \bigr) = O \bigl( \rho(h) \bigr) .$$

 \pref{TU-proj-onlyphiyyxx}
 From the proof of~\fullref{HinN-nosquareproj}{onlyphiyyxx}, we know
 $\rho(u) \asymp 1 + \Delta(u) \asymp \|u\|^{3/2}$. The proof is completed
as in~\pref{TU-proj-phi=0(dim1&2)}.

 \pref{TU-proj-dim3}
 Because $\phi_u \asymp y_u$, it is easy to see that
 $$h \asymp a_1 \bigl( 1 + |\phi_u|^2 |y_u|^2 + |\xx_u| \bigr)
 \asymp a_1 \bigl( 1 + |\phi_u|^4 + |\xx_u| \bigr) $$
 and
 $$ \Delta(h) \asymp a_1a_2 \bigl( |y_u|^4 |\phi_u|^2 + |\xx_u| |y_u|^2 \bigr)
 \asymp a_1^{3/2} \bigl( |\phi_u|^6 + |\xx_u| |\phi_u|^2 \bigr) 
 = O \bigl( \|h\|^{3/2} \bigr) .$$
 Then it is not difficult to see that
 $\rho(h) = O \bigl( \|h\|^{3/2} \bigr)$ for every $h \in H$.
 So $\mu(H) \approx \muH{\|h\|}{\|h\|^{3/2}}$.

 \pref{TU-proj-alpha-b} We have $\rho(a) \asymp
a$ for $a \in T$ and $\rho(u) \asymp \|u\|^2$ for $u \in U$, so $H$ is a
Cartan-decomposition subgroup.

\pref{TU-proj-alpha-a+2b} $H$ is conjugate (via an element
of~$U_{\alpha+2\beta}$) to $T \ltimes U_\alpha$. From
\cite[Prop.~3.18]{OhWitte-CDS}, we have $\rho(h) \asymp \|h\|^2$ for every
$h \in T \ltimes U_\alpha$. Therefore $\rho(h) \asymp \|h\|^2$ for every $h \in
H$.
 \end{proof}

\begin{Lemma}[{\cite[Lem.~2.4]{OhWitte-CDS}}] \label{not-semi}
 Assume $G = \SU(2,n)$, and let $H$ be a closed, connected subgroup of~$AN$
that is compatible with~$A$. Then either
 \begin{enumerate}
 \item $H = (H \cap A) \ltimes (H \cap N)$; or
 \item \label{not-semi-not}
 there is a positive root~$\omega$, a nontrivial group homomorphism
$\psi\colon \ker \omega \to U_\omega U_{2\omega}$, and a closed, connected
subgroup~$U$ of~$N$, such that
 \begin{enumerate}
 \item \label{not-semi-codim1}
 $H = \{\, a \psi(a) \mid a \in \ker \omega \,\} U$;
 \item \label{not-semi-disjoint}
 $U \cap \psi( \ker\omega) = e$; and
 \item \label{not-semi-normal}
 $U$ is normalized by both $\ker\omega$ and~$\psi( \ker\omega)$.
 \end{enumerate}
 \end{enumerate}
 \end{Lemma}

\begin{Proposition} \label{HnotTUproj}
 Assume that $G = \SU(2,n)$. Let $H$ be a closed, connected, nontrivial
subgroup of~$AN$, that is compatible with~$A$, such that
 \begin{itemize}
 \item $H \cap N$ is not a Cartan-decomposition subgroup;
 \item $H \neq (H \cap A) (H \cap N)$; and
 \item $\dim H > 1$.
 \end{itemize}
 Then there are positive roots~$\omega$ and~$\sigma$, and a one-dimensional
subspace~$\Lie X$ of $(\ker \omega) + \Lie U_\omega + \Lie U_{2\omega}$, such
that $\Lie H = \Lie X +
(\Lie H \cap \Lie N)$, $\Lie H \cap \Lie N \subset \Lie U_\sigma + \Lie
U_{2\sigma}$, and either:
 \begin{enumerate}
 \item \label{HnotTUproj-a,a+b}
  $\omega =\alpha$, $\sigma = \alpha+\beta$, and 
 $\mu(H) \approx \muH{\|h\|}{\|h\|^2/(\log \|h\|)}$; or
 \item \label{HnotTUproj-a,a+2b}
  $\omega =\alpha$ $\sigma = \alpha+2\beta$, and 
 $\mu(H) \approx \muH{\|h\|^2/(\log \|h\|)^2}{\|h\|^2}$; or
 \item \label{HnotTUproj-b,a+2b}
  $\omega =\beta$, $\sigma = \alpha+2\beta$, and 
 $\mu(H) \approx \muH{\|h\| \bigl( \log \|h\| \bigr)^{r/2}}{\|h\|^2}$,
 where 
 $$r = 
 \begin{cases}
 1 & \text{if $\Lie X \subset \Lie U_{2\beta}$} \\
 2 & \text{otherwise}
 \end{cases}
 $$
 or
 \item \label{HnotTUproj-b,a+b}
  $\omega =\beta$, $\sigma = \alpha+\beta$, and 
 $\mu(H) \approx \muH{\|h\|}{\|h\|  \bigl( \log \|h\| \bigr)^{r}}$,
 where $r$ is defined as above;
 or 
 \item  \label{HnotTUproj-CDS}
 $\Lie U \cap (\Lie U_{\omega} + \Lie U_{2\omega}) \neq 0$, in which case
$H$ is a Cartan-decomposition subgroup.
 \end{enumerate}
 \end{Proposition}

\begin{proof}
 We use the notation of Lemma~\ref{not-semi}: $T = \ker \omega$, $U
= H \cap N$, $\psi \colon T \to U_\omega U_{2\omega}$, and $H = \{a \psi(a) \}
\ltimes U$.

 We need only consider the cases in Corollary~\ref{HinN-N_A(H)} for which
$H$ (now called~$U$) is normalized by the kernel of some (reduced) positive
root. Here is a list of them.
 \begin{enumerate}

 \item $N_A(U) = \ker(\beta)$: \fullref{HinN-N_A(H)}{phi=0&x=0-notA-b},
\fullref{HinN-N_A(H)}{y&yy=0-b}, and \fullref{HinN-N_A(H)}{alpha-b}.

 \item $N_A(U) = \ker(\alpha+\beta)$: \fullref{HinN-N_A(H)}{phi=0&x=0-notA-a+b} and
\fullref{HinN-N_A(H)}{y&yy=0-a+b}.

 \item $N_A(U) = \ker(\alpha)$: \fullref{HinN-N_A(H)}{H=Z-a}, \fullref{HinN-N_A(H)}{phi=0},
and \fullref{HinN-N_A(H)}{phi=0(dim2)-a}.

 \item $N_A(U) = \ker(\alpha+2\beta)$: \fullref{HinN-N_A(H)}{alpha-a+2b}.

 \item $N_A(U) = A$: \fullref{HinN-N_A(H)}{H=Z-A}, \fullref{HinN-N_A(H)}{phi=0&x=0-A},
\fullref{HinN-N_A(H)}{y&yy=0-A}, \fullref{HinN-N_A(H)}{phi=0(dim2)-A}, and
\fullref{HinN-N_A(H)}{alpha-A}.

 \end{enumerate}
 Note that in each of the cases with $N_A(U) = A$, there is a (reduced)
positive root~$\sigma$, such that $\Lie U \subset \Lie U_\sigma + \Lie
U_{2\sigma}$.

\setcounter{case}{0}

\begin{case}
 Assume $\omega = \beta$.
 \end{case}
 
\begin{subcase}
 Assume \fullref{HinN-N_A(H)}{phi=0&x=0-notA-b}.
 \end{subcase}
 From \pref{HnotTU-omega+2omega}, we know that $H$ is a Cartan-decomposition
subgroup.
 
\begin{subcase} \label{HnotTU-beta-almostsame}
 Assume \fullref{HinN-N_A(H)}{y&yy=0-b}.
 \end{subcase}
 There is some $u \in U$, such that $\phi_u \neq 0$. Then, because $\psi(T)
\subset U_\beta U_{2\beta}$ normalizes~$U$, we must have $U \cap
U_{\alpha+2\beta} \neq e$. This is a contradiction.

\begin{subcase}
 Assume \fullref{HinN-N_A(H)}{alpha-b}.
 \end{subcase}
 Let $u \in \Lie U$. Because $U$ is normalized by $\psi(T)$, there is some
nonzero $v \in \Lie U_{\beta} + \Lie U_{2\beta}$, such that $v$
normalizes~$\Lie U$; thus, $[u,v] \in \Lie U$. Then, because $\phi_{[u,v]} =
0$, but $\phi_h \neq 0$ for every nontrivial $h \in U$, we conclude that
$[u,v] = 0$. However, $\phi_u \neq 0$, and either $y_v \neq 0$ or~$\yy_v \neq
0$, so either $x_{[u,v]} \neq 0$ or $\eta_{[u,v]} \neq 0$. This is a
contradiction.

\begin{subcase}
 Assume $N_A(U) = A$.
 \end{subcase}
 There is a positive root~$\sigma$, such that
$\Lie U \subset \Lie U_\sigma + \Lie U_{2\sigma}$.

If $\sigma = \beta$, then, from \pref{HnotTU-omega+2omega}, we know that $H$
is a Cartan-decomposition subgroup.

Suppose $\sigma = \alpha + 2\beta$. Clearly $\|h\| \asymp a_1 |\eta_u|$. Also,
 $$\rho(h) \asymp a_1|\eta_u|^2 + a_1 ( \log a_1 )^r ,$$
 where $r = 1$ if $\psi(T) \subset U_{2\beta}$ (i.e., if $y_h = 0$ for every
$h \in H$) and $r = 2$ if $\psi(T) \not\subset U_{2\beta}$.
 The smallest value of $\|\rho(h)\|$ relative to~$\|h\|$ is obtained by
taking $\eta_u \asymp (\log a_1)^{r/2}$, resulting in $\rho(h) \asymp \|h\|
\bigl( \log \|h\| \bigr)^{r/2}$. Then, since $\rho(u) \asymp \|u\|^2$ for $u
\in U$, we conclude that $\mu(H) \approx \muH{\|h\| \bigl( \log \|h\|
\bigr)^{r/2}}{\|h\|^2}$.

Because $U$ is normalized by the nontrivial subgroup $\psi(T)$ of $U_\beta
U_{2\beta}$, we know that $\sigma \neq \alpha$. Therefore, we may now assume $\sigma = \alpha + \beta$. We show that $\mu(H)
\approx \muH{\|h\|}{\|h\| (\log \|h\|)^r}$.
 For $u \in U$, we have $\rho(u) \asymp u$. For $a \in T$, we have
 $$\rho
\bigl( a \psi(a) \bigr) \asymp \|a\| \bigl( \log \|a\| \bigr)^r \asymp \|a
\psi(a)\| \bigl( \log \|a \psi(a)\| \bigr)^r .$$

All that remains is to show that $\rho(h) = O \bigl[ \|h\| \bigl( \log \|h\|
\bigr)^r \bigr]$  for every $h \in H$. Because $\rho(au) \asymp au$ for
every $au \in TU$ (see \cite[Cor.~3.18]{OhWitte-CDS}) and
 $ \| \psi(a) \| \asymp \| \psi(a)^{-1} \| \asymp (\log \|h\|)^r$, we have 
 \begin{eqnarray*}
 \rho(h)
 &=& \rho \bigl( \psi(a) \bigr) \rho(au) 
 = O \bigl[ \| \rho \bigl( \psi(a) \bigr) \| \| \rho(au) \| \bigr] \\
 &=& O \bigl[ \bigl( \log \| a \| \bigr)^r  \| au \| \bigr]
 = O \bigl[ \bigl( \log \|h\| \bigr)^r \|h\| \bigr]
 .
 \end{eqnarray*}

\begin{case}
 Assume $\omega = \alpha + \beta$.
 \end{case}
 The Weyl reflection corresponding to the root~$\alpha$ conjugates each of
\fullref{HinN-N_A(H)}{phi=0&x=0-notA-a+b} and \fullref{HinN-N_A(H)}{y&yy=0-a+b} to a
subgroup with $\omega = \beta$.

Thus, we may now assume $N_A(U) = A$. If $\sigma \neq \alpha$, then the Weyl
reflection corresponding to the root~$\alpha$ conjugates $H$ to a subgroup
with $\omega = \beta$. If $\sigma = \alpha$, then the Weyl reflection
corresponding to the root~$\beta$ does not change~$\omega$, but conjugates
$H$ to a subgroup~$H_1$ with $\sigma = \alpha+2\beta$. Then (as we already
observed) the Weyl reflection corresponding to the root~$\alpha$ conjugates
$H_1$ to a subgroup with $\omega = \beta$.

\begin{case}
 Assume $\omega = \alpha$.
 \end{case}
 Because $U$ must be normalized by the nontrivial subgroup $\psi(T)$
of~$U_\alpha$, we see that $U$ cannot be of type \fullref{HinN-N_A(H)}{H=Z-a} or
\fullref{HinN-N_A(H)}{phi=0}.

\begin{subcase}
 Assume \fullref{HinN-N_A(H)}{phi=0(dim2)-a}.
 \end{subcase}
  Because $U$ must be normalized by the nontrivial subgroup $\psi(T)$
of~$U_\alpha$, we see that $y_u = 0$ for every $u \in U$, so $\Lie U =
\Lie Z$. Thus, again using the fact that $U$ is normalized by $\psi(T)$, we
see that $\Lie U \subset \Lie U_{\alpha+2\beta} + \Lie U_{2\alpha+2\beta}$, 
and the projection of~$\Lie U$ to~$\Lie U_{\alpha+2\beta}$ is
one-dimensional. For every $z \in \Lie U$, we see that $\eta_z \neq 0$
(because $|\eta_z|^2 \neq \xx_z \yy_z$). Thus, we conclude that $\dim \Lie U
= 1$. Therefore $H$ is conjugate under~$U_\alpha$ to a subgroup of
type~\fullref{HinN-N_A(H)}{phi=0(dim2)-A} (considered in
Subsubcase~\ref{HnotTU-w=a,s=a+2b-pf} below).

\begin{subcase}
 Assume $N_A(U) = A$.
 \end{subcase}
 If $\sigma = \alpha$, then \pref{HnotTU-omega+2omega} implies that $H$ is a
Cartan-decomposition subgroup. Because $U$ is normalized by the nontrivial
subgroup $\psi(T)$ of~$U_\alpha$, we know that $\sigma \neq \beta$. 

\begin{subsubcase}
 Assume $\sigma = \alpha+\beta$.
 \end{subsubcase}
 We have
 $$ 
 h = a \psi(a) u = \begin{pmatrix}
 a_1 & a_1 \phi_{\psi(a)} & a_1 x_u & 0 & - \frac{1}{2} a_1 |x_u|^2 + i a_1
\xx_u \\
  & a_1 & 0 & 0 & 0 \\
 & & \cdots 
 \end{pmatrix}
 .
 $$
 We have  $\|h\| \asymp a_1 \log a_1 + a_1 |x_u|^2 + a_1 |\xx_u|$
 and, for $i > 1$, we have $ h_{i,j} = O \bigl( a_1 + |x_u| \bigr) $.
 The largest value of $\|\rho(h)\|$ relative to~$\|h\|$ is obtained by
taking $\log a_1 \asymp |x_u|^2$ (and $\xx_u$ small), which yields
 $\rho(h) \asymp a_1^2 \log a_1 \asymp \|h\|^2/\log\|h\|$. Because $\rho(u)
\asymp u$ for $u \in U$, we conclude that 
 $\mu(H) \approx \muH{\|h\|}{\|h\|^2/\log\|h\|}$.

\begin{subsubcase} \label{HnotTU-w=a,s=a+2b-pf}
 Assume $\sigma = \alpha+2\beta$.
 \end{subsubcase}
 We have
 $$ 
 h = a \psi(a) u = \begin{pmatrix}
 a_1 & a_1 \phi_{\psi(a)} & 0 & a_1 \eta_u & - a_1 \phi_{\psi(a)}
\cjg{\eta_u} \\
  & a_1 & 0 & 0 & - a_1 \cjg{\eta_u} \\
 & & \cdots 
 \end{pmatrix}
 .
 $$
 We have $h \asymp \bigl(1 + a_1 \|\psi(a)\| \bigr) \bigl( 1+ |\eta_u|
\bigr)$ and
 $\rho(h) \asymp a_1^2 \bigl( 1 + |\eta_u|^2 \bigr)$
 (note that 
 $\det \begin{pmatrix}
 h_{1,2} & h_{1,n+2} \\ 
 h_{2,2} & h_{2,n+2}
 \end{pmatrix}
 = 0$).
 The smallest value of $\|\rho(h)\|$ relative to~$\|h\|$ is obtained by
taking $\eta_u = O(1)$, which results in
 $\rho(h) \asymp a_1^2 \asymp \|h\|^2/ \bigl( \log \|h\| \bigr)^2$.
 Because $\rho(u) \asymp \|u\|^2$ for $u \in U$, we conclude that 
 $\mu(H) \approx \muH{\|h\|^2/ \bigl( \log \|h\| \bigr)^2}{\|h\|^2}$.

\begin{case}
 Assume $\omega = \alpha + 2\beta$.
 \end{case}
 The Weyl reflection corresponding to the root~$\beta$ conjugates
\fullref{HinN-N_A(H)}{alpha-a+2b} to a subgroup~$H'$ with $\omega = \alpha$
(of type~\fullref{HinN-N_A(H)}{phi=0(dim2)-a} with $\Lie H' = \Lie Z' \subset
\Lie U_{\alpha+2\beta} + \Lie U_{2\alpha+2\beta}$).

Thus, we may now assume $N_A(U) = A$. If $\sigma \neq \beta$, then the Weyl
reflection corresponding to the root~$\beta$ conjugates $H$ to a subgroup
with $\omega = \alpha$. Now assume $\sigma = \beta$. The Weyl reflection
corresponding to the root~$\alpha$ does not change~$\omega$, but conjugates
$H$ to a subgroup~$H_1$ with $\sigma = \alpha+2\beta$. Then (as we already
observed) the Weyl reflection corresponding to the root~$\beta$ conjugates
$H_1$ to a subgroup with $\omega = \alpha$.
 \end{proof}

\begin{Lemma}  \label{HnotTU-omega+2omega}
 Assume $G$ is a connected, almost simple, linear, real Lie group of real rank
two.
 Let $H$ be a closed, connected, nontrivial subgroup of $AN$, such that $H$ is
compatible with~$A$, and $H \neq (H \cap A) (H \cap N)$.  We use the notation
of \cite[Lem.~2.4]{OhWitte-CDS}: $T = \ker \omega$, $H = T \ltimes U$, $\psi
\colon T \to U_\omega U_{2\omega}$, and $H = \{a \psi(a) \} \ltimes U$.

 If $\Lie U \cap (\Lie U_{\omega} + \Lie U_{2\omega}) \neq 0$, then
$H$ is a Cartan-decomposition subgroup.
 \end{Lemma}

\begin{proof}
 By passing to a subgroup of~$H$, there is no harm in assuming $\Lie U \cap
(\Lie U_{\omega} + \Lie U_{2\omega})$.
 We use the notation of the proof of \cite[Prop.~3.17]{OhWitte-CDS}.
 For each $a \in T$, clearly
 $\mu_{MA} \bigl( a \psi(a) U \bigr) \supset \mu_{MA} \bigl( a \psi(a) \bigr)
A_\omega^+$, so $\mu_{MA}(H) \supset \mu_{MA} \bigl( \{ a \psi(a) \} \bigr)
A_\omega^+$. Beause $\mu_{MA}(T) = T$ is a line perpendicular to~$A_\omega$,
and $ \mu_{MA} \bigl( a \psi(a) \bigr)$ is logarithmically close to this
line, it is clear that $\mu_{MA} \bigl( a \psi(a) \bigr) A_\omega^+$ contains
all but a bounded subset of the region~$\mathord{\mathcal{C}}$. Therefore
$\mu(H)$ contains all but a bounded subset of~$A^+$, so $H$ is a
Cartan-decomposition subgroup.
 \end{proof}

\begin{Lemma}[{(cf.\ \cite[Prop.~3.16(3)]{OhWitte-CDS})}] \label{dim1-mu(H)}
 Assume that $G = \SU(2,n)$, and let $H$ be a nontrivial one-parameter
subgroup of~$AN$, such that $H$ is compatible with~$A$, but $H \neq (H \cap
A)(H \cap N)$.

Then there is a ray~$R$ in~$A^+$, a ray~$R'$ in~$A$ that is
perpendicular to~$R$, and a positive number~$k$, such that  
 $$\mu(H)
\approx  \{\, rs \mid r \in R, ~ s \in R', ~ \|s\| = (\log\|r\|)^k \,\} .$$
 \end{Lemma}

\section{Maximum dimensions of the subgroups} \label{CDSsummary}

 \begin{figure}
 \renewcommand{\arraystretch}{1.1}
 \begin{centering}
 \begin{tabular}{ccc}
 reference & Cartan projection & maximum dimension \\

\fullref{HinN-notCDS}{H=Z} & $\rho(h) \asymp h$ & $1$ \\

{\fullref{HinN-notCDS}{phi&y=0}} & $\mu(H) \approx \muH{\|h\|}{\|h\|^{3/2}}$
& $2n-3$ \\ 

{\fullref{HinN-notCDS}{phi&y=0}*} & $\rho(h) \asymp \|h\|^{3/2}$ & $1$ \\ 

\fullref{HinN-notCDS}{phi=0-neq0} & $\mu(H) \approx
\muH{\|h\|}{\|h\|^{3/2}}$ & $2n-3$ \\

\fullref{HinN-notCDS}{phi=0-neq0}* & $\rho(h)
\asymp \|h\|^{3/2}$ & $1$ \\

\fullref{HinN-notCDS}{phi=0-=0} & $\rho(h) \asymp h$ & $2n-3$ \\

{\fullref{HinN-notCDS}{y&yy=0}} & $\rho(h) \asymp h$ & $2n-1$ \\

\fullref{HinN-notCDS}{phi=yy} & $\mu(H) \approx
\muH{\|h\|}{\|h\|^{4/3}}$ & $2n-3$ \\

\fullref{HinN-notCDS}{phi=yy}* & $\rho(h)
\asymp \|h\|^{4/3}$ & $1$ \\

\fullref{HinN-notCDS}{phi=0(dim2)} & $\rho(h) \asymp \|h\|^2$ & 
 $\begin{cases}
 2n-1 & \text{$n$ even} \\
 2n-3 & \text{$n$ odd} 
 \end{cases}
 $
 \\

{\fullref{HinN-notCDS}{phi=0(dim1&2)}} & $\mu(H) \approx
\muH{\|h\|^{3/2}}{\|h\|^2}$ & 
 $\begin{cases}
 n+1 &  n \ge 4 \\
 \hfil 3 & n=3
 \end{cases}
 $
 \\

\fullref{HinN-notCDS}{onlyphiyyxx} & $\rho(h) \asymp \|h\|^{3/2}$ & 
 $\begin{cases}
 3 &  n \ge 4 \\
 2 & n=3
 \end{cases}
 $
 \\

\fullref{HinN-notCDS}{dim3} & $\mu(H) \approx \muH{\|h\|}{\|h\|^{3/2}}$ &
$2$ \\

\fullref{HinN-notCDS}{alpha} & $\rho(h) \asymp \|h\|^2$ & $2$ \\

\fullref{HinN-notCDS}{u+z} & $\mu(H) \approx
\muH{\|h\|^{5/4}}{\|h\|^2}$ & 2 \\

 \end{tabular}
 \caption{The subgroups of $N$ that are not Cartan-decomposition subgroups.
 \label{CprojN-summaryfigure}}
 \end{centering}
 
 \end{figure}

 \begin{figure}
 \renewcommand{\arraystretch}{1.1}
 \begin{centering}
 \begin{tabular}{ccc}
 reference & Cartan projection & maximum dimension \\

 \fullref{TU-proj}{H=Z-A} & $\mu(H) \approx \muH{\|h\|}{\|h\|^s}$ & 2 \\

 {\fullref{TU-proj}{phi&y=0}} & $\mu(H) \approx \muH{\|h\|}{\|h\|^{3/2}}$ &
$2n-2$ \\

 {\fullref{TU-proj}{phi&y=0}*} & $\rho(h) \asymp \|h\|^{3/2}$ & 2 \\

 {\fullref{TU-proj}{phi=0&x=0-A}} & $\mu(H) \approx \muH{\|h\|}{\|h\|^s}$ &
$2n-2$ \\

 {\fullref{TU-proj}{phi=0&x=0-notA-a+b}} & $\rho(h) \asymp h$ & $2n-2$\\

 {\fullref{TU-proj}{phi=0&x=0-notA-2a+b}} & 
$\mu(H) \approx \muH{\|h\|}{\|h\|^{3/2}}$ & $2n-2$ \\

 {\fullref{TU-proj}{phi=0&x=0-notA-2a+b}*} & $\rho(h) \asymp \|h\|^{3/2}$ & 2
\\

 {\fullref{TU-proj}{y&yy=0-A}} &  $\mu(H) \approx \muH{\|h\|}{\|h\|^s}$ &
$2n-2$ \\

 {\fullref{TU-proj}{y&yy=0-b}} & $\rho(h) \asymp h$ & $2n$ \\

 \fullref{TU-proj}{phi=yy} & 
$\mu(H) \approx \muH{\|h\|}{\|h\|^{4/3}}$ & $2n-2$ \\

 \fullref{TU-proj}{phi=yy}* & $\rho(h) \asymp \|h\|^{4/3}$ & 2 \\

 \fullref{TU-proj}{phi=0(dim2)-A} & $\mu(H) \approx
\muH{\|h\|^s}{\|h\|^2}$ & $3$ \\

 \fullref{TU-proj}{phi=0(dim2)-a} & $\rho(h) \asymp \|h\|^2$ & 
 $\begin{cases}
 2n & \text{$n$ even} \\
 2n-2 & \text{$n$ odd} 
 \end{cases}
 $
 \\

 {\fullref{TU-proj}{phi=0(dim1&2)-x}} & 
$\mu(H) \approx \muH{\|h\|^{3/2}}{\|h\|^2}$ & $4$ \\

 {\fullref{TU-proj}{phi=0(dim1&2)-y}} & $\mu(H) \approx
\muH{\|h\|^{3/2}}{\|h\|^2}$ & $4$ \\

 \fullref{TU-proj}{onlyphiyyxx} & $\rho(h) \asymp \|h\|^{3/2}$  & 
 $\begin{cases}
 4 & n \ge 4 \\
 3 & n = 3 
 \end{cases}
 $
 \\

 \fullref{TU-proj}{dim3} & $\mu(H) \approx
\muH{\|h\|}{\|h\|^{3/2}}$ & $3$ \\

 \fullref{TU-proj}{alpha-A} & $\mu(H) \approx \muH{\|h\|^s}{\|h\|^2}$ & $3$ \\

 \fullref{TU-proj}{alpha-a+2b} & $\rho(h)
\asymp \|h\|^2$ &$3$ \\

 \end{tabular}
 \caption{The subgroups of $AN$ that are not Cartan-decomposition subgroups,
and are a nontrivial semidirect product $T \ltimes U$.
 \label{Cprojsemi-summaryfigure}}
 \end{centering}
 
 \end{figure}

 \begin{figure}
 \renewcommand{\arraystretch}{1.1}
 \begin{centering}
 \begin{tabular}{ccc}
 reference & Cartan projection & maximum dimension \\

 \fullref{HnotTUproj}{a,a+b} & 
 $\mu(H) \approx \muH{\|h\|}{\|h\|^2/(\log \|h\|)}$ & $2n-2$ \\

 \fullref{HnotTUproj}{a,a+2b} & 
 $\mu(H) \approx \muH{\|h\|^2/(\log \|h\|)^2}{\|h\|^2}$ & $2$ \\

 \fullref{HnotTUproj}{b,a+2b} ($r=1$) & 
 $\mu(H) \approx \muH{\|h\| \bigl( \log \|h\| \bigr)^{1/2}}{\|h\|^2}$ & $3$ \\

 \fullref{HnotTUproj}{b,a+2b} ($r=2$) & 
 $\mu(H) \approx \muH{\|h\| \bigl( \log \|h\| \bigr)}{\|h\|^2}$ & $3$ \\

 \fullref{HnotTUproj}{b,a+b} ($r=1$) & 
 $\mu(H) \approx \muH{\|h\|}{\|h\|  \bigl( \log \|h\| \bigr)}$ & $2n-2$
\\

 \fullref{HnotTUproj}{b,a+b} ($r=2$) & 
 $\mu(H) \approx \muH{\|h\|}{\|h\|  \bigl( \log \|h\| \bigr)^{2}}$ & $2n-3$
\\

 \ref{dim1-mu(H)} & 
 $\rho(h) \asymp \|h\|^s (\log \|h\|)^{\pm k}$ & $1$
\\

 \end{tabular}
 \caption{The subgroups of $AN$ that
are not Cartan-decomposition subgroups \textbf{and} are not a semidirect product
of a torus and a unipotent subgroup.
 \label{Cprojnotsemi-summaryfigure}}
 \end{centering}
 
 \end{figure}

For convenience of reference, Tables~\ref{CprojN-summaryfigure},
\ref{Cprojsemi-summaryfigure} and~\ref{Cprojnotsemi-summaryfigure} list the
(approximate) Cartan projection of each subgroup of~$AN$ that is not a
Cartan-decomposition subgroup. The maximum possible dimension for a subgroup of
each type is also listed. (These dimensions are used in applications to the
existence of tessellations.)

\begin{Remark}
 Here are brief justifications of the dimensions listed in
Tables~\ref{CprojN-summaryfigure}, \ref{Cprojsemi-summaryfigure}
and~\ref{Cprojnotsemi-summaryfigure}.

\fullref{HinN-notCDS}{H=Z} By assumption, we have $\dim H = 1$.

{\fullref{HinN-notCDS}{phi&y=0}} Let $p \colon \Lie H \to \Lie
U_{\alpha+\beta}$ be the natural projection. Then $\ker p = \Lie Z \subset
\Lie U_{2\alpha+2\beta}$, so
 $$ \dim \Lie H \le (\dim \Lie U_{\alpha+\beta}) + (\dim \Lie
U_{2\alpha+2\beta})
 = 2(n-2) + 1 = 2n-3 .$$

\fullref{HinN-notCDS}{phi=0} We may assume $\lambda = 0$. Then $\Lie H
\subset \Lie U_\beta + \Lie U_{2\beta}$. So  
 $$ \dim \Lie H \le (\dim \Lie U_\beta) + (\dim \Lie U_{2\beta})
 = 2(n-2) + 1 = 2n-3 .$$
 It is easy to construct an algebra of this dimension, with or without an
element~$u$ as described in~\pref{HinN-notCDS-phi=0-neq0}.

{\fullref{HinN-notCDS}{y&yy=0}} Let $V$ be the projection of~$\Lie H$ to
$\Lie U_\alpha + \Lie U_{\alpha+\beta} + \Lie U_{\alpha+2\beta}$. Because $\phi
\cjg{\eta}$ is a form of signature $(2,2)$ on $\Lie U_\alpha + \Lie
U_{\alpha+2\beta}$, we know that $\dim \bigl( V \cap (\Lie U_\alpha + \Lie
U_{\alpha+2\beta}) \bigr) \le 2$. Thus we have
  \begin{eqnarray*}
 \dim \Lie H
 &\le& \dim V + \dim \Lie U_{2\alpha+2\beta}
 \le \bigl( \dim \Lie U_{\alpha+\beta} + 2 \bigr) + \dim \Lie
U_{2\alpha+2\beta} \\
 &=& \bigl( 2(n-2) + 2 \bigr) + 1 = 2n-1 .
 \end{eqnarray*}

\fullref{HinN-notCDS}{phi=yy} Consider $p \colon \Lie H \to \Lie
U_{\alpha}$. Because $\Lie Z = 0$, we have $\dim \ker p \le \dim \Lie
U_{\alpha+\beta} = 2n-4$. Because $p(\Lie H) \subset \real \phi_0$, we have
$\dim p(\Lie H) \le 1$. Thus,
 $\dim \Lie H \le 2n-3$.

\fullref{HinN-notCDS}{phi=0(dim2)} See Lemma~\ref{max-phi=0(dim2)} below.

{\fullref{HinN-notCDS}{phi=0(dim1&2)}} See Lemma~\ref{max-phi=0(dim1&2)}
below.

{\fullref{HinN-notCDS}{onlyphiyyxx}} See Lemma~\ref{max-onlyphiyyxx}
below.

\fullref{HinN-notCDS}{dim3}, 
\fullref{HinN-notCDS}{alpha},
\fullref{HinN-notCDS}{u+z} are obvious from the statements.

 \fullref{TU-proj}{H=Z-A} Because $\dim \Lie U = 1$, we have $\dim \Lie H =
\dim \Lie T + \dim \Lie U = 2$.

 {\fullref{TU-proj}{phi&y=0}} The kernel of the projection from~$\Lie U$
to~$\Lie U_{\alpha+\beta}$ is~$\Lie Z$, so 
 $\dim \Lie H = 1 + \dim U \le 1 + (1 + \dim \Lie U_{\alpha+\beta}) = 2n-2$.

 {\fullref{TU-proj}{phi=0&x=0}} $\dim \Lie H = 1 + \dim \Lie U
 \le 1 + \bigl( \dim \Lie U_{\beta} + \dim \Lie Z \bigr) = 2n-2$. 

 {\fullref{TU-proj}{y&yy=0-A}} $\dim \Lie H
 \le \dim \Lie T + \dim \Lie U_{\alpha+\beta} + \dim \Lie Z
 \le 1 + (2n-4) + 1  = 2n-2$.

 {\fullref{TU-proj}{y&yy=0-b}} Add~$1$ (the dimension of~$T$) to the bound
in~{\fullref{HinN-notCDS}{y&yy=0}}.

 \fullref{TU-proj}{phi=yy} Add~$1$ (the dimension of~$T$) to the bound
in~{\fullref{HinN-notCDS}{phi=yy}}.

 \fullref{TU-proj}{phi=0(dim2)-A} $\dim \Lie H \le \dim \Lie T + \dim \Lie
U_{\alpha+2\beta} = 1+2 = 3$.

 \fullref{TU-proj}{phi=0(dim2)-a}  Add~$1$ (the dimension of~$T$) to the
bound in~{\fullref{HinN-notCDS}{phi=0(dim2)}}.

 {\fullref{TU-proj}{phi=0(dim1&2)-x}} Because $\yy_u \neq 0$ for every
nonzero $u \in \Lie U \cap ( \Lie U_{\alpha+\beta} + \Lie
U_{2\beta})$, we have $\dim \bigl( \Lie U \cap ( \Lie U_{\alpha+\beta} + \Lie
U_{2\beta}) \bigr) \le 1$. Therefore
 $\dim \Lie H  \le \dim \Lie T + 1 + \dim \Lie U_{\alpha+2\beta} = 4$.

 {\fullref{TU-proj}{phi=0(dim1&2)-y}} This is conjugate
to~{\fullref{TU-proj}{phi=0(dim1&2)-x}}, via the Weyl reflection
corresponding to the root~$\alpha$.

 \fullref{TU-proj}{onlyphiyyxx}  Add~$1$ (the dimension of~$T$) to the bound
in~{\fullref{HinN-notCDS}{onlyphiyyxx}}. (To achieve this bound for $n \ge
4$, choose $u, \tilde u \in \Lie U \cap ( \Lie U_{\alpha} + \Lie
U_\beta)$ in the proof of Lemma~\ref{max-onlyphiyyxx}.

 \fullref{TU-proj}{dim3} and
 \fullref{TU-proj}{alpha} are obvious from the statements.

 \fullref{HnotTUproj}{a,a+b} $\dim \Lie H \le 1 + \dim( \Lie U_{\alpha+\beta}
+ \Lie U_{2\alpha+2\beta}) = 2n-2$.

 \fullref{HnotTUproj}{a,a+2b} Because $\psi(T)$ normalizes (hence
centralizes)~$U$, the subgroup~$U$ cannot be all of~$U_{\alpha+2\beta}$, so
$\dim U \le 1$. Therefore
 $\dim H = 1 + \dim U \le 2$.

 \fullref{HnotTUproj}{b,a+2b} 
 $\dim \Lie H \le 1 + \dim \Lie U_{\alpha+2\beta} = 3$.

 \fullref{HnotTUproj}{b,a+b} Because $\psi(T)$ normalizes (hence
centralizes)~$U$, the projection of~$\Lie U$ to~$\Lie U_{\alpha+\beta}$ cannot
be all of~$\Lie U_{\alpha+\beta}$ if $\psi(T) \not\subset U_{2\beta}$, that
is, if $r = 2$. Therefore $\dim U \le \dim(\Lie U_{\alpha+\beta} + \Lie
U_{2\alpha+2\beta}) - (r-1) =2n-2 - r$. Therefore
 $\dim \Lie H = 1 + \dim U \le 2n-1-r$.

 \end{Remark}

\begin{Lemma} \label{max-phi=0(dim2)}
 The maximum dimension of a subalgebra of
type~\fullref{HinN-notCDS}{phi=0(dim2)} is as stated in
Table~\ref{CprojN-summaryfigure}.
 \end{Lemma}

\begin{proof}
 We begin by showing that $\dim \Lie H \le 2n-1$ (cf.
\cite[Lem.~5.8]{OhWitte-CK}). Let $V$ be the projection of $\Lie H$ to $\Lie
U_\beta + \Lie U_{\alpha+\beta}$. Because $\dim \Lie Z \le 3$, we just need to
show that $\dim V \le 2n-4$. Because $V$ does not intersect $\Lie U_\beta$ (or
$\Lie U_{\alpha+\beta}$, either, for that matter), and $\Lie U_\beta$ has
codimension $2n-4$ in $\Lie U_\beta + \Lie U_{\alpha+\beta}$, this is immediate.

When $n$ is even, there is a subgroup of dimension $2n-1$. (For example,  the
$N$ subgroup of $\SP(1,n/2)$. More general examples are constructed in
\cite[\S4]{OhWitte-eg}.)

Let us show that if $n$ is odd, then $\dim H \le 2n-3$. (Our proof is
topological; we do not know an algebraic proof.)  Suppose that $\dim H
\ge 2n-2$ (this will lead to a contradiction). Because $\dim \Lie Z \le 3$,
we have $\dim \Lie H / \Lie Z \ge 2n-5$. Thus, there is a
$(2n-5)$-dimensional real subspace~$X$ of~$\complex^{n-2}$ and a real linear
transformation $T \colon X \to \complex^{n-2}$, such that $x$ and~$Tx$ are
linearly independent over~$\complex$, for every nonzero $x \in X$ (cf.\
\cite[Cor.~5.9]{OhWitte-CK}).
 Thus, if we define $U \colon X \to \complex^{n-2}$ by $Ux = ix$; then $x$,
$Tx$, and $Ux$ are linearly independent over~$\real$, for every nonzero $x
\in X$. Thus (writing $n = 2k+3$): there is a $(4k+1)$-dimensional real
subspace~$X$ of~$\real^{4k+2}$ and real linear transformations $T,U \colon X
\to \real^{4k+2}$, such that $x$, $Tx$, and $Ux$ are linearly independent
over~$\real$, for every nonzero $x \in X$. There is no harm in assuming $X =
\real^{4k+1}$ (under its natural embedding in~$\real^{4k+2}$).

Let $E = (S^{4k} \times \real^{4k+2})/ {\sim}$, where $(x,v) \sim (-x,-v)$,
and define a continuous map $\zeta \colon E \to \real P^{4k}$ by $\zeta(x,v)
= [x]$, so $(E, \zeta)$ is a vector bundle over $\real P^{4k}$. Then $(E,
\zeta) \cong \tau \oplus \epsilon^1 \oplus \gamma^1_{4k}$, where $\tau$ is the
tangent bundle of $\real P^{4k}$, $\epsilon^1$ is a trivial line bundle, and
$\gamma^1_{4k}$ is the canonical bundle of $\real P^{4k}$. (To see this, note
that the subbundle 
 $$ \{\, (x,v) \in S^{4k} \times \real^{4k+1} \mid v \perp x \,\}/ {\sim}$$
 is the total space of~$\tau$ \cite[pf.\ of Lem.~4.4,
pp.~43--44]{MilnorStasheff}, the subbundle 
 $$ \{\, (x,v) \in S^{4k} \times \real^{4k+1} \mid v \in \real x \,\}/
{\sim}$$
 has the  obvious section $x \mapsto (x,x)$,  and the subbundle $\bigl(
S^{4k} \times (0 \times \real) \bigr)/ {\sim}$ is isomorphic to
$\gamma^1_{4k}$ via the bundle map $\bigl( x, (0,t) \bigr) \mapsto (x,tx)$.)
Therefore, letting $a$ be a generator of the cohomology ring $H^*(\real
P^{4k}; \integer_2)$, we see that the total Stiefel-Whitney class of
$(E,\zeta)$ is
 $ w = (1+a)^{4k+1}(1)(1+a) = (1+a)^{4k+2}$
 \cite[Eg.~2, p.~43, and Thm.~4.5, p.~45]{MilnorStasheff},
 so 
 $$ w_{(4k+2)-3+1} = w_{4k} =
 \left( \begin{matrix}4k+2 \\ 4k \end{matrix} \right)
 a^{4k} =
(2k+1)(4k+1) a^{4k} \neq 0 $$
 (because $(2k+1)(4k+1)$ is odd). Therefore, there do not exist three
pointwise linearly independent sections of $(E, \zeta)$ \cite[Prop.~4,
p.~39]{MilnorStasheff}.

Any linear transformation $Q \colon \real^{4k+1} \to \real^{4k+2}$ induces a
continuous function $\hat Q \colon S^{4k} \to \real^{4k+2}$, such that $\hat
Q(-x) = - \hat Q(x)$ for all $x \in S^{4k}$; that is, a section of $(E,
\zeta)$. Thus, $\Id$, $T$, and~$U$ each define a section of $(E, \zeta)$.
Furthermore, these three sections are pointwise linearly independent, because
$x$, $Tx$, and $Ux$ are linearly independent over~$\real$, for every $x \in
S^{4k}$. This contradicts the conclusion of the preceding paragraph.
 \end{proof}

\begin{Lemma} \label{max-phi=0(dim1&2)}
 The maximum dimension of a subalgebra of
type~\fullref{HinN-notCDS}{phi=0(dim1&2)} is as stated in
Table~\ref{CprojN-summaryfigure}.
 \end{Lemma}

\begin{proof}
 Replacing $H$ by a conjugate under
$\langle U_\alpha, U_{-\alpha} \rangle$, we may assume $x_v = 0$. Therefore
$\xx_z = 0$ for every $z \in \Lie Z$. (Thus, in particular, we have $\dim
\Lie Z \le 2$.)

For the projection $p \colon \Lie H \to \Lie U_{\alpha+\beta}$, we have
$\ker p = \real v + \Lie Z$. (There cannot exist a linearly independent
$v'$; otherwise, replacing $v'$ by some linear combination with~$v$, we
could assume $\xx_{v'} = 0$, which is impossible.) Thus, $\dim \ker p \le 3$.

Because $\xx_z = 0$ for every $z \in \Lie Z$, $p(\Lie H)$ must be a totally
isotropic subspace for the symplectic form $\Im(x \tilde x^{\dagger})$, so
$\dim p(\Lie H) \le n-2$. Therefore $\dim \Lie H \le (n-2) + 3 = n+1$.

For $n \ge 4$, here is an example that achieves this bound:
 $$ \Lie H =
 \bigset{ \begin{pmatrix}
 0 & 0 &  x_1 & x_2 & \cdots & x_{n-2} & \eta & i \xx \\
 0 & 0 & i\xx & x_1 & \cdots & x_{n-3} & i x_{n-2} & - \cjg{\eta} \\
 & & & \cdots 
 \end{pmatrix}
 }{
 \begin{matrix}
 \xx, x_1,\ldots,x_{n-2} \in \real, \\
  \eta \in \complex
 \end{matrix}
 }
 .$$
 For $v \in \Lie H$, we claim that $\dim _{\complex} \langle x_v, y_v \rangle
= 1$ only if either $x_v = 0$ or $y_v = 0$. (In either case, it is clear from
the definition of~$\Lie H$ that either $\xx_u |y_u|^2 \neq 0$ or $\yy_u
|x_u|^2 \neq 0$, respectively.) Suppose $\dim _{\complex} \langle x_v, y_v
\rangle = 1$, with $x_v \neq 0$ and $y_v \neq 0$. There is some nonzero
$\lambda \in \complex$, such that $y_v = \lambda x_v$. We must have $x_1 \neq
0$. (Otherwise, let $i \in \{1,2,\ldots,n-2\}$ be minimal with $x_i \neq 0$.
Then $x_{i-1} = y_i = \lambda x_i \neq 0$, contradicting the minimality
of~$i$.)
 Because $y_1 = i \xx$ is pure imaginary, but $x_1$ is real, we see that
$\lambda$ is pure imaginary. On the other hand, $y_2 = x_1$ is real (and
nonzero), and $x_2$ is also real, so $\lambda$ is real. Because $\lambda
\neq 0$, this is a contradiction.

Now let $n = 3$, and suppose $\dim \Lie H = 4$. (This will lead to a
contradiction.) Because equality is attained in the proof above, we must have
$\dim p(\Lie H) = n-2 = 1$ and $\dim \Lie Z = 2$. In particular, there exists
$w \in \Lie H$ with $x_w \neq 0$. For $t \in \real$, let $w_t = w + tv$. Then
 $$ \xx_{w_t} |y_{w_t}|^2 + \yy_{w_t}|x_{w_t}|^2 + 2 \Im(
x_{w_t}y_{w_t}^{\dagger} \eta_{w_t})
 = t^3 \xx_v |y_v|^2 + O(t^2)
 \to
 \begin{cases}
 +\xx_v \infty & \text{as $t \to \infty$} \\
 -\xx_v \infty & \text{as $t \to -\infty$}
 .
 \end{cases}
 $$
 Thus, this expression changes sign, so it must vanish for some~$t$. This is a
contradiction, because
  $\dim _{\complex} \langle x_{w_t}, y_{w_t} \rangle = 1$ for every~$t$.
 \end{proof}

\begin{Lemma} \label{max-onlyphiyyxx}
 The maximum dimension of a subalgebra of
type~\fullref{HinN-notCDS}{onlyphiyyxx} is as stated in
Table~\ref{CprojN-summaryfigure}.
 \end{Lemma}

\begin{proof}
 For $n \ge 4$, here is the construction of  3-dimensional subalgebras of~$\Lie
N$ of this type. Let $\phi = 1$ and $\tilde \phi = i$. Choose $y,\tilde y, x,
\tilde x \in \complex^{n-2}$, $\eta, \tilde \eta \in \complex$, and $\xx,\tilde
\xx \in \real$, such that 
 \begin{equation} \label{eg3D(1,i)}
 |y|^2 = |\tilde y|^2 = 3i y \tilde y^{\dagger} \neq 0 .
 \end{equation}
 Now, choose $\yy,\tilde\yy \in \real$, such that
 \begin{equation} \label{xx[u,v,u]}
 \Im( \tilde y x^{\dagger} - i y x^{\dagger}
 + \tilde y x^{\dagger} - y \tilde x^{\dagger} + i \tilde \yy) = 0 
 \end{equation}
 and
 \begin{equation} \label{xx[u,v,v]}
 \Im( \tilde y \tilde x^{\dagger} - i y \tilde x^{\dagger}
 + i \tilde y x^{\dagger} - i y \tilde x^{\dagger} + i \yy) = 0
 .
 \end{equation}
 Define $u, \tilde u$ as in \eqref{umatrix}, and let $v = [u, \tilde u]$. Then
$\yy_v \neq 0$ and $x_v \neq 0$, but, from
\eqref{eg3D(1,i)}, \eqref{xx[u,v,u]} and~\eqref{xx[u,v,v]}, we have $[v, u] =
[v, \tilde u] = 0$. Thus, we may let $\Lie H$ be the subalgebra generated by
$u$ and~$\tilde u$. (So $\{u, \tilde u, v\}$ is a basis of~$\Lie H$
over~$\real$.)

 Note that, because $|y \tilde y^{\dagger}| = |y|^2/3 \neq |y|^2$, we know that
$y$ and~$\tilde y$ must be linearly independent over~$\complex$. Thus, these 
3-dimensional examples do not exist when $n = 3$.
 \end{proof}

\references
\lastpage
\end{document}